\documentclass[12pt,a4paper]{amsart}
\usepackage[english]{babel}
\usepackage[utf8]{inputenc}
\usepackage{fancyhdr}
\usepackage{xcolor}

\usepackage{amssymb}
\usepackage{amsmath} 
\usepackage{amsthm} 
\usepackage{enumitem}
\usepackage{bbm}
\usepackage{bbold} 
\usepackage{graphicx}

\newtheorem{thm}{Theorem}
\newtheorem{conj}{Conjecture}
\newtheorem{prop}{Proposition}
\newtheorem{lemma}{Lemma}

\theoremstyle{proposition}
\newtheorem*{prop*}{Proposition}
\theoremstyle{proposition}
\newtheorem*{rem*}{Remark}

\numberwithin{thm}{section}
\numberwithin{conj}{section}
\numberwithin{lemma}{section}
\numberwithin{equation}{section}
\numberwithin{prop}{section}
\numberwithin{re}{section}

\title{A weighted one-level density of families of $L$-functions}
\author{Alessandro Fazzari}
\address{Universit\`a  di Genova,  Dipartimento di Matematica. 
 Via Dodecaneso 35, 16146 Genova, Italy}
\address{Charles University, Faculty of Mathematics and Physics, Department of Algebra, Sokolovska 83, 18600 Praha 8, Czech Republic}
\email{fazzari@dima.unige.it}

\DeclareMathOperator{\supp}{ supp}
\DeclareMathOperator{\sinc}{ sinc}

\newmuskip\pFqmuskip
\newcommand*\pFq[6][8]{%
  \begingroup % only local assignments
  \pFqmuskip=#1mu\relax
  \mathchardef\normalcomma=\mathcode`,
  \mathcode`\,=\string"8000
  \begingroup\lccode`\~=`\,
  \lowercase{\endgroup\let~}\pFqcomma
  {}_{#2}F_{#3}{\left[\genfrac..{0pt}{}{#4}{#5};#6\right]}%
  \endgroup
}
\newcommand{\pFqcomma}{{\normalcomma}\mskip\pFqmuskip}

\begin{document}
\maketitle

\begin{abstract} 
This paper is devoted to a weighted version of the one-level density of the non-trivial zeros of $L$-functions, tilted by a power of the $L$-function evaluated at the central point. Assuming the Riemann Hypothesis and the ratio conjecture, for some specific families of $L$-functions we prove that the same structure suggested by the density conjecture holds also in this weighted investigation, if the exponent of the weight is small enough. Moreover we speculate about the general case, conjecturing explicit formulae for the weighted kernels.
\end{abstract}
%%%%%%%%%%%%%%%%%%%%%%%%%%%%%%%%%%%%%%%%%%%%%%%%%%%%%%%%%%%%%%%%%%%%%%%%%%%%

\section{A weighted version of the one-level density}\label{C4S1*}

Let us assume the Riemann Hypothesis for all the $L$-functions that arise. The classical one-level density considers a smooth localization at the central point of the counting function of the non-trivial zeros of an $L$-function, averaged over a \lq\lq natural'' family of $L$-functions in the Selberg class\footnote[1]{We refer e.g. to \cite{KP1} for the definitions and the basic properties of the Selberg class.}. More specifically, given an even and real-valued function $f$ in the Schwartz space\footnote[2]{In practice we will see that this condition can be weakened and a decay like $f(x)\ll1/(1+x^2)$ at infinity will suffice.} and an $L$-function $L(s)$ in a family $\mathcal F$, we consider the quantity
\begin{equation}\label{intro4.0} \sum_{\gamma_L} f(c(L)\gamma_L)\end{equation}
where $\gamma_L$ denotes the imaginary part of a generic non-trivial zero of $L$ and $c(L)$ the log-conductor of $L(s)$ at the central point. We recall that $1/c(L)$ is the mean spacing of the non-trivial zeros of $L(s)$ around $s=\frac{1}{2}$. The one-level density for the family $\mathcal F$ is the average of the above quantity over the family, i.e.
\begin{equation}\label{intro4.1} \lim_{X\to\infty}\frac{1}{\sum_{L\in\mathcal F_X}1}\sum_{L\in\mathcal F_X}\sum_{\gamma_L}f(c(L)\gamma_L), \end{equation}
with $$\mathcal F_X:=\{L\in\mathcal F: c(L)\leq \log X\}.$$
In the literature this is also referred to as the \lq\lq low-lying zeros'' density, as the sum (\ref{intro4.0}) gives information on the distribution of the zeros of $L$ which are close to the central point. Indeed if a zero is substantially more than $1/c(L)$ away from the central point, then it does not contribute significantly to the sum (see e.g. \cite{ILS} for a complete overview). 

Katz and Sarnak \cite{KatzSarnak2} studied a wide variety of families and attached to each of these families of $L$-functions a symmetry type (i.e. unitary, symplectic, or orthogonal, hereafter identified by a group $G$), which should govern the one-level density of the considered family. Namely, the density conjecture predicts that
\begin{equation}\label{intro4.2} \lim_{X\to\infty}\frac{1}{\sum_{L\in\mathcal F_X}1}\sum_{L\in\mathcal F_X}\sum_{\gamma_L}f(c(L)\gamma_L)=\int_{-\infty}^{+\infty}f(x)W_{\mathcal F}(x)dx \end{equation}
where $W_{\mathcal F}$ equals the one-level density function $W_{G}$ for the (scaled) limit of $G\in\{U(N),$ $U\!Sp(2N), O(N), SO(2N), SO(2N+1)\}$\footnote[3]{Note that the scaled limit of $SO(2N)$ (resp. $SO(2N+1)$) is commonly denoted by $SO^+$ (resp. $SO^-$) in the literature and also in the rest of this paper.}, i.e. the kernel appearing in the analogous average in the corresponding random matrix theory setting. In particular, the kernel $W_{\mathcal F}$ is predicted to depend on $G$ only. We recall that the function $W_G$ is known for all the classical compact groups, being
\begin{equation}\begin{split}\notag
& W_U(x)=1,
\\& W_{U\!Sp}(x)=1-\frac{\sin(2\pi x)}{2\pi x},
\\& W_{O}(x)=1+\frac{1}{2}\delta_0(x),
\\& W_{SO^+}(x)=1+\frac{\sin(2\pi x)}{2\pi x},
\\& W_{SO^-}(x)=\delta_0(x)+1-\frac{\sin(2\pi x)}{2\pi x},
\end{split}\end{equation}
with $\delta_0$ the Dirac $\delta$-function centered at 0. Examples of one-level density theorems which prove (\ref{intro4.2}) in specific cases can be found e.g. in \cite{ILS, HR, HR1, OS, MSConrey, CSapplications}. \\

In this paper, we investigate a weighted analogue of the one-level density. In particular we consider a tilted average over the family $\mathcal F$ of the quantity (\ref{intro4.0}), multiplied by a power of $L$ evaluated at the central point. 
%\textcolor{red}{This approach links the results on moments to those on the one-level density, making the connection between non-trivial zeros and the size of $L(\tfrac{1}{2})$ explicit. In the same spirit as in \cite{1., 2.}, the weighted average we consider amplifies the contribution coming from the $L$-functions that are large at the central point, near which zeros are expected to be rarer.}
The philosophy of this tilted average is similar to that of \cite{1., 2.}; the weight has the effect of giving more relevance to the $L$-functions which are large at the central point, near which zeros are expected to be rarer.

More specifically, given $k\in\mathbb N$, we are interested in
\begin{equation}\label{intro4.3}
\mathcal D^{\mathcal F}_k(f)=\mathcal D^{\mathcal F}_k(f,X):=\frac{1}{\sum\limits_{L\in\mathcal F_X}V(L(\tfrac{1}{2}))^k}\sum_{L\in\mathcal F_X}\sum_{\gamma_L}f(c(L)\gamma_L)V(L(\tfrac{1}{2}))^k 
\end{equation}
in the limit $X\to\infty$, where $V$ depends on the symmetry type of the family; in particular $V(z)=|z|^2$ in the unitary case and $V(z)=z$ for the symplectic and orthogonal cases. 
The quantity $\mathcal D^{\mathcal F}_k(f)$ links the moments to the one-level density, making the connection between non-trivial zeros and the size of $L(\tfrac{1}{2})$ explicit. Indeed, $\mathcal D^{\mathcal F}_k(f)$ can be seen as a special case of 
\begin{equation}\label{Sandro.Ltaliche1} \frac{1}{\sum\limits_{L\in\mathcal F_X}V(L(\tfrac{1}{2}))^k}\sum_{L\in\mathcal F_X}g(L)V(L(\tfrac{1}{2}))^k \end{equation}
with $g(L)$ a function over the $L$-functions of a given family $\mathcal F$. 
In the unitary case, for example, we know from Soundararajan's work \cite{Sound} that the dominant contribution to the $2k$-th moment comes from those $L$-functions such that the size of $|L(\tfrac{1}{2})|$ is about $(\log X)^{k+o(1)}$, which form a thin subset of size about $\#\mathcal F_X/(\log T)^{k^2+o(1)}$. Thus, if the function $g$ has size 1, then only these $L$-functions contribute to the main term of the sum in (\ref{Sandro.Ltaliche1}). With the choice we made in (\ref{intro4.3}), we have $g(L)=\sum_{\gamma_L}f(c(L)\gamma_L)$, which is not bounded but only $\ll c(L)$, by the Riemann-Von Mangoldt formula. However, the standard $n$-th level density \cite{RudSar2} implies that $g(L)\ll c(L)^\varepsilon$ for all but $\#\mathcal F_X/(\log X)^A$ $L$-functions in the family, for every $A>0$. Therefore, also in (\ref{intro4.3}), we have that only the $L$-functions such that $|L(\tfrac{1}{2})|\asymp (\log X)^{k\pm\varepsilon}$ contribute significantly to the main term of the sum. 
%For this reason, $\mathcal D_k^{\mathcal F}$ is a weighted one-level density, focused on the $L$-functions in the family which are responsible to the $k$-th moment (or $2k$-th moment in the unitary case).
For this reason, for unitary families, $\mathcal D_k^{\mathcal F}$ can be interpreted as a (weighted) one-level density for the thin subset $\{ L \in \mathcal F : (\log X)^{k-\varepsilon}\ll |L(\frac{1}{2})| \ll (\log X)^{k+\epsilon} \}$. Similarly, in the symplectic and orthogonal cases, $\mathcal D_k^{\mathcal F}$ is a weighted one-level density, focused on the $L$-functions in the family which are responsible to the $k$-th moment.

From the computations we perform throughout this paper in some specific cases, we speculate that the structure suggested by the density conjecture  (\ref{intro4.2}) holds also in the weighted case. Namely, we expect that
\begin{equation}\label{intro4.4} 
\mathcal D^{\mathcal F}_k(f)=\int_{-\infty}^{+\infty}f(x)W_G^k(x)dx+O\Big(\frac{1}{\log X}\Big)
\end{equation}
where the weighted one-level density function $W_G^k$ only depends on $k$ and on the symmetry type of the family $\mathcal F$. Note that the superscript $^k$ is an index, indicating that we are weighting with the $k$-th power of $V(L(\tfrac{1}{2}))$; in particular $W_G^k$ is not the $k$-th power of $W_G$.\\

This kind of weighting naturally appears also in other contexts, such as Kowalski, Saha and Tsimerman's paper \cite{KowSahaTsim}. Given a Siegel modular form $F$ of genus 2, the authors compute the one-level density of the spinor $L$-functions of $F$, with a weight $\omega ^F$ which is essentially the modulus square of the first Fourier coefficient\footnote[4]{I.e. the Fourier coefficient corresponing to the identity matrix.} of $F$. This family is expected to be orthogonal, but with this weight one does not obtain the usual kernel $W_O$. This discrepancy can be explained by B\"ocherer's conjecture \cite{Bocherer, surveyBocherer} (now proved by Furusawa and Morimoto \cite{FM}), which claims that $\omega^F$ is proportional to the central value $L(\frac{1}{2},F)$. To be more precise, it says that $\omega^F\approx L(\frac{1}{2},F)L(\frac{1}{2},F\times\chi_4)$. Since $L(\frac{1}{2},F\times\chi_4)$ is \lq\lq uncorrelated'' with $L(s,F)$ and with its zeros, then the kernel they obtained is indeed $W_{SO^+}^1$ (see e.g. equation (\ref{aggiuntaperkowalski})\footnote[5]{In \cite{KowSahaTsim} the kernel is written as $1-\frac{\delta_0}{2}$, which is equivalent to $W_{SO^+}^1$ for test functions whose Fourier transforms are supported in $[-1,1]$, which is an assumption in \cite{KowSahaTsim}.} and note that weighting with $L(\frac{1}{2},F)^k$ the odd part of the family does not contribute, if $k>0$).  Moreover, they notice that this kernel is the one that arises from symplectic symmetry types. Thus, the symmetry of the family jumps from $O$ to $U\!Sp$, after weighting with $\omega^F$ (see also \cite{KnightlyReno, Sugiyama} for other examples where this phenomenon of change of symmetry type is observed). This transition can be seen as a particular case of equation (\ref{jonmie1}) below, which conjecturally predicts a relation between the weighted one-level density functions of different symmetry types. 

%subject to a specific weighting $\omega_k^F$ comparable to the one in (\ref{intro4.3}), i.e. defined in terms of values of $L$-functions at special points (in \cite{KowSahaTsim}, see (1.1.1) for the definition of $\omega_k^F$ and Theorem 1.2 for the distribution of low-lying zeros). Their analysis brings out that the weight $\omega_k^{F}$ \emph{contains arithmetic information related to central $L$-values of the Siegel cusp forms} and gives global evidence in favour of B\"ocherer's conjecture (see \cite{Bocherer} or e.g. \cite{surveyBocherer}).\\

\section{Statement of main results}\label{C4S1.5}

In the following, we focus on three specific families of $L$-functions, each with a different symmetry type; first we consider the unitary family $\boldsymbol\zeta:=\{\zeta(\tfrac{1}{2}+it):t\in\mathbb R\}$, i.e. the continuous family of the Riemann zeta function parametrized by a vertical shift. Then we study the symplectic family $\boldsymbol{L}_\chi$ of quadratic Dirichlet $L$-functions. Finally we look at the orthogonal family $\boldsymbol{L}_{\Delta,\chi}$ of the quadratic twists of the $L$-function associated with the discriminant modular form $\Delta$. For these families, under the assumption of the relevant Riemann Hypothesis and ratio conjecture, we perform an asymptotic analysis of $\mathcal D^{\mathcal F}_k(f)$. %; in the unitary case we also require that $k$ is even, which is a typical condition for this symmetry type 
Our results confirm our prediction  (\ref{intro4.4}), for small values of $k$. We recall that the case $k=0$ is already known in the literature for all of these families, both assuming the ratio conjecture (see \cite{CSapplications}) and without (for restricted ranges for $f$, see e.g. \cite{HR, MSConrey, OS}).%, as well as for many other families (see e.g. \cite{ILS, DueMill, FH}).
%\begin{equation}\label{intro4.4} 
%S^{\mathcal F}_k(f)=\int_{-\infty}^{+\infty}f(x)W_{G}^{(k)}(x)dx+O\Big(\frac{1}{\log X}\Big)
%\end{equation}
%for $k\leq 4$ (and even in the unitary case). Note that the superscript $^{(k)}$ indicates that we weight with the second power of $L(1/2)$; it particular it does not indicate a derivative.

We start with the unitary family. Note that, since this is a continuous family, the average over the family in the definition of $\mathcal D_k^{\boldsymbol\zeta}(f)$ is given by an integration over $t\in[T,2T]$ instead of the sum in (\ref{intro4.3}). In this case, setting
\begin{equation}\notag{W}_U^{0}(x):= W_U(x)=1,\end{equation}
\begin{equation}\notag{W}_U^{1}(x):= 1-\frac{\sin^2(\pi x)}{(\pi x)^2}\end{equation}
and
\begin{equation}\notag
{W}_U^{2}(x):=1-\frac{2+\cos(2\pi x)}{(\pi x)^2}+\frac{3\sin(2\pi x)}{(\pi x)^3}+\frac{3(\cos(2\pi x)-1)}{2(\pi x)^4},
\end{equation}
we prove the following theorem.
\begin{thm}\label{thm1LDunitary}
Let us assume the Riemann Hypothesis and the ratio conjecture (see Conjecture \ref{ratioconj}). Let us consider a test function $f$, which is holomorphic throughout the strip $|\Im(z)|<2$, real on the real line, even and such that $f(x)\ll1/(1+x^2)$ as $x\to\infty$. Then, for $k\in\mathbb N$ and $k\leq 2$, we have
\begin{equation}\notag
\mathcal D_k^{\boldsymbol\zeta}(f)=\int_{-\infty}^{+\infty}f(x)W_U^k(x)dx+O\Big(\frac{1}{\log T}\Big).
\end{equation}
\end{thm}

For this unitary family, in \cite{3.} we also develop an alternative method built on Hughes-Rudnick's technique in \cite{HR}, which allows us to show (\ref{intro4.4}) unconditionally\footnote[6]{Neither the Riemann Hypothesis nor the ratio conjecture is required. However, this unconditional strategy works only for test functions whose Fourier transform's support is small enough.}. Moreover, the analogue of Theorem \ref{thm1LDunitary} can be proved in the random matrix theory setting without any assumptions, since the formula for the ratios of characteristic polynomials averaged over the unitary group is known unconditionally (see \cite[Theorem 4.1]{CFZratioconj} and also \cite{CFZsupersymmetry, HuckPutt}). Therefore, denoting 
\begin{equation}\label{v2.1} Z=Z(A,\theta)=\det(I-Ae^{i\theta}) \end{equation}
for the characteristic polynomial of $N\times N$ matrices $A$ and $\theta_1,\dots,\theta_N$ for the phases of the eigenvalues of $A$, we prove a result that is the random matrix analogy of Theorem 2.1, essentially with the same proof. When we work on the random matrix theory side, to ensure that the one-level density is well-defined, we need the test function to be $2\pi$-periodic. Hence, given $f:\mathbb R\to \mathbb R$ an even Schwartz function, we define
\begin{equation}\label{well-defined}
F_N(x) =\sum_{h\in\mathbb Z}f\left( \frac{N}{2\pi}(x-2\pi h) \right)
\end{equation}  
and we prove the following.
%\begin{thm}\label{analogoRMTunitario}
%Let us consider a test function $f$, which is holomorphic throughout the strip $|\Im(z)|<2$, real on the real line, even and such that $f(x)\ll1/(1+x^2)$ as $x\to\infty$. 
%Then, for $k\in\mathbb N$ and $k\leq 2$, we have
%\begin{equation}\notag
%\frac{1}{\int_{U(N)}|Z|^{2k}d_{Haar}}\int_{U(N)}\sum_{j=1}^N f\Big(\frac{N}{2\pi}\theta_j\Big)|Z|^{2k}d_{Haar}
%\stackrel{N\to\infty}{\longrightarrow}
%\int_{-\infty}^{+\infty}f(x)W_{U}^{k}(x)dx.
%\end{equation}
%\end{thm}
\begin{thm}\label{analogoRMTunitario}
Let us consider $f:\mathbb R\to \mathbb R$ an even Schwartz function and $F_N$ as in \eqref{well-defined}.
Then, for $k\in\mathbb N$ and $k\leq 2$, we have
\begin{equation}\notag
\frac{1}{\int_{U(N)}|Z|^{2k}d_{Haar}}\int_{U(N)}\sum_{j=1}^N F_N(\theta_j)|Z|^{2k}d_{Haar}
\stackrel{N\to\infty}{\longrightarrow}
\int_{-\infty}^{+\infty}f(x)W_{U}^{k}(x)dx.
\end{equation}
\end{thm}

In the symplectic case, we compute the weighted one-level density functions for any non-negative integer $k\leq 4$. We set
%\begin{equation}\notag
%W_{U\!Sp}^{1}(x):=1+\frac{\sin(2\pi x)}{2\pi x}-\frac{2\sin^2(\pi x)}{(\pi x)^2},
%\end{equation}
%\begin{equation}\notag
%W_{U\!Sp}^{2}(x):=
%1-\frac{\sin(2\pi x)}{2\pi x}
%-\frac{24(1-\sin^2(\pi x))}{(2\pi x)^2}
%+\frac{48\sin(2\pi x)}{(2\pi x)^3}
%-\frac{96\sin^2(\pi x)}{(2\pi x)^4},
%\end{equation}
%\begin{equation}\begin{split}\notag
%W_{U\!Sp}^{(3)}(x)=
%1&+\frac{\sin(2\pi x)}{2\pi x}-\frac{12\sin^2(\pi x)}{(\pi x)^2}-\frac{240\sin(2\pi x)}{(2\pi x)^3}
%\\&-\frac{15(6-10\sin^2(\pi x))}{(\pi x)^4}+\frac{2880\sin(2\pi x)}{(2\pi x)^5}-\frac{90\sin^2(\pi x)}{(\pi x)^6},
%\end{split}\end{equation}
%\begin{equation}\begin{split}\notag
%W_{U\!Sp}^{(4)}(x)=
%1&-\frac{\sin(2\pi x)}{2\pi x}-\frac{10(1+\cos(2\pi x))}{(\pi x)^2}+\frac{90\sin(2\pi x)}{(\pi x)^3}
%\\&-\frac{15(3-31\cos(2\pi x))}{(\pi x)^4}-\frac{1470\sin(2\pi x)}{(\pi x)^5}
%\\&-\frac{315(1+9\cos(2\pi x))}{(\pi x)^6}+\frac{3150\sin(2\pi x)}{(\pi x)^7}-\frac{1575(1-\cos(2\pi x))}{(\pi x)^8}
%\end{split}\end{equation}
\begin{equation}\begin{split}\notag
&W_{U\!Sp}^{0}(x):=W_{U\!Sp}(x)=1-\frac{\sin(2\pi x)}{2\pi x},
\\&W_{U\!Sp}^{1}(x):=1+\frac{\sin(2\pi x)}{2\pi x}-\frac{2\sin^2(\pi x)}{(\pi x)^2},
\\&
W_{U\!Sp}^{2}(x):=
1-\frac{\sin(2\pi x)}{2\pi x}
-\frac{24(1-\sin^2(\pi x))}{(2\pi x)^2}
+\frac{48\sin(2\pi x)}{(2\pi x)^3}
-\frac{96\sin^2(\pi x)}{(2\pi x)^4},
\\&
W_{U\!Sp}^{3}(x):=
1+\frac{\sin(2\pi x)}{2\pi x}-\frac{12\sin^2(\pi x)}{(\pi x)^2}-\frac{240\sin(2\pi x)}{(2\pi x)^3}
\\&\hspace{2cm}-\frac{15(6-10\sin^2(\pi x))}{(\pi x)^4}+\frac{2880\sin(2\pi x)}{(2\pi x)^5}-\frac{90\sin^2(\pi x)}{(\pi x)^6},
\\&
W_{U\!Sp}^{4}(x):=
1-\frac{\sin(2\pi x)}{2\pi x}-\frac{10(1+\cos(2\pi x))}{(\pi x)^2}+\frac{90\sin(2\pi x)}{(\pi x)^3}
\\&\hspace{2cm}-\frac{15(3-31\cos(2\pi x))}{(\pi x)^4}-\frac{1470\sin(2\pi x)}{(\pi x)^5}
\\&\hspace{2cm}-\frac{315(1+9\cos(2\pi x))}{(\pi x)^6}+\frac{3150\sin(2\pi x)}{(\pi x)^7}-\frac{1575(1-\cos(2\pi x))}{(\pi x)^8}
\end{split}\end{equation}
and we prove the following result.
\begin{thm}\label{thm1LDsymplectic}
Let us assume the Riemann Hypothesis and the ratio conjecture for the $L$-functions in the family $\boldsymbol L_\chi$ (see Conjecture \ref{ratioconjsympl}). Let us consider a test function $f$, which is holomorphic throughout the strip $|\Im(z)|<2$, real on the real line, even and such that $f(x)\ll1/(1+x^2)$ as $x\to\infty$. Then, for $k\in\mathbb N$ and $k\leq 4$, we have
\begin{equation}\notag
\mathcal D_k^{\boldsymbol{L}_\chi}(f)=\int_{-\infty}^{+\infty}f(x)W_{U\!Sp}^k(x)dx+O\Big(\frac{1}{\log X}\Big).
\end{equation}
\end{thm}
Also in the symplectic case, with the same proof we also get the corresponding result in the random matrix theory setting unconditionally, as \cite[Theorem 4.2]{CFZratioconj} provides the analogue of Conjecture \ref{ratioconjsympl}. Note that, in the symplectic (resp. orthogonal) case, $Z$ defined as in \eqref{v2.1} is the characteristic polynomial of $2N\times 2N$ symplectic (resp. orthogonal) matrices.
%\begin{thm}\label{analogoRMTsimplettico}
%Let us consider a test function $f$, which is holomorphic throughout the strip $|\Im(z)|<2$, real on the real line, even and such that $f(x)\ll1/(1+x^2)$ as $x\to\infty$. Then, for $k\in\mathbb N$ and $k\leq 4$, we have
%\begin{equation}\notag
%\frac{1}{\int_{U\!Sp(2N)}Z^kd_{Haar}}\int_{U\!Sp(2N)}\sum_{j=1}^N f\Big(\frac{N}{\pi}\theta_j\Big)Z^kd_{Haar}
%\stackrel{N\to\infty}{\longrightarrow}
%\int_{-\infty}^{+\infty}f(x)W_{U\!Sp}^{k}(x)dx.
%\end{equation}
%\end{thm}
\begin{thm}\label{analogoRMTsimplettico}
Let us consider $f:\mathbb R\to \mathbb R$ an even Schwartz function and $F_N$ as in \eqref{well-defined}.
Then, for $k\in\mathbb N$ and $k\leq 4$, we have
\begin{equation}\notag
\frac{1}{\int_{U\!Sp(2N)}Z^kd_{Haar}}\int_{U\!Sp(2N)}\sum_{j=1}^N F_{2N}(\theta_j)\;Z^kd_{Haar}
\stackrel{N\to\infty}{\longrightarrow}
\int_{-\infty}^{+\infty}f(x)W_{U\!Sp}^{k}(x)dx.
\end{equation}
\end{thm}

Finally, for the (even) orthogonal family $\boldsymbol L_{\Delta,\chi}$, we denote
\begin{equation}\begin{split}\notag
&W_{SO^+}^{0}(x):={W}_{SO^+}(x) = 1+\frac{\sin(2\pi x)}{2\pi x},
\\&W_{SO^+}^{1}(x):=1-\frac{\sin(2\pi x)}{2\pi x},
\\&
W_{SO^+}^{2}(x):=1+\frac{\sin(2\pi x)}{2\pi x}-\frac{2\sin^2(\pi x)}{(\pi x)^2},
\\&
W_{SO^+}^{3}(x):=1-\frac{\sin(2\pi x)}{2\pi x}
-\frac{24(1-\sin^2(\pi x))}{(2\pi x)^2}
+\frac{48\sin(2\pi x)}{(2\pi x)^3}
-\frac{96\sin^2(\pi x)}{(2\pi x)^4},
\\&
W_{SO^+}^{4}(x):=
1+\frac{\sin(2\pi x)}{2\pi x}-\frac{12\sin^2(\pi x)}{(\pi x)^2}-\frac{240\sin(2\pi x)}{(2\pi x)^3}
\\&\hspace{2cm}-\frac{15(6-10\sin^2(\pi x))}{(\pi x)^4}+\frac{2880\sin(2\pi x)}{(2\pi x)^5}-\frac{90\sin^2(\pi x)}{(\pi x)^6}.
\end{split}\end{equation}
%\begin{equation}\notag W_{SO^+}^{0}(x):= {W}_{SO^+}(x) = 1+\frac{\sin(2\pi x)}{2\pi x}\end{equation}
%and 
%\begin{equation}\notag W_{SO^+}^{k}(x):= {W}_{U\!Sp}^{k-1}(x)\end{equation}
%for $k\in\mathbb Z_+$, $k\leq 4$.
Notice that there are strong similarities with the symplectic kernels; we will discuss these analogies below.
With these notations, we prove the following theorem.
\begin{thm}\label{thm1LDorthogonal}
Let us assume the Riemann Hypothesis and the ratio conjecture for the $L$-functions in the family $\boldsymbol L_{\Delta,\chi}$ (see Conjecture \ref{RCorthogonal}). Let us consider a test function $f$, which is holomorphic throughout the strip $|\Im(z)|<2$, real on the real line, even and such that $f(x)\ll1/(1+x^2)$ as $x\to\infty$. Then, for $k\in\mathbb N$ and $k\leq 4$, we have
\begin{equation}\notag
\mathcal D_k^{\boldsymbol{L}_{\Delta,\chi}}(f)=\int_{-\infty}^{+\infty}f(x)W_{SO^+}^k(x)dx+O\Big(\frac{1}{\log X}\Big).
\end{equation}
\end{thm}
Again the analogous result in random matrix theory is instead unconditional (relying on \cite[Theorem 4.3]{CFZratioconj} in place of Conjecture \ref{RCorthogonal}).
%\begin{thm}\label{analogoRMTortogonale}
%Let us consider a test function $f$, which is holomorphic throughout the strip $|\Im(z)|<2$, real on the real line, even and such that $f(x)\ll1/(1+x^2)$ as $x\to\infty$. Then, for $k\in\mathbb N$ and $k\leq 4$, we have
%\begin{equation}\notag
%\frac{1}{\int_{SO(2N)}Z^kd_{Haar}}\int_{SO(2N)}\sum_{j=1}^N f\Big(\frac{N}{\pi}\theta_j\Big)Z^kd_{Haar}
%\stackrel{N\to\infty}{\longrightarrow}
%\int_{-\infty}^{+\infty}f(x)W_{SO^+}^{k}(x)dx.
%\end{equation}
%\end{thm}
\begin{thm}\label{analogoRMTortogonale}
Let us consider $f:\mathbb R\to \mathbb R$ an even Schwartz function and $F_N$ as in \eqref{well-defined}. 
Then, for $k\in\mathbb N$ and $k\leq 4$, we have
\begin{equation}\notag
\frac{1}{\int_{SO(2N)}Z^kd_{Haar}}\int_{SO(2N)}\sum_{j=1}^N F_{2N}(\theta_j)\;Z^kd_{Haar}
\stackrel{N\to\infty}{\longrightarrow}
\int_{-\infty}^{+\infty}f(x)W_{SO^+}^{k}(x)dx.
\end{equation}
\end{thm}

%%%%%%%%%%%%%%%%%%%%%%%%%%%%%%%%%%%%%%%%%%%%%%%%%%%%%%%%%%%%%%%%%%%%%%%%%%%%%%%%%%%%%%%%%%%%%%%%%

\subsection{A general Conjecture for $W_G^k$}\label{C4S2sottosezione1}

Thanks to the explicit expressions we get for the kernels $W_G^{k}$ in the range $k\leq 4$, we can speculate about what happens for any $k\in\mathbb N$.
First of all, we notice that the Fourier transform of the kernels $W_G^{k}$ exhibits a structure. From the explicit formulae for $W_G^{k}$ we get in the range $k\leq 4$, $\widehat W_G^{k}$ turns out to be an even function, supported on $[-1,1]$, uniquely determined by a polynomial on $[0,1]$. %The degree of this polynomial and its value at 0 depend on $k$ and on the symmetry type of the family. 
More precisely, we conjecture that 
\begin{equation}\label{ftconj1000}
\widehat W_G^k(y)=\delta_0(y)+P_G^k(|y|)\chi_{[-1,1]}(y)
\end{equation}
where $P_G^k$ is a polynomial depending on $k$ and $G$ only. In particular, in the unitary case and with $k\geq 1$, we expect the degree of $P_U^k$ to be $2k-1$ and $P_U^k(0)=-k$, $P_U^k(1)=0$. For the symplectic family, if $k\geq 1$, we predict $P_{U\!Sp}^k$ with degree $2k-1$ and $P_{U\!Sp}^k(0)=-(2k+1)/2$, $P_{U\!Sp}^k(1)=(-1)^{k+1}/2$. Finally for the orthogonal symmetry type, we conjecture the degree of $P_{SO^+}^k$ to be $2k-3$ and $P_{SO^+}^k(0)=-(2k-1)/2$, $P_{SO^+}^k(1)=(-1)^k/2$ for any $k\geq 2$ (the case $k=1$ yields $\widehat W_{SO^+}^{1}(y)=\delta_0(y)-1/2$).
We collect into a table all the values of $P_G^k$ we obtained for $k$ small, which support our speculations. 
Note that the case $k=0$, corresponding to the first row in the table, was already known in the literature, while all other results are new.% We should also specify that the last row of the table is instead conjectural.  
\\

\noindent
\begin{tabular}{p{1.1cm}|p{3.5cm}|p{3.5cm}|p{3.5cm}|}
\rule[-5mm]{0mm}{1.2cm}
\hfil $P_G^k $& \hfil$G=U$ & \hfil$G=U\!Sp$ & \hfil$G=SO^+$   
\\ \hline 
\rule[-6mm]{0mm}{1.5cm}
\hfil$k=0$ 
& \hfil $0$
& \hfil $\displaystyle -\frac{1}{2}$
& \hfil $\displaystyle \frac{1}{2}$
\\ \hline 
\rule[-6mm]{0mm}{1.5cm}
\hfil$k=1$ 
& \hfil$\displaystyle y-1$
& \hfil$\displaystyle 2y-\frac{3}{2}$
& \hfil$\displaystyle -\frac{1}{2}$
\\ \hline 
\rule[-6mm]{0mm}{1.5cm}
\hfil$k=2$ 
& \hfil $\displaystyle -2y^3+4y-2$ 
& \hfil $\displaystyle -4y^3+6y-\frac{5}{2}$
& \hfil $\displaystyle 2y-\frac{3}{2}$
\\ \hline 
\rule[-6mm]{0mm}{1.5cm}
\hfil$k=3$ 
& 
& \hfil $\displaystyle 12y^5-20y^3+12y-\frac{7}{2}$ 
& \hfil $\displaystyle -4y^3+6y-\frac{5}{2}$
\\ \hline 
\rule[-1cm]{0mm}{1.7cm}
\hfil$k=4$ 
& 
& \hfil $\displaystyle -40y^7+84y^5-60y^3+20y-\frac{9}{2}$
& \hfil $\displaystyle 12y^5-20y^3+12y-\frac{7}{2}$ 
\\ \hline
%\rule[-1cm]{0mm}{1.5cm}
%conj.  \normalsize
%&\hfil $\deg(P_{U}^k)=2k-1$\;\;\; \par $P_U^k(0)=-k$ \par $P_U^k(1)=0$ \par ($k\in\mathbb Z_+$) 
%&\hfil $\deg(P_{U\!Sp}^k)=2k-1$ \par $P_{U\!Sp}^k(0)=-\frac{2k+1}{2}$  \par $P_{U\!Sp}^k(1)=\frac{(-1)^{k+1}}{2}$ ($k\in\mathbb Z_+$) 
%&\hfil $\deg(P_{SO^+}^k)=2k-3$ \par $P_{SO^+}^k(0)=-\frac{2k-1}{2}$ \par $P_{SO^+}^k(1)=\frac{(-1)^{k}}{2}$ ($k\in\mathbb N, k\geq 2$)
%\\ \hline
\end{tabular}\\ \\

\begin{figure}[h]
\includegraphics[width=13.7cm]{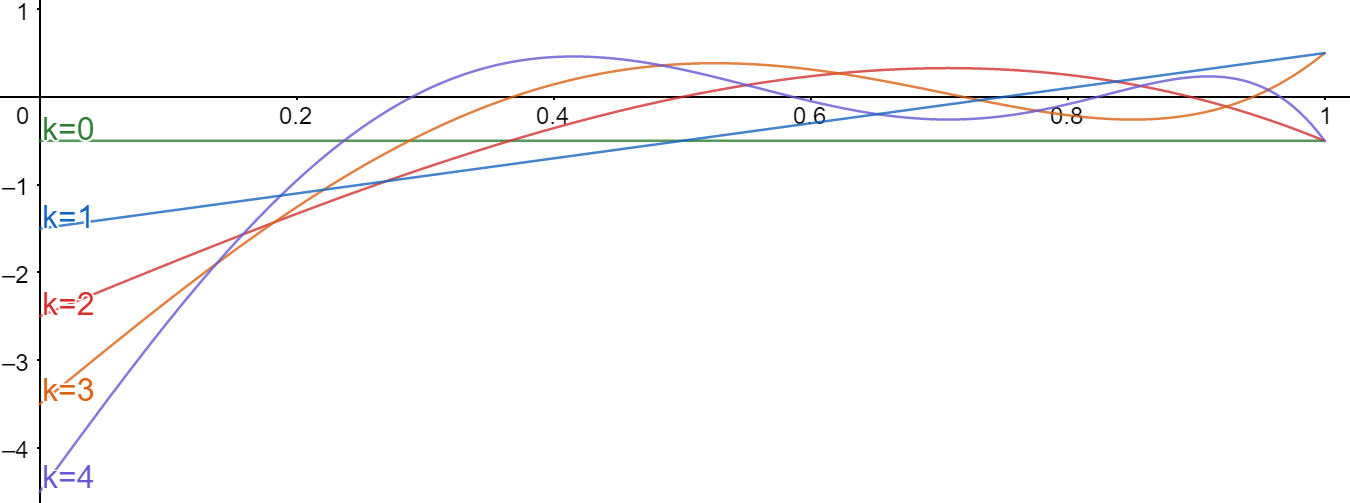}
\caption{$P_{U\!Sp}^k(y)$, for $y\in[0,1]$.}
\end{figure}

Looking at the table, we can detect relations between the weighted one-level density functions with different symmetry types. In particular, from the above discussion, it seems natural to expect that 
\begin{equation}\label{jonmie1} W_{SO^+}^{k}(x)=W_{U\!Sp}^{k-1}(x) \end{equation}
for any $k\in\mathbb Z_+$. Moreover, the Fourier transforms of $W_G^k$ suggest that the weighted one-level density function in the unitary case is the average of the symplectic and orthogonal cases; namely we conjecture that
\begin{equation}\label{jonmie2} W_{U}^{k}(x)=\frac{W_{U\!Sp}^{k}(x)+W_{SO^+}^k(x)}{2}. \end{equation}
We note that also the leading order moment coefficients $f_G(k)$ for the three compact groups $U,U\!Sp,SO^+$ satisfy relations linking them with each other, being (see \cite{MSJon}, Equations (6.10) and (6.11))
$$  f_{SO^+}(k)=2^kf_{U\!Sp}(k-1) \quad \text{and}\quad 
2^{k^2}f_{U}(k)=f_{U\!Sp}(k)f_{SO^+}(k).$$
Equations (\ref{jonmie1}) and (\ref{jonmie2}) can be seen as the analogue of the above formulae, in the context of the weighted one-level density.

Finally we conjecture an explicit formula for the polynomials $P_G^k$, which together with (\ref{ftconj1000}) provides a precise conjecture for the weighted kernels $W_G^k$. In view of equations (\ref{jonmie1}) and (\ref{jonmie2}), it suffices to focus on the symplectic case only. 
Looking at what happens for $k\leq 4$, we speculate that for every positive integer $k$ we have 
%\begin{equation}\label{poly101}\frac{d}{dy} P_{U\!Sp}^k(y) = -k(k+1)\sum_{j=1}^k t(k,j)y^{2j-2}\end{equation}
%where the coefficient $t(k,j)$ is defined by 
%$$ t(k,j)=(-1)^{j}T(k+2,j-1) $$
%for $k\geq 1$ and $j=1,\dots,k$ and 
%$$ T(m,n)=\frac{1}{n+1}\binom{m-3}{n}\binom{m+n-1}{n} $$
%is the number of diagonal dissections of a convex $m$-gon into $n+1$ regions.
% Since we expect $P_{U\!Sp}^k(0)=-\frac{2k+1}{2}$, equation (\ref{poly101}) leads to believe that
%\begin{equation}\label{poly102}\notag
%P_{U\!Sp}^k(y) =-\frac{2k+1}{2}- k(k+1)\sum_{j=1}^k t(k,j)\frac{y^{2j-1}}{2j-1}.
%\end{equation}
\begin{equation}\label{poly101}P_{U\!Sp}^k(y) =-\frac{2k+1}{2} -k(k+1)\sum_{j=1}^k(-1)^j c_{j,k}\frac{y^{2j-1}}{2j-1}\end{equation}
where the coefficient $c_{j,k}$ is defined by 
$$c_{j,k}=\frac{1}{j}\binom{k-1}{j-1}\binom{k+j}{j-1}. $$
We note that the sequence of the $c_{j,k}$'s appears in OEIS\footnote[7]{https://oeis.org/A033282.}, as the number of diagonal dissections of a convex $(k+2)$-gon into $j$ regions.
By Fourier inversion, from (\ref{ftconj1000}) and (\ref{poly101}), we get an explicit conjectural formula for $W_{U\!Sp}^k$, being
$$ W_{U\!Sp}^k(x)=1-(2k+1)\frac{\sin(2\pi x)}{2\pi x}+\sum_{j=1}^k\frac{k(k+1)}{2^{2j-2}\pi^{2j-1}}\frac{c_{j,k}}{2j-1}\frac{d^{2j-1}}{dx^{2j-1}}\bigg[\frac{1-\cos(2\pi x)}{2\pi x}\bigg]. $$

\begin{figure}[b]
\includegraphics[width=13.7cm]{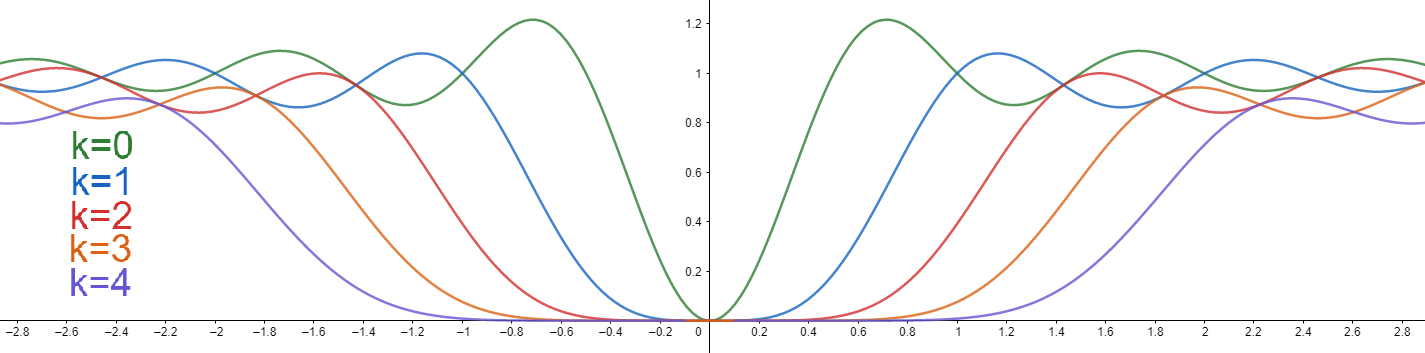}
\caption{$W_{U\!Sp}^k(x)$ for $k=0,\dots,4$.}
\end{figure}

From all these discussions, we can formulate the following conjecture.
\begin{conj}\label{congetturadecisiva1}
Let us consider a test function $f$, holomorphic in the strip $|\Im(z)|<2$, even, real on the real line and such that $f(x)\ll 1/(1+x^2)$ as $x\to\infty$, then for any $k\in\mathbb N$. Given a family $\mathcal F$ of $L$-functions with symmetry type $G\in\{U,U\!Sp,SO^+\}$, we have
\begin{equation}\notag
\mathcal D_k^{\mathcal F}(f)=\int_{-\infty}^{+\infty}f(x)W_G^k(x)dx+O\Big(\frac{1}{\log X}\Big)
\end{equation}
as $X\to\infty$, where the weighted one-level density function $W_G^k$ depends on $k$ and $G$ only. In addition the following relations hold
$$W_{SO^+}^{k}(x)=W_{U\!Sp}^{k-1}(x)  \quad \text{and} \quad
W_{U}^{k}(x)=\frac{W_{U\!Sp}^{k}(x)+W_{SO^+}^k(x)}{2}  $$
for any $k\in\mathbb Z_+$ and $k\in\mathbb N$ respectively. Moreover, for every $k\in\mathbb Z_+$, in the symplectic case (the others can be recovered by the above relations), we have that
$$ \widehat W_{U\!Sp}^k(y)=\delta_0(y)+P_{U\!Sp}^k(|y|)\chi_{[-1,1]}(y) $$
where $P_{U\!Sp}^k$ is a polynomial of degree $2k-1$, given by
$$P_{U\!Sp}^k(y)=
-\frac{2k+1}{2}
-k(k+1)\sum_{j=1}^{k}(-1)^jc_{j,k}\frac{y^{2j-1}}{2j-1},$$ 
with
$$c_{j,k}=\frac{1}{j}\binom{k-1}{j-1}\binom{k+j}{j-1}.$$ 
%In particular, as $x\to0$, we have
%$$ W_{U\!Sp}^k(x)\sim \frac{2\pi^{2(k+1)}}{(2k+1)!!(2k+3)!!}x^{2(k+1)}.$$
\end{conj}

\subsection{An expression for $W_G^k(x)$ in terms of hypergeometric functions and its vanishing at $x=0$}\label{C4S2sottosezione2}
 
We now focus on the behaviour of the weighted kernels $W_G^k(x)$ at $x= 0$. For all symmetry types, it seems clear that the order of vanishing of $W_{G}^{k}(x)$ for $x\to0$ increases as $k$ grows. This phenomenon reflects the effect of the weight $V(L(1/2))^k$ in the average over the family, which gives more and more relevance to those $L$-functions that are large at the central point, as $k$ increases. % and, on the other hand, kills the contribution of those which are small therein. 
More precisely, for the unitary family we conjecture that 
\begin{equation}\label{asintoticacongetturale}W_U^{k}(x)\sim \frac{\pi^{2k}x^{2k}}{(2k-1)!!(2k+1)!!}\end{equation}
as $x\to0$, $k\in\mathbb N$. In particular, together with (\ref{intro4.4}), this suggests that, on weighted average over the considered family, the number of normalized zeros which are less than $\varepsilon$ away from the central point is typically $\asymp_k\varepsilon^{2k+1}$. 
Analogously, the asymptotic behaviour of the symplectic and orthogonal kernels can be deduced from (\ref{asintoticacongetturale}) by equations (\ref{jonmie1}) and (\ref{jonmie2}).
%Analogously, in the symplectic and orthogonal cases we conjecture that 
%$$W_{U\!Sp}^{k}(x)\sim \frac{2\pi^{2(k+1)}x^{2(k+1)}}{(2k+1)!!(2k+3)!!}$$
%and
%$$W_{SO^+}^{k}(x) \sim \frac{2\pi^{2k}x^{2k}}{(2k-1)!!(2k+1)!!}$$
%as $x\to0$, $k\in\mathbb N$. 
%questi levali, ci metti solo che seguono dalle equazioni (...) e (...). Poi nella conjettura li scrivi bene tutti i 3
For small values of $k$, the behaviour of $W_G^k(x)$ at $x=0$ is outlined in the following table; the first row was already known in the literature, all the others are new.%, while the last one is conjectural.
\\ \\
\noindent
\begin{tabular}{p{1.1cm}|p{3.5cm}|p{3.5cm}|p{3.5cm}|}
\rule[-5mm]{0mm}{1.2cm}
$W_G^k$\par $x\to0$ & \hfil$G=U$ & \hfil$G=U\!Sp$ & \hfil$G=SO^+$   
\\ \hline 
\rule[-7mm]{0mm}{1.7cm}
\hfil$k=0$ 
& \hfil 1
& \hfil $\displaystyle\frac{2\pi^2x^2}{3}$
& \hfil 2
\\ \hline 
\rule[-7mm]{0mm}{1.7cm}
\hfil$k=1$ 
& \hfil$\displaystyle\frac{\pi^2x^2}{3}$
& \hfil$\displaystyle\frac{2\pi^4x^4}{45}$
& \hfil$\displaystyle\frac{2\pi^2x^2}{3}$
\\ \hline 
\rule[-7mm]{0mm}{1.7cm}
\hfil$k=2$ 
& \hfil $\displaystyle\frac{\pi^4x^4}{45}$
& \hfil$\displaystyle\frac{2\pi^6x^6}{1575}$
& \hfil$\displaystyle\frac{2\pi^4x^4}{45}$
\\ \hline 
\rule[-7mm]{0mm}{1.7cm}
\hfil$k=3$ 
& 
& \hfil $\displaystyle\frac{2\pi^8x^8}{99225}$ 
& \hfil $\displaystyle\frac{2\pi^6x^6}{1575}$
\\ \hline 
\rule[-7mm]{0mm}{1.7cm}
\hfil$k=4$ 
& 
& \hfil $\displaystyle\frac{2\pi^{10}x^{10}}{9823275}$
& \hfil $\displaystyle\frac{2\pi^8x^8}{99225}$
\\ \hline
%\rule[-7mm]{0mm}{1.7cm}
%conj.
%& \hfil$\displaystyle\frac{\pi^{2k}x^{2k}}{(2k-1)!!(2k+1)!!}$
%& \hfil$\displaystyle\frac{2\pi^{2(k+1)}x^{2(k+1)}}{(2k+1)!!(2k+3)!!}$
%& \hfil$\displaystyle\frac{2\pi^{2k}x^{2k}}{(2k-1)!!(2k+1)!!}$
%\\ \hline
\end{tabular}\\ \\

%\begin{figure}[h]
%\includegraphics[width=13.7cm]{at0(3)_sympl.png}
%\caption{$W_{U\!Sp}^k(x)$, as $x\to0$.}
%\end{figure}

In the following conjecture, we condense all the speculations about the behaviour of the weighted kernels $W_G^k(x)$ as $x\to0$.

\begin{conj}\label{congetturadecisiva2}
For $G\in\{U,U\!Sp,SO^+\}$, $k\in\mathbb N$, the weighted kernels $W_G^k$ defined in Conjecture \ref{congetturadecisiva1} satisfy the following asymptotic relations as $x\to0$:
\begin{equation}\begin{split}\notag
W_U^{k}(x)&\sim \frac{\pi^{2k}x^{2k}}{(2k-1)!!(2k+1)!!}
\\W_{U\!Sp}^{k}(x)&\sim \frac{2\pi^{2(k+1)}x^{2(k+1)}}{(2k+1)!!(2k+3)!!}
\\W_{SO^+}^{k}(x)& \sim \frac{2\pi^{2k}x^{2k}}{(2k-1)!!(2k+1)!!}.
\end{split}\end{equation}
\end{conj}

Finally, assuming Conjecture \ref{congetturadecisiva1}, we obtain the expansion of $W_G^k(x)$ at $x=0$. In particular, we show that the asymptotic behaviour of $W_G^k(x)$ can be deduced from the explicit formulae that we conjectured in Section \ref{C4S2sottosezione1}. In view of equations (\ref{jonmie1}) and (\ref{jonmie2}), it suffices to consider the symplectic case only.
\begin{thm}\label{conj1implicaconj2}  Let us assume Conjecture \ref{congetturadecisiva1}. Then for any $k\in\mathbb N$ we have
$$ W_{U\!Sp}^{k}(x)=\sum_{m=1}^\infty \beta_{m,k} x^{2m} $$
with
$$ \beta_{m,k}=(-1)^{m+1}\frac{(2\pi)^{2m}}{(2m+1)!}\bigg((-1)^k+\frac{k(k+1)}{m+1} \pFq{3}{2}{1-k,k+2,m+1}{m+2,2}{1}\bigg) $$ 
where ${}_3F_{2}$ denotes the generalized hypergeometric function. Moreover, we have 
$$ \pFq{3}{2}{1-k,k+2,m+1}{m+2,2}{1}=\begin{cases} \frac{(m+1)(-1)^{k+1}}{k(k+1)} &\text{if }1\leq m \leq k \\ \frac{2(-1)^{k+1}(k-1)!(k+2)!}{(2k+2)!}\big(\binom{2k+1}{k+1}-1\big) &\text{if } m=k+1.\end{cases} $$ 
In particular, Conjecture \ref{congetturadecisiva2} follows.
\end{thm}

\section{Proof of Theorem~\ref{thm1LDunitary} and Theorem \ref{analogoRMTunitario}}\label{C4S2}

We first tilt the Lebesgue measure multiplying by $|\zeta(1/2+it)|^2$ and denote 
\begin{equation}\label{tiltedmean}
\langle g\rangle_{|\zeta|^2}
%\mathcal D_1^{\boldsymbol\zeta}(f)
:= \frac{1}{T\log T}\int_T^{2T}g(t)|\zeta(1/2+it)|^2dt.
\end{equation}
then we consider $f$ an even test function and its Fourier transform 
\begin{equation}\notag\widehat{f}(y):=\int_{-\infty}^{+\infty}f(x)e^{-2\pi ixy}dx.\end{equation}

We recall that Conrey, Farmer and Zirnbauer \cite{CFZratioconj} applied a modification of the recipe for integral moments to the case of ratios getting the following statement, called the ratio conjecture (here we state this conjecture in a slightly weaker form than in \cite{CFZratioconj}, as far as the shifts are concerned).

\begin{conj}[\cite{CFZratioconj}, Conjecture 5.1]\label{ratioconj}
Let us denote $\chi(s)$ the explicit factor in the functional equation $\zeta(s)=\chi(s)\zeta(1-s)$. For any positive integers $K,L,Q,R$ and for any $\alpha_1,\dots,\alpha_{K+L},\gamma_1,\dots,\gamma_Q,\delta_1,\dots,\delta_R$ complex shifts with real part $\asymp(\log T)^{-1}$ and imaginary part $\ll_\varepsilon T^{1-\varepsilon}$ for every $\varepsilon>0$, then
\begin{equation}\begin{split}\label{2.3.17}\notag
\frac{1}{T}\int_T^{2T}&\frac{\prod_{k=1}^K\zeta(s+\alpha_k)\prod_{l=K+1}^{K+L}\zeta(1-s-\alpha_l)}{\prod_{q=1}^Q\zeta(s+\gamma_q)\prod_{r=1}^R\zeta(1-s+\delta_r)}dt\\&
=\frac{1}{T}\int_T^{2T}\sum_{\sigma\in\Xi_{K,L}}\prod_{k=1}^K\frac{\chi(s+\alpha_k)}{\chi(s-\alpha_{\sigma(k)})}Y_UA_\zeta(...)dt+O(T^{1/2+\varepsilon})\\&
\text{with } (...)=(\alpha_{\sigma(1)}, \dots,\alpha_{\sigma(K)};-\alpha_{\sigma(K+1)},\dots,-\alpha_{\sigma(K+L)};\gamma;\delta)
\end{split}\end{equation}
where
\begin{equation}\label{2.3.18}\notag
Y_U(\alpha;\beta;\gamma;\delta):=\frac{\prod_{k=1}^K\prod_{l=1}^L\zeta(1+\alpha_k+\beta_l)\prod_{q=1}^Q\prod_{r=1}^R\zeta(1+\gamma_q+\delta_r)}{\prod_{k=1}^K\prod_{r=1}^R\zeta(1+\alpha_k+\delta_r)\prod_{l=1}^L\prod_{q=1}^Q\zeta(1+\beta_l+\gamma_q)}
\end{equation}
and $A_\zeta$ is an Euler product, absolutely convergent for all of the variables in small disks
around 0, which is given by
\begin{equation}\begin{split}\label{2.3.19}\notag
A_\zeta(\alpha;\beta;&\gamma;\delta):=\prod_p\frac{\prod_{k=1}^K\prod_{l=1}^L(1-1/p^{1+\alpha_k+\beta_l})\prod_{q=1}^Q\prod_{r=1}^R(1-1/p^{1+\gamma_q+\delta_r})}{\prod_{k=1}^K\prod_{r=1}^R(1-1/p^{1+\alpha_k+\delta_r})\prod_{l=1}^L\prod_{q=1}^Q(1-1/p^{1+\beta_l+\gamma_q})} \\&
\sum_{\sum a_k+\sum c_q=\sum b_l+\sum d_r}\frac{\prod\mu(p^{c_q})\prod\mu(p^{d_r})}{p^{\sum(1/2+\alpha_k)a_k+\sum(1/2+\beta_l)b_l+\sum(1/2+\gamma_q)c_q+\sum(1/2+\delta_r)d_r}}
\end{split}\end{equation}
while $\Xi_{K,L}$ denotes the subset of permutations $\sigma \in S_{K+L}$ of $\{1, 2, \dots , K + L\}$ for which $\sigma(1)<\sigma(2) <\dots<\sigma(K)$ and $\sigma(K+1)<\sigma(K +2)<\dots<\sigma(K+L).$
\end{conj}
%With an assumption about the moments of the Riemann zeta function we can remove the extra condition about the support of $\widehat{f}$. In particular, on RH, if we assume the ratio conjecture (see Conjecture~\ref{ratioconj}), we can also handle the general case with $f$ a decaying function, without any additional condition on its Fourier transorm's support:
By assuming this conjecture about the moments of zeta, 
denoting $$ N_f(t)=\sum_\gamma f\left(\frac{\log T}{2\pi}(\gamma-t)\right), $$
we can prove the following result.
\begin{prop}\label{congettura}
Let us assume Conjecture~\ref{ratioconj} and the Riemann Hypothesis.
We consider a test function $f(z)$ which is holomorphic throughout the strip $|\Im(z)| < 2$, real on the real line, even and such that $f(x) \ll 1/(1 + x^2)$ as $x\to\infty$.
Then 
\begin{equation}\notag
\mathcal D_1^{\boldsymbol\zeta}(f):=\langle N_f\rangle_{|\zeta|^2}=\int_{-\infty}^{+\infty}{f}(x){W}_U^{1}(x)dx+O\Big(\frac{1}{\log T}\Big)
\end{equation}
with
\begin{equation}\notag{W}_U^{1}(x):= 1-\sinc^2(x)=1-\frac{\sin^2(\pi x)}{(\pi x)^2}.\end{equation}
\end{prop}
%
%The right hand side in Theorem~\ref{teorema} and Proposition~\ref{congettura} is revealing and it can be easily compared to the density conjecture in~\eqref{intro4.2} and to the classical mean~\eqref{1ldzetanew}. Indeed, ${W}_U^{1}(x)\sim\frac{\pi^2}{3}x^2$ vanishes at $x=0$ of order 2, showing a repulsion of zeros at height $t$ which does not occur in the classical case. This repulsion can be explained by the fact that the measure $|\zeta(1/2+it)|^2dt$ gives more weight to the large values of zeta, around which is more unlikely to find a zero.\par
In addition, with the same strategy as in the proof of Proposition~\ref{congettura} (but much longer computations, which can be done by using Sage\footnote[8]{SageMath, the Sage Mathematics Software System (Version 0.6.3), The Sage Developers, 2021, https://www.sagemath.org.}) we can also study the \lq\lq fourth moment\rq\rq case. Namely, we denote
\begin{equation}\label{Zonta.1}
\langle g\rangle_{|\zeta|^4}:= \frac{1}{\frac{1}{2\pi^2}T(\log T)^4}\int_T^{2T}g(t)|\zeta(1/2+it)|^4dt.
\end{equation}
and we prove the following result.

\begin{prop}\label{casozeta^4}
Let us assume Conjecture~\ref{ratioconj} and the Riemann Hypothesis.
We consider a test function $f(z)$ which is holomorphic throughout the strip $|\Im(z)| < 2$, real on the real line, even and such that $f(x) \ll 1/(1 + x^2)$ as $x\to\infty$.
Then 
\begin{equation}\label{aux400}
\mathcal D_{2}^{\boldsymbol\zeta}(f):=\langle N_f\rangle_{|\zeta|^4}=\int_{-\infty}^{+\infty}{f}(x){W}_U^{2}(x)dx+O\Big(\frac{1}{\log T}\Big)
\end{equation}
with
\begin{equation}\notag
{W}_U^{2}(x):=1-\frac{2+\cos(2\pi x)}{(\pi x)^2}+\frac{3\sin(2\pi x)}{(\pi x)^3}+\frac{3(\cos(2\pi x)-1)}{2(\pi x)^4}.
\end{equation}
%and in particular ${W}_U^{(4)}(x)\sim \frac{\pi^4}{45}x^4$ as $x\to0$. 
\end{prop}
%We note that Propositions~\ref{congettura} and~\ref{casozeta^4} together imply Theorem~\ref{thm1LDunitary}. 
%Finally we record the Fourier transform of the weighted kernels, being
%\begin{equation}\notag
%\widehat{W}_U^{1}(y)=
%\begin{cases}  0 & \text{if }|y|>1\\ -y-1 &\text{if }-1\leq y<0 \\ 0  &\text{if }y=0 \\ y-1 &\text{if }0< y\leq 1\end{cases}
%\end{equation}
%and
%\begin{equation}\notag
%\widehat{W}_U^{2}(y)=\begin{cases}  0 &\text{if }|y|>1\\ 2y^3-4y-2 &\text{if } -1\leq y<0 \\ -1 &\text{if }y=0 \\ -2y^3+4y-2 &\text{if }0<y\leq 1. \end{cases}
%\end{equation}
%We refer to the next chapter for a complete analysis of these kernels and for speculations about the general case with weight $|\zeta(1/2+it)|^{2k}$, $k\in\mathbb N$.

%%%%%%%%%%%%%%%%%%%%%%%%%%%%%%%%%%%%%%%%%%%%%%%%%%%%%%%%%%%%%%%%%%%%%%%%%%%%

\subsection{Proof of Proposition~\ref{congettura}}\label{C4S2S3}

To prove Proposition~\ref{congettura} we strongly rely on Conjecture~\ref{ratioconj}, which allows us to perform a similar computation as in Section 3 of \cite{CSapplications}. We introduce two parameters $\alpha,\beta\in\mathbb R$ of size $\asymp1/\log T$, we denote
\begin{equation}\label{aux4.00}
\boldsymbol\zeta_{\alpha,\beta}(t):=\zeta(1/2+\alpha+it)\zeta(1/2+\beta-it)
\end{equation}
and we look at 
\begin{equation}\label{aux4.0}
\langle N_f \rangle_{|\zeta|^2}^{\alpha,\beta}:=\frac{1}{T\log T}\int_T^{2T}\sum_\gamma f\Big(\frac{\log T}{2\pi}(\gamma-t)\Big)\boldsymbol\zeta_{\alpha,\beta}(t)dt
\end{equation}
with $\gamma\in\mathbb R$ since we are assuming RH (we recall that $\rho=1/2+i\gamma$ are the non-trivial zeros of $\zeta$). By the residue theorem we have that
\begin{equation}\begin{split}\label{aux4.1}
\langle N_f\rangle_{|\zeta|^2}^{\alpha,\beta}=\frac{1}{T\log T}\int_T^{2T}&\frac{1}{2\pi i}\bigg(\int_{(c)}-\int_{(1-c)}\bigg) \frac{\zeta'}{\zeta}(s+it)\\&
\cdot f\Big(\frac{-i\log T}{2\pi}(s-1/2)\Big)ds\;\boldsymbol\zeta_{\alpha,\beta}(t)dt
\end{split}\end{equation}
where $c\in(\frac{1}{2},1)$ and $\int_{(c)}$ denotes the integral over the vertical line of those $s$ such that $\Re(s)=c$. We select $c=\frac{1}{2}+\delta$ with $\delta\asymp(\log T)^{-1}$ and we first consider the integral over the $c$-line
\begin{equation}\begin{split}\notag
\mathcal I:&
=\frac{1}{T\log T}\int_T^{2T}\frac{1}{2\pi i}\int_{(c)}\frac{\zeta'}{\zeta}(s+it)f\Big(\frac{-i\log T}{2\pi}\Big(s-\frac{1}{2}\Big)\Big)ds\; \boldsymbol\zeta_{\alpha,\beta}(t)dt\\&
=\frac{1}{2\pi}\int_{-\infty}^{+\infty}f\Big(\frac{\log T}{2\pi}(y-i\delta)\Big)\frac{d}{d\gamma}\bigg[\frac{I(\alpha;\beta;\delta+iy+\gamma;\delta+iy)}{T\log T} \bigg]_{\gamma=0} dy
\end{split}\end{equation}
where 
\begin{equation}\label{aux4.01}
I(A;B;C;D):=\int_T^{2T}\frac{\zeta(1/2+A+it)\zeta(1/2+B-it)\zeta(1/2+C+it)}{\zeta(1/2+D+it)}dt.
\end{equation}
Moments like~\eqref{aux4.01} can be computed thanks to Conjecture~\ref{ratioconj} and it turns out to be
\begin{equation}\begin{split}\label{aux4.3}
I(A;B;C;D)&=
\int_T^{2T}\bigg\{\frac{\zeta(1+A+B)\zeta(1+B+C)}{\zeta(1+B+D)}\\&
\hspace{1cm}+\Big(\frac{t}{2\pi}\Big)^{-A-B}\frac{\zeta(1-A-B)\zeta(1-A+C)}{\zeta(1-A+D)}\\& 
\hspace{1cm}+\Big(\frac{t}{2\pi}\Big)^{-B-C}\frac{\zeta(1+A-C)\zeta(1-B-C)}{\zeta(1-C+D)}\bigg\}dt\\&
\hspace{.5cm}+O(T^{1/2+\varepsilon})
\end{split}\end{equation}
for suitable shifts $A,B,C,D$, i.e. with real part $\asymp(\log T)^{-1}$ and imaginary part $\ll_\varepsilon T^{1-\varepsilon}$, for every $\varepsilon>0$ (see e.g. \cite[Section 2.1]{CSapplications}). Notice that the arithmetical factor $A_{\zeta}(\alpha;\beta;\gamma;\delta)$ from Conjecture~\ref{ratioconj} equals $1$ in our case, with $K=2$, $L=1$, $Q=1$, $R=0$ (this can be easily proven by direct computation or deduced by \cite[Corollary 2.6.2]{CFKRSrecipe}).
%\begin{equation}\begin{split}\notag
%&A_\zeta(A,B,C,D)\\&
%=\prod_p\frac{(1-\frac{1}{p^{1+A+B}})(1-\frac{1}{p^{1+C+D}})}{(1-\frac{1}{p^{1+B+D}})}\sum_{a+c+d=b}\frac{\mu(p^d)}{p^{(1/2+A)a+(1/2+C)c+(1/2+B)b+(1/2+D)d}}\\&
%=\prod_p\frac{(1-\frac{1}{p^{1+A+B}})(1-\frac{1}{p^{1+C+D}})}{(1-\frac{1}{p^{1+B+D}})}\sum_{a,c,d}\frac{\mu(p^d)}{p^{(1+A+B)a+(1+C+B)c+(1+D+B)d}}=1. 
%\end{split}\end{equation}
We now want to apply~\eqref{aux4.3} with $A=\alpha$, $B=\beta$, $C=\delta+iy+\gamma$, $D=\delta+iy$ and to do so we need that the imaginary parts of all the shifts are $\ll_\varepsilon T^{1-\varepsilon}$. A standard technique to avoid this issue is splitting the integral over $y$ in two pieces; the contribution to $\mathcal I$ coming from $|y|>T^{1-\varepsilon}$ is $\ll T^{-1+\varepsilon}$, thanks to the good decaying of $f$ and to RH, since 
\begin{equation}\begin{split}\notag
\frac{1}{T\log T}\int_T^{2T}|\boldsymbol\zeta_{\alpha,\beta}(t)&|\int_{|y|>T^{1-\varepsilon}}|f(y\log T)|\bigg|\frac{\zeta'}{\zeta}(1/2+\delta+iy+it)\bigg|dydt\\&\hspace{1cm}
\ll \frac{T^{\varepsilon/100}}{T}\int_T^{2T}\int_{|y|>T^{1-\varepsilon}}\frac{\log(|y|t)}{|y|^2}dydt\ll T^{-1+\varepsilon}.
\end{split}\end{equation}
%%%%%%%%%%%%%%%%%%%%%%%%%%%%%%%%%%%%%%%%%%%%%%%%%%%%%%%%%
%qui rifaccio le cose con l'asintotica per zeta(1+z), che non va bene nel caso di f di schwartz, ma serve il supporto compatto
%Therefore we can truncate the integral over $y$ at height $T^{1-\varepsilon}$, apply (\ref{aux4.3}) and then re-extend the integration over $y$ to infinity with a small error term. Thus, by the asymptotic expansion $\zeta(1+s)=\frac{1}{s}+O(1)$ at $s=0$, differentiating with respect to $\gamma$ at $\gamma=0$, setting $\delta=0$ (we are allowed to so so since now the integral is regular at $\delta=0$) we get
%\begin{equation}\begin{split}\label{aux4.4}
%I_c^f=
%\Big(1&+O\Big(\frac{1}{\log T}\Big)\Big)\frac{1}{2\pi}\int_{-\infty}^{+\infty}f\Big(\frac{\log T}{2\pi}y\Big)\frac{1}{T\log T}\int_T^{2T}f_{\alpha,\beta}(y;t)dt\;dy%+O(T^{1/2+\varepsilon}).
%\end{split}\end{equation}
%with
%\begin{equation}\begin{split}\label{aux4.4aggiunta}
%f_{\alpha,\beta}(y;t):=-\frac{1}{(\alpha+\beta)(\beta+iy)}
%&-\Big(\frac{t}{2\pi}\Big)^{-\alpha-\beta}\frac{1}{(-\alpha-\beta)(-\alpha+iy)}\\&
%-\Big(\frac{t}{2\pi}\Big)^{-\beta-iy} \frac{1}{(\alpha-iy)(-\beta-iy)}.
%\end{split}\end{equation}
%%%%%%%%%%%%%%%%%%%%%%%%%%%%%%%%%%%%%%%%%%%%%%%%%%%%%%%%%
%qui invece faccio preciso
Therefore we can truncate the integral over $y$ at height $T^{1-\varepsilon}$, apply~\eqref{aux4.3} and then re-extend the integration over $y$ to infinity with a small error term. Thus, differentiating with respect to $\gamma$ at $\gamma=0$, moving the path of integration to $\delta=0$ (we are allowed to so so since now the integral is regular at $\delta=0$) we get
\begin{equation}\label{aux4.4}
\mathcal I=
\frac{1}{2\pi}\int_{-\infty}^{+\infty}f\Big(\frac{\log T}{2\pi}y\Big)\frac{1}{T\log T}\int_T^{2T}g_{\alpha,\beta}(y;t)dtdy+O(T^{-1/2+\varepsilon})
\end{equation}
with
\begin{equation}\begin{split}\label{aux4.4aggiunta}
g_{\alpha,\beta}(y;t):=&\frac{\zeta(1+\alpha+\beta)\zeta'(1+\beta+iy)}{\zeta(1+\beta+iy)}
\\&\hspace{.5cm}+\Big(\frac{t}{2\pi}\Big)^{-\alpha-\beta}\frac{\zeta(1-\alpha-\beta)\zeta'(1-\alpha+iy)}{\zeta(1-\alpha+iy)}\\&\hspace{.5cm}-\Big(\frac{t}{2\pi}\Big)^{-\beta-iy}\zeta(1+\alpha-iy)\zeta(1-\beta-iy).
\end{split}\end{equation}
We notice that, when computing this derivative, it is useful to observe that if $f(z)$ is analytic at $z = 0$, then (see \cite[Equation (2.13)]{CSapplications}) 
\begin{equation}\notag
\frac{d}{d\gamma}\bigg[\frac{f(\gamma)}{\zeta(1-\gamma)}\bigg]_{\gamma=0}=-f(0).
\end{equation}
%%%%%%%%%%%%%%%%%%%%%%%%%%%%%%%%%%%%%%%%%%%%%%%%%%%%%%%%%
Similarly we deal with the integral over the $(1-c)$-line in~\eqref{aux4.1}
\begin{equation}\begin{split}\notag
\mathcal J:&
=\frac{1}{T\log T}\int_T^{2T}\frac{1}{2\pi i}\int_{(1-c)}\frac{\zeta'}{\zeta}(s+it)f\Big(\frac{-i\log T}{2\pi}\Big(s-\frac{1}{2}\Big)\Big)ds\;\boldsymbol\zeta_{\alpha,\beta}(t)dt\\&
=\frac{1}{2\pi i}\int_{(c)}f\Big(\frac{-i\log T}{2\pi}\Big(s-\frac{1}{2}\Big)\Big)\frac{1}{T\log T}\int_T^{2T}\frac{\zeta'}{\zeta}(1-s+it)\boldsymbol\zeta_{\alpha,\beta}(t)dtds.
\end{split}\end{equation}
Using the functional equation 
$$ \frac{\zeta'}{\zeta}(1-z) =\frac{X'}{X}(z)-\frac{\zeta'}{\zeta}(z) $$
where
$$ \frac{X'}{X}(z):=\log \pi - \frac{1}{2}\frac{\Gamma'}{\Gamma}\Big(\frac{z}{2}\Big) - \frac{1}{2}\frac{\Gamma'}{\Gamma}\Big(\frac{1-z}{2}\Big) $$
we express $\mathcal J$ as a sum of two terms
\begin{equation}\label{aux4.a.1}
\mathcal J=\mathcal J_1-\mathcal J_2
\end{equation}
with
\begin{equation}\notag
\mathcal J_1:=\frac{1}{2\pi i}\int_{(c)}f\Big(\frac{-i\log T}{2\pi}(s-1/2)\Big)\frac{1}{T\log T}\int_T^{2T}\frac{X'}{X}(s-it)\boldsymbol\zeta_{\alpha,\beta}(t)dtds
\end{equation}
and
\begin{equation}\notag
\mathcal J_2:=\frac{1}{2\pi i}\int_{(c)}f\Big(\frac{-i\log T}{2\pi}(s-1/2)\Big)\frac{1}{T\log T}\int_T^{2T}\frac{\zeta'}{\zeta}(s-it)\boldsymbol\zeta_{\alpha,\beta}(t)dt ds.
\end{equation}
With $c=1/2+\delta$, $\delta\to 0$, it is easy to see that
\begin{equation}\begin{split}\label{aux4.a.4}
\mathcal J_1&
=\frac{1}{2\pi}\int_{-\infty}^{+\infty}f\Big(\frac{\log T}{2\pi}y\Big)\frac{1}{T\log T}\int_T^{2T}\frac{X'}{X}(\tfrac{1}{2}+iy-it)\boldsymbol\zeta_{\alpha,\beta}(t)dtdy\\&
=-\frac{1}{2\pi }\int_{-\infty}^{+\infty}f\Big(\frac{\log T}{2\pi}y\Big)\frac{\log T+O(1)}{T\log T}\int_T^{2T}\boldsymbol\zeta_{\alpha,\beta}(t)dt dy
\end{split}\end{equation}
since, using Stirling's approximation to estimate the gamma-factors, we have (again we can assume $y\ll T^{1-\varepsilon}$ because of the great decaying of $f$)
\begin{equation}\begin{split}\notag
\frac{X'}{X}(1/2+iy-it)&=-\frac{1}{2}\frac{\Gamma'}{\Gamma}\bigg(\frac{1}{4}+\frac{iy}{2}-\frac{it}{2}\bigg)-\frac{1}{2}\frac{\Gamma'}{\Gamma}\bigg(\frac{1}{4}-\frac{iy}{2}+\frac{it}{2}\bigg)+O(1)\\&
=-\frac{1}{2}\log\Big(-\frac{it}{2}\Big)-\frac{1}{2}\log\Big(\frac{it}{2}\Big)+O(1)=-\log T+O(1).
\end{split}\end{equation}
Moreover, with the same choice of $c$ as before, if we set $\alpha=\beta$ we get
\begin{equation}\label{aux4.a.5}
\mathcal J_2=\mathcal I.
\end{equation}
Then~\eqref{aux4.1},~\eqref{aux4.a.1} and~\eqref{aux4.a.5} imply that 
\begin{equation}\label{aux4.10} 
\langle N_f \rangle_{|\zeta|^2}^{\alpha,\alpha}=\mathcal I+\mathcal J=-\mathcal J_1+2\mathcal I 
\end{equation}
and the function $\mathcal J_1=\mathcal J_1(\alpha)$ is regular at $\alpha=0$, then we can take the limit in~\eqref{aux4.a.4}, getting
\begin{equation}\begin{split}\label{aux4.11} 
\lim_{\alpha\to 0}\mathcal J_1&
=-\frac{1}{2\pi}\int_{-\infty}^{+\infty}f\Big(\frac{\log T}{2\pi}y\Big)\frac{\log T+O(1)}{T\log T}\int_T^{2T}|\zeta(1/2+it)|^2dt\; dy\\&
=-\frac{\log T+O(1)}{2\pi}\int_{-\infty}^{+\infty}f\Big(\frac{\log T}{2\pi}y\Big) dy
=-\int_{-\infty}^{+\infty}f(x)dx+O\Big(\frac{1}{\log T}\Big)
\end{split}\end{equation}
with the change of variable $\frac{\log T}{2\pi}y=x$. Lastly, we study the remaining term $\mathcal I$, from~\eqref{aux4.4} and~\eqref{aux4.4aggiunta}. We set $\alpha=\beta=a/\log T$ with $0<a<1$, we perform the same change of variable $\frac{\log T}{2\pi}y=x$ as before and we get
%%%%%%%%%%%%%%%%%%%%%%%%%%%%%%%%%%%%%%%%%%%%%%%%%%%%%%%%%
%vecchio modo, con le asintotiche
%\begin{equation}\begin{split}\notag
%I_c^f&
%=\Big(1+O\Big(\frac{1}{\log T}\Big)\Big)\int_{-\infty}^{+\infty}f(x)\frac{1}{T(\log T)^2}\int_T^{2T}g_{\frac{a}{\log T},\frac{a}{\log T}}\Big(\frac{2\pi x}{\log T};t\Big)dtdx\\&
%=\Big(1+O\Big(\frac{1}{\log T}\Big)\Big)\int_{-\infty}^{+\infty}f(x)\frac{1}{T}\int_T^{2T}\bigg\{
%-\frac{1}{2a(a+2\pi ix)}-\Big(\frac{t}{2\pi}\Big)^{-\frac{2a}{\log T}}\\&\hspace{1cm}
%\cdot\frac{1}{-2a(-a+2\pi i x)}
%-\Big(\frac{t}{2\pi}\Big)^{-\frac{a+2\pi ix}{\log T}} \frac{1}{(a-2\pi ix)(-a-2\pi ix)}
%\bigg\}dt\;dy
%\end{split}\end{equation}
%and since $\log\frac{t}{2\pi}=\log T+O(1)$ as $t\in[T,2T]$, then
%\begin{equation}\notag
%I_c^f=\Big(1+O\Big(\frac{1}{\log T}\Big)\Big)\int_{-\infty}^{+\infty}f(x)\mathcal P(a,x)dy 
%\end{equation}
%%%%%%%%%%%%%%%%%%%%%%%%%%%%%%%%%%%%%%%%%%%%%%%%%%%%%%%%%
%
%%%%%%%%%%%%%%%%%%%%%%%%%%%%%%%%%%%%%%%%%%%%%%%%%%%%%%%%%
%nuovo modo, con le asintotiche
\begin{equation}\notag
\mathcal I
=%\Big(1+O\Big(\frac{1}{\log T}\Big)\Big)
\int_{-\infty}^{+\infty}f(x)\frac{1}{T(\log T)^2}\int_T^{2T}g_{\frac{a}{\log T},\frac{a}{\log T}}\Big(\frac{2\pi x}{\log T};t\Big)dtdx
+O(T^{1/2+\varepsilon})
\end{equation}
and since $\log\frac{t}{2\pi}=\log T+O(1)$ as $t\in[T,2T]$, then
\begin{equation}\begin{split}\label{aggiunta**}
\mathcal I=\bigg(\frac{1}{(\log T)^2}+O\Big(\frac{1}{(\log T)^3}\Big)\bigg)&\int_{-\infty}^{+\infty}f(x) 
\bigg(\zeta(1+\tfrac{2a}{\log T})\frac{\zeta'}{\zeta}(1+\tfrac{a+2\pi ix}{\log T})
\\&+e^{-2a}\zeta(1-\tfrac{2a}{\log T})\frac{\zeta'}{\zeta}(1-\tfrac{a-2\pi ix}{\log T})
\\&-e^{-a-2\pi ix}\zeta(1+\tfrac{a-2\pi ix}{\log T})\zeta(1-\tfrac{a+2\pi ix}{\log T})
\bigg)dx 
\end{split}\end{equation} 
where the error term is uniform in $a$.
%%%%%%%%%%%%%%%%%%%%%%%%%%%%%%%%%%%%%%%%%%%%%%%%%%%%%%%%%%%%%%%%%%%%%%%%%%%%%%%%%%%%%%%%%%%%%%%%%%%%%%%%%%%%%%%%%%%%%%%%%%%%%%%%%%%%%%%%%%%%%%%%%%%%%%%%%%%%%%%%%%%%%%%%%%%%%%%%%%%%%%%%%%%%%%%%%%%%%%%%
%For $T\to\infty$, $\mathcal I=\mathcal I(a)$ can then be evaluated as $a\to0$, getting:
Now, we will prove that the above expression is regular at $a=0$, showing that
\begin{equation}\label{cucuzzona}
\lim_{a\to0}\mathcal I=\int_{-\infty}^{+\infty}f(x) \mathcal P(x)dx +O\Big(\frac{1}{\log T}\Big)
\end{equation}
as $T\to\infty$, where
\begin{equation}\notag
\mathcal P(x):=\frac{-1+2\pi ix+e^{-2\pi ix}}{4\pi^2x^2}.
\end{equation}
Intuitively, if we replace each zeta function with its leading term in the expansion at the point 1 given by $\zeta(1+z)\sim\frac{1}{z}$, we have
\begin{equation}\begin{split}\notag
\mathcal I&\approx\frac{1}{(\log T)^2}\int_{-\infty}^{+\infty}f(x) 
\bigg(\frac{-(\log T)^2}{2a(a+2\pi ix)}+e^{-2a}\frac{-(\log T)^2}{2a(a-2\pi ix)}
\\&\hspace{6.65cm}+e^{-a-2\pi ix}\frac{(\log T)^2}{(a-2\pi ix)(a+2\pi ix)}
\bigg)dx 
\\&=\int_{-\infty}^{+\infty}f(x) 
\bigg(-\frac{1}{2a(a+2\pi ix)}-\frac{e^{-2a}}{2a(a-2\pi ix)}+\frac{e^{-a-2\pi ix}}{(a-2\pi ix)(a+2\pi ix)}
\bigg)dx 
\end{split}\end{equation} 
and the function inside the parentheses above equals
\begin{equation}\notag
\frac{-a(1+e^{-2a})+2\pi ix(1-e^{-2a})+2ae^{-a-2\pi ix}}{2a(a^2+4\pi^2x^2)}
=\frac{-1+2\pi ix+e^{-2\pi ix}+O(a)}{4\pi^2x^2+O(a^2)}
\end{equation}
and then tends to $\mathcal P(x)$ as $a\to0$.\\
To show~\eqref{cucuzzona} rigorously, we split the integral over $x$ into two parts. We start with the case $x\ll \log T$; from Taylor approximation $f(1+s\pm y)= f(1+s)\pm yf'(1+s)+O_s(y^2)$ we get
\begin{equation}\begin{split}\notag 
\frac{\zeta'}{\zeta}\Big(1+\frac{2\pi ix}{\log T}\pm \frac{a}{\log T}\Big)&=  
\frac{\zeta'}{\zeta}\Big(1+\frac{2\pi ix}{\log T}\Big)\pm \frac{a}{\log T}\Big(\frac{\zeta'}{\zeta}\Big)'\Big(1+\frac{2\pi ix}{\log T}\Big)+O_T(a^2)
\\&=: c_1(x)\pm\frac{a}{\log T}c_2(x)+O_T(a^2)
\end{split}\end{equation}
and
\begin{equation}\notag \zeta\Big(1-\frac{2\pi ix}{\log T}\pm \frac{a}{\log T}\Big)=  
\zeta\Big(1-\frac{2\pi ix}{\log T}\Big)+O_T(a)=:k(x)+O_T(a)\end{equation}
as $a\to0$, with the notations $c_1(x)=c_1(x,T):=\frac{\zeta'}{\zeta}(1+\frac{2\pi ix}{\log T})$, $c_2(x)=c_1(x,T):=(\frac{\zeta'}{\zeta})'(1+\frac{2\pi ix}{\log T})$ and $k(x)=k(x,T):=\zeta(1-\frac{2\pi ix}{\log T})$.
Moreover we use the asymptotic expansion 
\begin{equation}\label{zetaat1}\zeta(1+z)=\frac{1}{z}+\gamma+O(z) \quad\quad z\to0,\end{equation}
and we get
\begin{equation}\begin{split}\notag
&\int_{x\ll \log T}f(x) 
\bigg(\zeta(1+\tfrac{2a}{\log T})\frac{\zeta'}{\zeta}(1+\tfrac{a+2\pi ix}{\log T})
\\&\hspace{3cm}+e^{-2a}\zeta(1-\tfrac{2a}{\log T})\frac{\zeta'}{\zeta}(1-\tfrac{a-2\pi ix}{\log T})
\\&\hspace{3cm}-e^{-a-2\pi ix}\zeta(1+\tfrac{a-2\pi ix}{\log T})\zeta(1-\tfrac{a+2\pi ix}{\log T})
\bigg)dx
\\&=\int_{x\ll \log T}f(x) 
\bigg(\Big[\frac{\log T}{2a}+\gamma+O_T(a)\Big]
\Big[c_1(x)+\frac{a}{\log T}c_2(x)+O_T(a^2)\Big]
\\&\hspace{1.8cm}+e^{-2a}
\Big[\frac{-\log T}{2a}+\gamma+O_T(a)\Big]
\Big[c_1(x)-\frac{a}{\log T}c_2(x)+O_T(a^2)\Big]
\\&\hspace{1.8cm}-e^{-a-2\pi ix}
\Big[k(x)+O_T(a)\Big]^2
\bigg)dx
\end{split}\end{equation} 
whose limit as $a\to0$ is
%\begin{equation}\begin{split}\notag
%&\int_{x\ll \log T}f(x)\bigg\{\lim_{a\to0} 
%\bigg(\Big[\frac{c_1(x)\log T}{2a}+\frac{c_2(x)}{2}+\gamma c_1(x)\Big]
%\\&\hspace{.8cm}+e^{-2a}
%\Big[\frac{-c_1(x)\log T}{2a}+\frac{c_2(x)}{2}+\gamma c_1(x)\Big]\bigg)
%-e^{-2\pi ix}
%k(x)^2\bigg\}
%dx
%\\&=\int_{x\ll \log T}f(x)\bigg\{\lim_{a\to0} 
%\bigg(c_1(x)\log T\Big(\frac{1-e^{-2a}}{2a}\Big)+c_2(x)\Big(\frac{1+e^{-2a}}{2}\Big)
%\\&\hspace{4.9cm}+\gamma c_1(x)(1+e^{-2a})\bigg)
%-e^{-2\pi ix}
%k(x)^2\bigg\}
%dx
%\\&=\int_{x\ll \log T}f(x)\bigg\{c_1(x)\log T+c_2(x)+
%2\gamma c_1(x)-e^{-2\pi ix}k(x)^2\bigg\}dx.
%\end{split}\end{equation}
$$ \int_{x\ll \log T}f(x)\bigg\{c_1(x)\log T+c_2(x)+
2\gamma c_1(x)-e^{-2\pi ix}k(x)^2\bigg\}dx. $$
By definition of $c_1(x),c_2(x), k(x)$, the asymptotic expansion~\eqref{zetaat1} yields $c_1(x)=-\frac{\log T}{2\pi ix}+O(1)$, $c_2(x)=\frac{(\log T)^2}{(2\pi ix)^2}+O(1)$ and $k(x)=\frac{(\log T)^2}{(2\pi ix)^2}-\frac{2\gamma\log T}{2\pi ix}+O(1)$, uniformly for $x\ll\log T$. Then the above is
\begin{equation}\notag
=\int_{x\ll \log T}f(x)\bigg\{
-\frac{(\log T)^2}{2\pi ix}+\frac{(\log T)^2}{(2\pi ix)^2}-e^{-2\pi ix}\frac{(\log T)^2}{(2\pi ix)^2}+O(\log T)
\bigg\}dx 
\end{equation}
(note that the sum $2\gamma c_1(x)-e^{-2\pi ix}k(x)^2$ gives the third term in the parentheses with an error $O(\log T)$, a possible pole at $x=0$ cancels out), which is
\begin{equation}\notag
=(\log T)^2\int_{x\ll \log T}f(x)\frac{-1+2\pi ix+e^{-2\pi ix}}{4\pi^2x^2}
dx+O(\log T).
\end{equation}  
Finally we can re-extend the range of integration with a small error term (being $f(x)\ll1/(1+x^2)$ and $\mathcal P(x)$ bounded), getting that the contribution of $x\ll\log T$ in the integral over $x$ in~\eqref{aggiunta**}, in the limit as $a\to0$, equals
\begin{equation}\notag
=(\log T)^2\int_{-\infty}^{+\infty}f(x)\mathcal P(x)dx+O(\log T).
\end{equation}
To prove~\eqref{cucuzzona}, we finally have to bound the contribution of $x\gg\log T$ in the integral on the right hand side~\eqref{aggiunta**}, as $a\to0$; to do so, we use the bounds $\zeta(1+iy)\ll\log y$ (see \cite[Theorem 3.5]{Tit}) and $\frac{\zeta'}{\zeta}(1+iy)\ll\log y$ (see \cite[Equation (3.11.9)]{Tit}) for $y\gg 1$, thus the contribution coming from $x\gg\log T$ is
\begin{equation}\begin{split}\notag
&=\lim_{a\to0}\int_{x\gg\log T}f(x)\bigg(\Big[\frac{\log T}{2a}+O(1)\Big]\Big[\frac{\zeta'}{\zeta}(1+\tfrac{2\pi ix}{\log T})+O\Big(\frac{a}{\log T}\Big)\Big]\\&
\;\;\quad+e^{-2a}\Big[-\frac{\log T}{2a}+O(1)\Big]\Big[\frac{\zeta'}{\zeta}(1+\tfrac{2\pi ix}{\log T})+O\Big(\frac{a}{\log T}\Big)\Big]+O\big((\log x)^2\big)\bigg)dx
\\& =\int_{x\gg\log T}f(x)\bigg(\log T\frac{\zeta'}{\zeta}(1+\tfrac{2\pi ix}{\log T})\lim_{a\to0}\Big[\frac{1-e^{-2a}}{2a}\Big]+O\big((\log x)^2\big)\bigg)dx
\\&\ll \int_{x\gg\log T}|f(x)|\Big(\log T\log x+(\log x)^2\Big)dx 
\ll \log T\int_{x\gg\log T}\frac{\log x}{x^2}dx,
\end{split}\end{equation}
then~\eqref{cucuzzona} follows, being $\int_{x\gg\log T}\frac{\log x}{x^2}dx\ll1$. Finally, if we decompose $\mathcal P(x)$ in even and odd parts
\begin{equation}\label{limitatezzaP}
\mathcal P(x)=-\frac{1}{2}\frac{\sin^2(\pi x)}{(\pi x)^2}-\frac{i\big(\sin(2\pi x)-2\pi x\big)}{4\pi^2x^2}
\end{equation}
since $f$ is even and $\mathcal P(x)$ bounded, we have
\begin{equation}\label{aux4.b}
\lim_{a\to0}\mathcal I= -\frac{1}{2}\int_{-\infty}^{+\infty}f(x)\frac{\sin^2(\pi x)}{(\pi x)^2}dy+O\Big(\frac{1}{\log T}\Big).
\end{equation}
Putting together~\eqref{aux4.10},~\eqref{aux4.11} and~\eqref{aux4.b} we finally get
\begin{equation}\notag
\langle N_f\rangle_{|\zeta|^2}= \int_{-\infty}^{+\infty}f(x)\bigg(1-\frac{\sin^2(\pi x)}{(\pi x)^2}\bigg)dx+O\Big(\frac{1}{\log T}\Big)
\end{equation}
as $T\to\infty$ and the theorem has been proved.

\subsection{Proof of Proposition~\ref{casozeta^4}}\label{C4S2S4}
This proof builds on the same ideas as that of Proposition~\ref{congettura}, even though we have to handle longer computations; 
to begin with, we introduce four parameters $\alpha,\beta,\nu,\eta\in\mathbb R$ of size $1/\log T$, we denote
\begin{equation}\notag
\boldsymbol\zeta_{\alpha,\beta,\nu,\eta}(t):=\zeta(1/2+\alpha+it)\zeta(1/2+\beta+it)\zeta(1/2+\nu-it)\zeta(1/2+\eta-it)
\end{equation}
and we look at 
\begin{equation}\label{aux4.3.3.0}
\langle N_f \rangle_{|\zeta|^4}^{\alpha,\beta,\nu,\eta}:=\frac{1}{\frac{1}{2\pi^2}T(\log T)^4}\int_T^{2T}\sum_\gamma f\Big(\frac{\log T}{2\pi}(\gamma-t)\Big)\boldsymbol\zeta_{\alpha,\beta,\nu,\eta}(t)dt
\end{equation}
with $\gamma\in\mathbb R$ since we are assuming RH. In analogy to Equation~\eqref{aux4.10}, the residue theorem yields
\begin{equation}\label{aux4.3.3.1}
\langle N_f \rangle_{|\zeta|^4}^{\alpha,\beta,\alpha,\beta}=-\mathcal J_1+2\mathcal I
\end{equation}
with
\begin{equation}\begin{split}\label{aux4.3.3.2}
\mathcal J_1=\mathcal J_1(\alpha,\beta)&=-\int_{-\infty}^{+\infty}f(x)\frac{\log T+O(1)}{\frac{1}{2\pi^2}T(\log T)^5}\int_T^{2T}\boldsymbol\zeta_{\alpha,\beta,\alpha,\beta}(t)dtdx
\\&\stackrel{\alpha,\beta\to0}{\longrightarrow}
-\int_{-\infty}^{+\infty}f(x)dx+O\Big(\frac{1}{\log T}\Big)
\end{split}\end{equation}
and  
\begin{equation}\begin{split}\notag
\mathcal I
&=\frac{1}{2\pi}\int_{-\infty}^{+\infty}f\Big(\frac{\log T}{2\pi}(y-i\delta)\Big)\frac{d}{d\gamma}\bigg[\frac{I(\alpha;\beta;\gamma;\delta+iy;\alpha;\beta)}{\frac{1}{2\pi^2}T(\log T)^4} \bigg]_{\gamma=\delta+iy} dy
\\&=\int_{-\infty}^{+\infty}f\Big(x-\frac{i\delta\log T}{2\pi}\Big)\frac{d}{d\gamma}\bigg[\frac{I(\frac{a}{\log T};\frac{b}{\log T};\gamma;\delta+\frac{2\pi ix}{\log T};\frac{a}{\log T};\frac{b}{\log T})}{\frac{1}{2\pi^2}T(\log T)^5} \bigg]_{\gamma=\delta+\frac{2\pi ix}{\log T}} dx
\end{split}\end{equation}
where $a,b\asymp 1$, $\delta\asymp 1/\log T$ and $I(A;B;C;D;F;G)$ is defined by
\begin{equation}\notag
\int_T^{2T}\frac{\zeta(\tfrac{1}{2}+A+it)\zeta(\tfrac{1}{2}+B+it)\zeta(\tfrac{1}{2}+C+it)\zeta(\tfrac{1}{2}+F-it)\zeta(\tfrac{1}{2}+G-it)}{\zeta(\tfrac{1}{2}+D+it)}dt.
\end{equation}
If the shifts satisfy the conditions prescribed by Conjecture~\ref{ratioconj} then such an integral can be evaluated by using the ratio conjecture. According to the recipe, up to an error $O(T^{1/2+\varepsilon})$, the above moment is a sum of ten pieces, the first being
\begin{equation}\notag
\int_T^{2T} \tfrac{\zeta(1+A+F)\zeta(1+A+G)\zeta(1+B+F)\zeta(1+B+G)\zeta(1+C+F)\zeta(1+C+G)}{\zeta(1+D+F)\zeta(1+D+G)}\mathcal A_{A,B,C,D,F,G} \;dt
\end{equation}
where
\begin{equation}\begin{split}\notag
\mathcal A_{A,B,C,D,F,G} 
=&\prod_p (1-\frac{1}{p^{1+A+F}})(1-\frac{1}{p^{1+A+G}})(1-\frac{1}{p^{1+B+F}})(1-\frac{1}{p^{1+B+G}})
\\&(1-\frac{1}{p^{1+C+F}})(1-\frac{1}{p^{1+C+G}})
(1-\frac{1}{p^{1+D+F}})^{-1}(1-\frac{1}{p^{1+D+G}})^{-1}
\\&\sum_{a+b+c+d=f+g}\frac{\mu(p^d)}{p^{(\frac{1}{2}+A)a+(\frac{1}{2}+B)b+(\frac{1}{2}+C)c+(\frac{1}{2}+D)d+(\frac{1}{2}+F)f+(\frac{1}{2}+G)g}}.
\end{split}\end{equation}
It will be useful to notice that if all the shifts equal zero, then 
\begin{equation}\notag
\mathcal A:=\mathcal A_{0,0,0,0,0,0}=\frac{1}{\zeta(2)};
\end{equation}
again this can be proven by direct computation or deduced by \cite[Corollary 2.6.2]{CFKRSrecipe}.
%indeed, since $\mu(p^d)$ forces $d$ to be either $0$ or $1$, we have that
%\begin{equation}\begin{split}\notag
%\mathcal A
%&=\prod_p \Big(1-\frac{1}{p}\Big)^4\bigg(\sum_{a+b+c=f+g}\frac{1}{p^{\frac{a+b+c+f+g}{2}}}+\sum_{a+b+c+1=f+g}\frac{-1}{p^{\frac{a+b+c+1+f+g}{2}}}\bigg)
%\\&=\prod_p \Big(1-\frac{1}{p}\Big)^4
%\sum_{n=0}^{\infty}\frac{1}{p^n}\bigg[\sum_{a+b+c=n}1\cdot\sum_{f+g=n}1
%-\sum_{a+b+c+1=n}1\cdot\sum_{f+g=n}1\bigg].
%\end{split}\end{equation}
%We recall that the number of representations of an integer $n$ as a sum of $k$ nonnegative integers is $\binom{n+k-1}{k-1}$ thus the sum over $n$ above equals
%$$\sum_{n=0}^{\infty}\frac{1}{p^n}\bigg[\frac{(n+2)(n+1)^2}{2}-\frac{n(n+1)^2}{2}\bigg]
%=\sum_{n=0}^{\infty}\frac{(n+1)^2}{p^n}.$$
%Differentiating the closed formula for the geometric series, one easily gets that $\sum_n (n+1)^2x^n=\frac{1+x}{(1-x^3)}$ for $|x|<1$, therefore
%$$ \mathcal A= \prod_p \Big(1-\frac{1}{p}\Big)^4\frac{1+\frac{1}{p}}{(1-\frac{1}{p})^3}=\prod_p \bigg(1-\frac{1}{p^2}\bigg)=\frac{1}{\zeta(2)}$$
%and~\eqref{aux4.3.3.3} is proven.\par
All the other nine terms can be recovered from the first one just by swapping the shifts as prescribed by the recipe; doing so yields a formula for $I(A;B;C;D;F;G)$ and differentiating with respect to $C$ at $C=D$ we get
\begin{equation}\begin{split}\label{aux4.3.3.4}
&\frac{d}{dC}\big[I(A;B;C;D;F;G)\big]_{C=D}
\\&=\int_T^{2T}\bigg(R_1+(\tfrac{t}{2\pi})^{-A-F}R_2+(\tfrac{t}{2\pi})^{-A-G}R_3+(\tfrac{t}{2\pi})^{-B-F}R_4+(\tfrac{t}{2\pi})^{-B-G}R_5
\\&\hspace{2cm}+(\tfrac{t}{2\pi})^{-D-F}R_6+(\tfrac{t}{2\pi})^{-D-G}R_7+(\tfrac{t}{2\pi})^{-A-B-F-G}R_8
\\&\hspace{2cm}+(\tfrac{t}{2\pi})^{-A-D-F-G}R_9+(\tfrac{t}{2\pi})^{-B-D-F-G}R_{10}\bigg)dt+O(T^{1/2+\varepsilon})
\end{split}\end{equation}
with
\begin{equation}\begin{split}\notag
&R_1=R_1(A,B,D,F,G)
=\tfrac{\mathcal A_{A,B,D,D,F,G}}{\zeta(1+A+F)\zeta(1+A+G)\zeta(1+B+F)\zeta(1+B+G)\zeta(1+D+F)\zeta(1+D+G)}
\\&\hspace{4.6cm}\cdot\Big[\tfrac{\zeta'}{\zeta}(1+D+F)+\tfrac{\zeta'}{\zeta}(1+D+G)+\frac{\mathcal A'_{A,B,D,D,F,G}}{A_{A,B,D,D,F,G}}\Big]
\\&R_2=R_1(-F,B,D,-A,G)
\\&R_3=R_1(-G,B,D,F,-A)
\\&R_4=R_1(A,-F,D,-B,G)
\\&R_5=R_1(A,-G,D,F,-B)
\\&R_6=R_6(A,B,D,F,G)\\&\hspace{.5cm}=-\tfrac{\zeta(1+A-D)\zeta(1+A+G)\zeta(1+B-D)\zeta(1+B+G)\zeta(1-F-D)\zeta(1-F+G)}{\zeta(1+D+G)}\mathcal A_{A,B,D,D,F,G}
\\&R_7=R_6(A,B,D,G,F)
\\&R_8=R_1(-F,-G,D,-A,-B)
\\&R_9=R_6(-F,B,D,G,-A)
\\&R_{10}=R_6(A,-F,D,G,-B).
\end{split}\end{equation}
If the shifts $A,B,D,F,G$ are $\ll 1/\log T$ the above formula simplifies a lot, since we have
\begin{equation}\begin{split}\notag
R_1&=\frac{(-2D-F-G)\mathcal A}{(A+F)(A+G)(B+F)(B+G)(D+F)(D+G)}+O\Big(\frac{(\log T)^5}{\log T}\Big)
%\\& =:f(A,B,D,F,G)+O\Big(\frac{(\log T)^5}{\log T}\Big)
\end{split}\end{equation}
and
\begin{equation}\begin{split}\notag
R_6&=\frac{-(D+G)\mathcal A}{(A-D)(A+G)(B-D)(B+G)(-F-D)(-F+G)}+O\Big(\frac{(\log T)^5}{\log T}\Big).
%\\& =:g(A,B,D,F,G)+O\Big(\frac{(\log T)^5}{\log T}\Big).
\end{split}\end{equation}
As in the proof of Proposition~\ref{congettura}, by a truncation of the integral over $x$ and Taylor approximations, we can use~\eqref{aux4.3.3.4} to evaluate $\mathcal I$; one can use Sage to carry out this massive computation, getting
\begin{equation}\begin{split}\label{aux4.3.3.6}
\lim_{\substack{a\to0 \\ b\to0}}\mathcal I
=\int_{-\infty}^{+\infty}f(x)\frac{(\log T)^5\mathcal A}{\frac{1}{2\pi^2}(\log T)^5}h(2\pi ix) dx+O\Big(\frac{1}{\log T}\Big)
\end{split}\end{equation}
with
\begin{equation}\begin{split}\notag
h(y):=-\frac{y^3-2y^2+6-e^{-y}(y^2+6y+6)}{6y^4}.
\end{split}\end{equation}
Note that, as in the last section, we moved the path of integration over $x$ to $\delta=0$, being the integral regular at $\delta=0$.
Therefore, putting together~\eqref{aux4.3.3.1},~\eqref{aux4.3.3.2} and~\eqref{aux4.3.3.6}, we get that
\begin{equation}\begin{split}\notag
\langle N_f\rangle_{|\zeta|^4}
&=\int_{-\infty}^{+\infty}f(x)\bigg(1+2\frac{2\pi^2}{\zeta(2)}h(2\pi ix)\bigg)dx+O\Big(\frac{1}{\log T}\Big)
\\&= \int_{-\infty}^{+\infty}f(x)\bigg(1+24h(2\pi ix)\bigg)dx+O\Big(\frac{1}{\log T}\Big)
\\&= \int_{-\infty}^{+\infty}f(x)W^{2}_U(x)dx+O\Big(\frac{1}{\log T}\Big)
\end{split}\end{equation}
since $f$ is even.

\subsection{Proof of Theorem \ref{analogoRMTunitario}}\label{RMTaggiuntaUnitary}
As mentioned in the introduction, the analogue of the ratio conjecture is a theorem in random matrix theory. Therefore, the same machinery described above proves Theorem \ref{analogoRMTunitario}, using \cite[Theorem 4.1]{CFZratioconj} in place of Conjecture \ref{ratioconj}.

%%%%%%%%%%%%%%%%%%%%%%%%%%%%%%%%%%%%%%%%%%%%%%%%%%%

\section{Proof of Theorem~\ref{thm1LDsymplectic} and Theorem \ref{analogoRMTsimplettico}}\label{C5S1}
The family $\{L(1/2,\chi_d):d>0$, $d$ fundamental discriminant$\}$ is a symplectic family% with conductor $\frac{d}{\pi}$
%(i.e. log-conductor $\log\frac{d}{2\pi}$ in (\ref{def1ldgenerale}))
, in the sense that it can be modeled by characteristic polynomials of symplectic matrices in the group $U\!Sp(2N)$, if we identify $2N\approx\log\frac{d}{\pi}$. 
Indeed $\frac{d}{\pi}$ is the analytic conductor of $L(s, \chi_d)$, thus $\log\frac{d}{\pi}$ (i.e. the density of zeros) plays the role of $2N$ in the random matrix theory setting\footnote[9]{See \cite[Conjecture 1.5.3]{CFKRSrecipe} and comments below for some clarification concerning the \lq\lq conductor\rq\rq).}.
%Analogous considerations can be done working with negative fundamental discriminants.

%\subsection{The ratio conjecture for $L(s,\chi_d)$}\label{C5S1S2}

We consider the moments of quadratic Dirichlet $L$-functions at the critical point $s=\frac{1}{2}$, i.e. the mean value
\begin{equation}\label{5.2.0} \sum_{d\leq X}L(\tfrac{1}{2},\chi_d)^k \end{equation}
in the limit $X\to\infty$, where the summation over $d$ has to be interpreted as the sum over all the positive fundamental discriminants $d$ below $X$, here and in the following. 
Also, we will denote by $X^*\sim\frac{1}{2\zeta(2)}X$ the number of fundamental discriminants below $X$.
We recall that Jutila \cite{JutilaDirichlet} proved asymptotic formulae for the first moment, showing that
\begin{equation}\label{5.2.1} \sum_{d\leq X}L(\tfrac{1}{2},\chi_d)\sim \frac{\mathcal A}{2}\frac{1}{2\zeta(2)}X\log X \end{equation}
where
\begin{equation}\label{5.2.2} \mathcal A=\prod_p\bigg(1-\frac{1}{p(p+1)}\bigg) \end{equation}
and also for the second moment, proving
\begin{equation}\label{5.2.3} \sum_{d\leq X}L(\tfrac{1}{2},\chi_d)^2\sim \frac{\mathcal B}{24}\frac{1}{2\zeta(2)}X(\log X)^3 \end{equation}
with
\begin{equation}\label{5.2.4} \mathcal B=\prod_p\bigg(1-\frac{4p^2-3p+1}{p^3(p+1)}\bigg). \end{equation}
It is believed that
\begin{equation}\label{5.2.5} \sum_{d\leq X}L(\tfrac{1}{2},\chi_d)^k\sim C_kX(\log X)^{k(k+1)/2} \end{equation}
and using analogies with random matrix theory, Keating and Snaith \cite{KSboh} also conjectured a precise value for the constant $C_k$. Moreover, the recipe produces a conjectural asymptotic formula with all the main terms for the moments~\eqref{5.2.0} with $k$ integer and also for ratios of products of quadratic Dirichlet $L$-functions (see \cite{CFZratioconj}), which is a symplectic analogue of Conjecture~\ref{ratioconj}.

\begin{conj}[\cite{CFZratioconj}, Conjecture 5.2]\label{ratioconjsympl}
Let $K,Q$ be two positive integers, $\alpha_1,\dots,\alpha_K$ and $\gamma_1,\dots,\gamma_Q$ complex shifts with real part $\asymp(\log T)^{-1}$ and imaginary part $\ll_\varepsilon T^{1-\varepsilon}$ for every $\varepsilon>0$, then
\begin{equation}\begin{split}\notag
&\sum_{d\leq X}\frac{\prod_{k=1}^KL(1/2+\alpha_k,\chi_d)}{\prod_{q=1}^QL(1/2+\gamma_q,\chi_d)}
\\&\hspace{1cm}=\sum_{d\leq X}\sum_{\epsilon\in\{-1,1\}^K}\Big(\frac{d}{\pi}\Big)^{\frac{1}{2}\sum_k(\epsilon_k\alpha_k-\alpha_k)}\prod_{k=1}^Kg_S\Big(\frac{1}{2}+\frac{\alpha_k-\epsilon_k\alpha_k}{2}\Big)Y_S\mathcal A_S(\cdots)
\\&\hspace{10.5cm}+O(X^{1/2+\varepsilon})
\\&\hspace{1cm}with\; (\cdots)=(\epsilon_1\alpha_1,\dots,\epsilon_K\alpha_K;\gamma)
\end{split}\end{equation}
where
\begin{equation}\notag Y_S(\alpha;\gamma):=\frac{\prod_{j\leq k\leq K}\zeta(1+\alpha_j+\alpha_k)\prod_{q<r\leq Q}\zeta(1+\gamma_q+\gamma_r)}{\prod_{k=1}^K\prod_{q=1}^Q\zeta(1+\alpha_k+\gamma_q)} \end{equation}
and $\mathcal A_S$ is an Euler product, absolutely convergent for all of the variables in small disks around 0, which is given by
\begin{equation}\begin{split}\notag
\mathcal A_S(\alpha;\gamma):=&\prod_p\frac{\prod_{j\leq k\leq K}(1-1/p^{1+\alpha_j+\alpha_k})\prod_{q<r\leq Q}(1-1/p^{1+\gamma_q+\gamma_r})}{\prod_{k=1}^K\prod_{q=1}^Q(1-1/p^{1+\alpha_k+\gamma_q})}
\\& \bigg(1+(1+1/p)^{-1}\sum_{0<\sum_ka_k+\sum_qc_q \text{ is even}}\frac{\prod_q\mu(p^{c_q})}{p^{\sum_ka_k(1/2+\alpha_k)+\sum_qc_q(1/2+\gamma_q)}}\bigg)
\end{split}\end{equation}
while \begin{equation}\notag g_S(z):=\frac{\Gamma(\frac{1-z}{2})}{\Gamma(\frac{z}{2})} .\end{equation}
\end{conj}
In particular, for our applications to the weighted one-level density, we are interested in the case $Q=1$, $2\leq K\leq 5$. 

\subsection{Conjecture \ref{ratioconjsympl} in the case K=2, Q=1.}
We start with
\begin{equation}\label{5.2.6} \sum_{d\leq X}\frac{L(1/2+A,\chi_d)L(1/2+C,\chi_d)}{L(1/2+D,\chi_d)}\end{equation}
with $A,C,D$ shifts, which satisfy the hypotheses prescribed by Conjecture~\ref{ratioconjsympl}; by the ratio conjecture, up to a negligible error $O(X^{1/2+\varepsilon})$, this is a sum of four terms and the first is
\begin{equation}\notag
\sum_{d\leq X}
\frac{\zeta(1+2A)\zeta(1+2C)\zeta(1+A+C)}{\zeta(1+A+D)\zeta(1+C+D)}\mathcal A(A,C;D)
\end{equation}
where
\begin{equation}\begin{split}\notag
&\mathcal A(A,C;D)=\prod_p(1-\tfrac{1}{p^{1+2A}})(1-\tfrac{1}{p^{1+2C}})(1-\tfrac{1}{p^{1+A+C}})(1-\tfrac{1}{p^{1+A+D}})^{-1}(1-\tfrac{1}{p^{1+C+D}})^{-1}
\\&\hspace{5.7cm}\cdot\bigg(1+\tfrac{p}{p+1}\sum_{0<a+c+d \text{ even}}\frac{\mu(p^d)}{p^{a(1/2+A)+c(1/2+C)+d(1/2+D)}}\bigg).
\end{split}\end{equation}
%By using $\sum_{m+n \text{ even}}x^my^n=\frac{1+xy}{(1-x^2)(1-y^2)}$ and $\sum_{m+n \text{ odd}}x^my^n=\frac{x+y}{(1-x^2)(1-y^2)}$, since only the cases $d=0$ and $d=1$ contribute, the last factor equals
%\begin{equation}\begin{split}\notag
%&1+\frac{p}{p+1}\Big[ \sum_{0<a+c \text{ even}}\frac{1}{p^{a(1/2+A)+c(1/2+C)}}+\sum_{a+c \text{ odd}}\frac{-1}{p^{a(1/2+A)+c(1/2+C)+(1/2+D)}} \Big]
%\\&=1+\frac{p}{p+1}\bigg[ \frac{1+\frac{1}{p^{1+A+C}}}{(1-\frac{1}{p^{1+2A}})(1-\frac{1}{p^{1+2C}})}
%-1-\frac{\frac{1}{p^{1+A+D}}+\frac{1}{p^{1+C+D}}}{(1-\frac{1}{p^{1+2A}})(1-\frac{1}{p^{1+2C}})} \bigg]
%\end{split}\end{equation}
%then
%\begin{equation}\begin{split}\notag
%\mathcal A(A&,C;D)\\&=\prod_p(1-\tfrac{1}{p^{1+2A}})(1-\tfrac{1}{p^{1+2C}})(1-\tfrac{1}{p^{1+A+C}})(1-\tfrac{1}{p^{1+A+D}})^{-1}(1-\tfrac{1}{p^{1+C+D}})^{-1}
%\\&\hspace{1cm}\bigg(1+\frac{p}{p+1}\Big[(1+\tfrac{1}{p^{1+A+C}})(1-\tfrac{1}{p^{1+2A}})^{-1}(1-\tfrac{1}{p^{1+2C}})^{-1}-1
%\\&\hspace{3cm}-(\tfrac{1}{p^{1+A+D}}+\tfrac{1}{p^{1+B+C}})(1-\tfrac{1}{p^{1+2A}})^{-1}(1-\tfrac{1}{p^{1+2C}})^{-1}\Big]\bigg).
%\end{split}\end{equation}
In the following, it will be relevant to notice that for small values of the shifts, then the arithmetical coefficient $\mathcal A(A,C;D)$ tends to $\mathcal A$, defined in~\eqref{5.2.2}; this essentially follows from \cite[Corollary 6.4]{CFZratioconj}.
%namely, if $A,C,D\to0$ then $\mathcal A(A,C;D)\sim\mathcal A(0,0;0)$ where
%\begin{equation}\begin{split}\notag
%\mathcal A(0,0;0)&= \prod_p(1-\tfrac{1}{p})^{3-2}\bigg(1+\tfrac{p}{p+1}\Big[(1+\tfrac{1}{p})(1-\tfrac{1}{p})^{-2}-1-\tfrac{2}{p}(1-\tfrac{1}{p})^{-2}\Big]\bigg)
%\\&=\prod_p(1-\tfrac{1}{p})\bigg(1+\tfrac{p^2}{(p-1)(p+1)}-\tfrac{p}{p+1}\bigg)
%=\prod_p\Big(\frac{p^2+p-1}{p(p+1)}\Big)=\mathcal A.
%\end{split}\end{equation}
All the other terms can be easily recovered from the first one, just by changes of sign of the shifts, as the recipe suggests. This yields a formula for~\eqref{5.2.6}, written as a sum of four pieces; by computing the derivative $\frac{d}{dC}[\cdots]_{C=D}$, we get
\begin{equation}\begin{split}\label{5.2.7}
&\sum_{d\leq X}\frac{L'}{L}(1/2+D,\chi_d)L(1/2+A,\chi_d)
\\&\hspace{1cm}=
\sum_{d\leq X}\bigg( 
Q_1+(\tfrac{d}{\pi})^{-A}g_S(\tfrac{1}{2}+A)Q_2+(\tfrac{d}{\pi})^{-D}g_S(\tfrac{1}{2}+D)Q_3
+
\\&\hspace{4.5cm}(\tfrac{d}{\pi})^{-A-D}g_S(\tfrac{1}{2}+A+D)Q_4
\bigg)+O(X^{1/2+\varepsilon})
\end{split}\end{equation}
with
\begin{equation}\begin{split}\notag
Q_1&
=\mathcal A(A,D;D)\tfrac{\zeta(1+2A)}{\zeta(1+A+D)}\bigg(\tfrac{2\zeta'(1+2D)\zeta(1+A+D)}{\zeta(1+2D)}+\tfrac{\zeta'(1+A+D)\zeta(1+2D)}{\zeta(1+2D)}
\\&\hspace{3.5cm}-\tfrac{\zeta'(1+2D)\zeta(1+A+D)}{\zeta(1+2D)}\bigg)+\mathcal A'(A,D;D)\zeta(1+2A)
\\&
=\mathcal A(A,D;D)\frac{\zeta(1+2A)}{\zeta(1+A+D)}\Big(\tfrac{\zeta'}{\zeta}(1+2D)\zeta(1+A+D)+\zeta'(1+A+D)\Big)
\\&\hspace{8cm}+\mathcal A'(A,D;D)\zeta(1+2A)
\end{split}\end{equation}
\begin{equation}\begin{split}\notag
Q_2&=
\mathcal A(-A,D;D)\frac{\zeta(1-2A)}{\zeta(1-A+D)}\Big(\tfrac{\zeta'}{\zeta}(1+2D)\zeta(1-A+D)+\zeta'(1-A+D)\Big)
\\&\hspace{8cm}+\mathcal A'(-A,D;D)\zeta(1-2A)
\end{split}\end{equation}
\begin{equation}\begin{split}\notag
Q_3=-\mathcal A(A,-D;D)\frac{\zeta(1+2A)\zeta(1-2D)\zeta(1+A-D)}{\zeta(1+A+D)}
\end{split}\end{equation}
\begin{equation}\begin{split}\notag
Q_4=-\mathcal A(-A,-D;D)\frac{\zeta(1-2A)\zeta(1-2D)\zeta(1-A-D)}{\zeta(1-A+D)}.
\end{split}\end{equation}
Moreover, we notice that if the shifts are $\ll(\log X)^{-1}$, then we can approximate the formula~\eqref{5.2.7}, getting
\begin{equation}\begin{split}\label{5.2.8}
&\sum_{d\leq X}\frac{L'}{L}(1/2+D,\chi_d)L(1/2+A,\chi_d)
\\&\hspace{.5cm}=
\mathcal AX^*\bigg(\frac{-A-3D}{(2A)(2D)(A+D)}
+X^{-A}\frac{A-3D}{(-2A)(2D)(-A+D)}
\\&\hspace{1cm}+X^{-D}\frac{A+D}{(2A)(2D)(A-D)}
+X^{-A-D}\frac{-A+D}{(-2A)(2D)(-A-D)}
\bigg)
\\&\hspace{1cm}+O(\log X)%\cdot\bigg(1+O\Big(\frac{1}{\log X}\Big)\bigg)
\end{split}\end{equation}
being that $\mathcal A(\pm A,\pm D,D)=\mathcal A+O(1/\log X)$ and $\zeta(1+z)=\frac{1}{z}+O(1)$ as $z\to0$.
%equals
%\begin{equation}\begin{split}\notag
%\sum_{d\leq X}\bigg(&
%Z_{A,C;D}\mathcal A(A,C;D)
%+\Big(\frac{d}{\pi}\Big)^{-A}g_S(\tfrac{A}{2})Z_{-A,C;D}\mathcal A(-A,C;D)
%\\&+\Big(\frac{d}{\pi}\Big)^{-C}g_S(\tfrac{C}{2})Z_{A,-C;D}\mathcal A(A,-C;D)
%+\Big(\frac{d}{\pi}\Big)^{-A-C}g_S(\tfrac{A}{2})g_S(\tfrac{C}{2})Z_{-A,-C;D}\mathcal A(-A,-C;D)
%\bigg)
%\end{split}\end{equation}
%with a small error $O(X^{1/2+\varepsilon})$, where
%\begin{equation}\notag

%\end{equation}

\subsection{Conjecture \ref{ratioconjsympl} in the case K=3, Q=1.}

Now we study in detail
\begin{equation}\label{5.2.10}
\sum_{d\leq X}\frac{L(1/2+A,\chi_d)L(1/2+B,\chi_d)L(1/2+C,\chi_d)}{L(1/2+D,\chi_d)}
\end{equation}
with $A,B,C,D$ as prescribed by Conjecture~\ref{ratioconjsympl}. This time, the asymptotic formula suggested by recipe is a sum of eight terms; the first is 
%\begin{equation}\notag
%\sum_{d\leq X}\tfrac{\zeta(1+2A)\zeta(1+2B)\zeta(1+2C)\zeta(1+A+B)\zeta(1+A+C)\zeta(1+B+C)}{\zeta(1+A+D)\zeta(1+B+D)\zeta(1+C+D)}\mathcal A(A,B,C;D)
%\end{equation}
\footnotesize
\begin{equation}\notag
\sum_{d\leq X}\frac{\zeta(1+2A)\zeta(1+2B)\zeta(1+2C)\zeta(1+A+B)\zeta(1+A+C)\zeta(1+B+C)}{\zeta(1+A+D)\zeta(1+B+D)\zeta(1+C+D)}\mathcal A(A,B,C;D)
\end{equation}
\normalsize
where the (rather horrible) arithmetical coefficient %, which can be recovered by noticing that $\sum_{m+n+h\text{ even}}x^my^nz^h=\frac{1+xy+zx+zy}{(1-x^2)(1-y^2)(1-z^2)}$ and similarly $\sum_{m+n+h\text{ odd}}x^my^nz^h=\frac{x+y+z+xyz}{(1-x^2)(1-y^2)(1-z^2)},$ 
is given by
\begin{equation}\begin{split}\notag
\mathcal A(A,B,C;D)=&
\prod_p(1-\tfrac{1}{p^{1+2A}})(1-\tfrac{1}{p^{1+2B}})(1-\tfrac{1}{p^{1+2C}})(1-\tfrac{1}{p^{1+A+B}})(1-\tfrac{1}{p^{1+A+C}})
\\&\hspace{.5cm}(1-\tfrac{1}{p^{1+B+C}})(1-\tfrac{1}{p^{1+A+D}})^{-1}(1-\tfrac{1}{p^{1+B+D}})^{-1}(1-\tfrac{1}{p^{1+C+D}})^{-1}
\\&\hspace{.5cm}\bigg(1+\tfrac{p}{p+1}
\Big[(1+\tfrac{1}{p^{1+A+B}}+\tfrac{1}{p^{1+A+C}}+\tfrac{1}{p^{1+B+C}})(1-\tfrac{1}{p^{1+2A}})^{-1}(1-\tfrac{1}{p^{1+2B}})^{-1}
\\&\hspace{.5cm}(1-\tfrac{1}{p^{1+2C}})^{-1}
-1-(\tfrac{1}{p^{1+A+D}}+\tfrac{1}{p^{1+B+D}}+\tfrac{1}{p^{1+C+D}}+\tfrac{1}{p^{2+A+B+C+D}})
\\&\hspace{.5cm}(1-\tfrac{1}{p^{1+2A}})^{-1}(1-\tfrac{1}{p^{1+2B}})^{-1}(1-\tfrac{1}{p^{1+2C}})^{-1}
\Big]\bigg).
\end{split}\end{equation}
We notice that, as in the proof of \cite[Corollary 6.4]{CFZratioconj}, we can prove that the arithmetical coefficient is convergent if all the variables are in small disk around 0, being
%\begin{equation}\begin{split}\notag
%\mathcal A(\underline 0)&=\prod_p(1-\tfrac{1}{p})^3\bigg(1+\tfrac{p}{p+1}
%\Big[(1+\tfrac{3}{p})(1-\tfrac{1}{p})^{-3}-1-(\tfrac{3}{p}+\tfrac{1}{p^2})(1-\tfrac{1}{p})^{-3}\Big]\bigg)
%\\& =\prod_p(1-\tfrac{1}{p})^{3}\Big(1+\tfrac{p^2}{(p-1)^2}-\tfrac{p}{p+1}\Big)
%=\prod_p\Big(1-\tfrac{4p^2-3p+1}{p^3(p+1)}\Big)=\mathcal B,
%\end{split}\end{equation}
$\mathcal A(\underline 0) = \mathcal B$, with $\mathcal B$ defined in \eqref{5.2.4}.
As in the previous example, this gives a formula for~\eqref{5.2.10} with all the main terms and error $O(X^{1/2+\varepsilon})$.
Differentiating this formula with respect to $C$ at $C=D$, we get
\begin{equation}\begin{split}\label{5.2.11}
\sum_{d\leq X}&\frac{L'}{L}(1/2+D,\chi_d)L(1/2+A,\chi_d)L(1/2+B,\chi_d)
\\&=\sum_{d\leq X}\bigg( 
R_1+(\tfrac{d}{\pi})^{-A}g_S(\tfrac{1}{2}+A)R_2+(\tfrac{d}{\pi})^{-B}g_S(\tfrac{1}{2}+B)R_3
\\&\hspace{1.5cm}+(\tfrac{d}{\pi})^{-D}g_S(\tfrac{1}{2}+D)R_4
+(\tfrac{d}{\pi})^{-A-B}g_S(\tfrac{1}{2}+A+B)R_5
\\&\hspace{1.5cm}+(\tfrac{d}{\pi})^{-A-D}g_S(\tfrac{1}{2}+A+D)R_6
+(\tfrac{d}{\pi})^{-B-D}g_S(\tfrac{1}{2}+B+D)R_7
\\&\hspace{1.5cm}+(\tfrac{d}{\pi})^{-A-B-D}g_S(\tfrac{1}{2}+A+B+D)R_8
\bigg)+O(X^{1/2+\varepsilon})
\end{split}\end{equation}
with
\begin{equation}\begin{split}\notag
&R_1=R_1(A,B,D)=\mathcal A(A,B,D;D)\tfrac{\zeta(1+2A)\zeta(1+2B)\zeta(1+A+B)}{\zeta(1+A+D)\zeta(1+B+D)}
\\&\hspace{2cm}\bigg(\tfrac{2\zeta'(1+2D)\zeta(1+A+D)\zeta(1+B+D)+\zeta(1+2D)\zeta'(1+A+D)\zeta(1+B+D)}{\zeta(1+2D)}
\\&\hspace{2cm}+\tfrac{\zeta(1+2D)\zeta(1+A+D)\zeta'(1+B+D)-\zeta(1+A+D)\zeta(1+B+D)\zeta'(1+2D)}{\zeta(1+2D)}
\bigg)
\\&\hspace{2cm}+{\zeta(1+2A)\zeta(1+2B)\zeta(1+A+B)}\mathcal A'(A,B,D;D)
\\&R_2=R_1(-A,B,D)
\\&R_3=R_1(A,-B,D)
\\&R_4=R_4(A,B,D)=-\tfrac{\zeta(1+2A)\zeta(1+2B)\zeta(1+A+B)\zeta(1-2D)\zeta(1+A-D)\zeta(1+B-D)}{\zeta(1+A+D)\zeta(1+B+D)}\\&\hspace{9cm}\cdot\mathcal A(A,B,-D;D)
\\&R_5=R_1(-A,-B,D)
\\&R_6=R_4(-A,B,D)
\\&R_7=R_4(A,-B,D)
\\&R_8=R_4(-A,-B,D).
\end{split}\end{equation}
If $A,B,D\ll(\log X)^{-1}$ the above formula simplifies a lot, since in this case
\begin{equation}\begin{split}\notag
R_1&=\frac{-AB-3AD-3BD-5D^2}{(2A)(2B)(2D)(A+B)(A+D)(B+D)}\mathcal B+O\Big(\frac{(\log X)^6}{(\log X)^3}\Big)%\bigg(1+O\Big(\frac{1}{\log X}\Big)\bigg)
\\&=:f(A,B,D)\mathcal B+O\Big((\log X)^3\Big)%\bigg(1+O\Big(\frac{1}{\log X}\Big)\bigg)
\end{split}\end{equation}
and
\begin{equation}\begin{split}\notag
R_4&=\frac{-(A+D)(B+D)}{(2A)(2B)(-2D)(A+B)(A-D)(B-D)}\mathcal B+O\Big(\frac{(\log X)^6}{(\log X)^3}\Big)%\bigg(1+O\Big(\frac{1}{\log X}\Big)\bigg)
\\&=:g(A,B,D)\mathcal B+O\Big((\log X)^3\Big)%\bigg(1+O\Big(\frac{1}{\log X}\Big)\bigg)
\end{split}\end{equation}
giving
\begin{equation}\begin{split}\label{5.2.12}
&\sum_{d\leq X}\frac{L'}{L}(1/2+D,\chi_d)L(1/2+A,\chi_d)L(1/2+B,\chi_d)
\\&=\mathcal BX^*\bigg( 
f(A,B,D)+X^{-A}f(-A,B,D)+X^{-B}f(A,-B,D)
\\&+X^{-D}g(A,B,D)+X^{-A-B}f(-A,-B,D)+X^{-A-D}g(-A,B,D)
\\&+X^{-B-D}g(A,-B,D)+X^{-A-B-D}g(-A,-B,D)
\bigg)+O\Big((\log X)^3\Big).%\bigg(1+O\Big(\frac{1}{\log X}\Big)\bigg).
\end{split}\end{equation}

Analogous (but longer) formulae can be obtained also in the cases $K=4,Q=1$ and $K=5,Q=1$. With exactly the same ideas (but much longer computations) also the case $K>5$, $Q=1$ can be dealt.

\subsection{The weighted one-level density for $\{L(\frac{1}{2},\chi_d)\}_d$}\label{C5S1S3}

%In this section, we want to perform similar computations as in Sections~\ref{C4S2S3} and~\ref{C4S2S4}, investigating the weighted one-level density of the non-trivial zeros of quadratic Dirichlet $L$-functions. To warm up the engines, we sketch what happens in the classical case; 
We recall that the one-level density for the symplectic family of quadratic Dirichlet $L$-functions has been studied originally by \"Ozluk and Snyder \cite{OS} and independently by Katz and Sarnak \cite{KatzSarnak2}\footnote[10]{See also \emph{\lq\lq Zeroes of Zeta Functions, their Spaces and their Spectral Nature\rq\rq} by Katz and Sarnak, the 1997 preprint  version of \cite{KatzSarnak2}.}, who proved that 
\begin{equation}\label{5.3.1}
\lim_{X\to\infty}\frac{1}{X^*}\sum_{d\leq X}\sum_{\gamma_d}f\Big(\frac{\log X}{2\pi}\gamma_d\Big)=\int_{-\infty}^{+\infty}f(x)\bigg(1-\frac{\sin(2\pi x)}{2\pi x}\bigg)dx
\end{equation}
under GRH, for any $f$ such that $\supp \widehat f\subset(-2,2)$. Moreover, Conrey and Snaith \cite{CSapplications} showed~\eqref{5.3.1} (also with lower order terms) with no constraint on the support of $\widehat f$, under the assumption of the ratio conjecture; namely, they consider $f$ a test function, holomorphic throughout the strip $|\Im(z)|<2$, even, real on the real line and such that $f(x)\ll1/(1+x^2)$ as $x\to \infty$ and they study
\begin{equation}\label{5.3.2}
\mathcal D_{0}^{\boldsymbol L_\chi}(f):=\frac{1}{X^*}\sum_{d\leq X}\sum_{\gamma_d}f\Big(\frac{\log X}{2\pi}\gamma_d\Big).
\end{equation}
As $X\to\infty$, they show that the above is asymptotic to
the right-hand side of \eqref{5.3.2},
which matches with the one-level density for the eigenvalues of the matrices from the symplectic group $U\!Sp(2N)$. In particular, we notice that the one-level density function $1-\frac{\sin(2\pi x)}{2\pi x}$ vanishes of order 2 at $x=0$, being $\sim\frac{2\pi^2}{3}x^2$ as $x\to0$. \par

Similarly to what we did in Section~\ref{C4S2}, we now want to compute the weighted one-level density in the symplectic case, tilted by $L(\frac{1}{2},\chi_d)$. We note that, differently from what happens in the Riemann zeta function case, here we are allowed to consider the first power as well, as $L(\frac{1}{2},\chi_d)$ is real. The analogue of~\eqref{tiltedmean} in this context is
\begin{equation}\label{5.3.3}
\mathcal D_{1}^{\boldsymbol L_\chi }(f):=\frac{1}{\sum_{d\leq X}L(\frac{1}{2},\chi_d)}\sum_{d\leq X}\sum_{\gamma_d}f\Big(\frac{\log X}{2\pi}\gamma_d\Big)L(\tfrac{1}{2},\chi_d)
\end{equation}
and via ratio conjecture in the form of Equation~\eqref{5.2.7} this can be studied asymptotically, as shown in the following result.

\begin{prop}\label{W1LD-S1}
Assume GRH and Conjecture~\ref{ratioconjsympl} for $K=2, Q=1$. For any test function $f$, holomorphic in the strip $\Im(z)<2$, even, real on the real line and such that $f(x)\ll 1/(1+x^2)$ as $x\to \infty$, we have
\begin{equation}\notag
\mathcal D_{1}^{\boldsymbol L_\chi }(f)
=\int_{-\infty}^{+\infty}f(x)W_{U\!Sp}^{1}(x)dx+O\Big(\frac{1}{\log X}\Big)
\end{equation}
 as $X\to\infty$, where
\begin{equation}\notag
W_{U\!Sp}^{1}(x):=1+\frac{\sin(2\pi x)}{2\pi x}-\frac{2\sin^2(\pi x)}{(\pi x)^2}
%\textcolor{red}{=1+\frac{\sin(2\pi x)}{2\pi x}-\frac{4(1-\cos(2\pi x))}{(2\pi x)^2}}%
. \end{equation}
%In particular we notice that $W_{USp}^{(1)}(x)\sim \frac{2\pi^4}{45}x^4$, as $x\to0$. 
\end{prop}

\proof
We start looking at 
\begin{equation}\label{5.3.4}
\frac{1}{\frac{\mathcal A}{2}X^*\log X}\sum_{d\leq X}\sum_{\gamma_d}f\Big(\frac{\log X}{2\pi}\gamma_d\Big)L(\tfrac{1}{2}+\alpha,\chi_d)
\end{equation}
with $\alpha\asymp(\log X)^{-1}$; note that, as $\alpha\to0$, $\sum_{d\leq X}L(\frac{1}{2}+\alpha,\chi_d)$ tends to $\frac{1}{2\zeta(2)}\frac{\mathcal A}{2}X\log X$ which is the normalization $\frac{\mathcal A}{2}X^*\log X$ we have in~\eqref{5.3.4}. As usual, we use the Cauchy theorem and the functional equation for $\frac{L'}{L}(s,\chi_d)$ to write 
\begin{equation}\label{5.3.5}
\frac{1}{\frac{\mathcal A}{2}X^*\log X}\sum_{d\leq X}\sum_{\gamma_d}f\Big(\frac{\log X}{2\pi}\gamma_d\Big)L(\tfrac{1}{2}+\alpha,\chi_d)
=-\mathcal J(\alpha)+2\mathcal I(\alpha)+O\Big(\frac{1}{\log X}\Big)
\end{equation}
where
\begin{equation}\begin{split}\notag
\mathcal J(\alpha):&=\frac{2}{\mathcal AX^*\log X}\frac{1}{2\pi}\int_{-\infty}^{+\infty}(-\log X)f\Big(\frac{y\log X}{2\pi}\Big)\sum_{d\leq X}L(\tfrac{1}{2}+\alpha,\chi_d)dy
\\& =\frac{-2}{\mathcal AX^*\log X}\int_{-\infty}^{+\infty}f(x)\sum_{d\leq X}\bigg(L(\tfrac{1}{2},\chi_d)+O(\tfrac{1}{\log X})\bigg)dx
\\&= -\int_{-\infty}^{+\infty}f(x)dx+O\Big(\frac{1}{\log X}\Big)
\end{split}\end{equation}
and
\begin{equation}\begin{split}\notag
\mathcal I(\alpha):&=
%\frac{2}{\mathcal AX^*\log X}\frac{1}{2\pi}\int_{-\infty}^{+\infty}f\Big(\frac{y\log X}{2\pi}\Big)\sum_{d\leq X}\frac{L'}{L}(\tfrac{1}{2}+\delta+iy,\chi_d)L(\tfrac{1}{2}+\alpha,\chi_d)dy
%\\&=
\frac{2}{\mathcal AX^*(\log X)^2}\int_{-\infty}^{+\infty}f(x)\sum_{d\leq X}\frac{L'}{L}(\tfrac{1}{2}+\delta+\tfrac{2\pi ix}{\log X},\chi_d)L(\tfrac{1}{2}+\alpha,\chi_d)dx
\end{split}\end{equation}
with $\delta\asymp(\log X)^{-1}$. Now we rely on the assumption of the ratio conjecture (in particular Equation~\eqref{5.2.7} and~\eqref{5.2.8}) to compute the sum over $d$; in particular, in the same way as in the the proof of Theorem~\ref{congettura}, by a truncation of the integral over $x$ and Taylor approximations, we get
\begin{equation}\notag
\mathcal I(\alpha)
=\frac{2}{(\log X)^2}\int_{-\infty}^{+\infty}f(x)g_{X}(\alpha,\delta+\tfrac{2\pi ix}{\log X})dx+O\Big(\frac{1}{\log X}\Big)
\end{equation}
where 
\begin{equation}\begin{split}\notag
g_{X}(\alpha,w):=
&\frac{-\alpha-3w}{(2\alpha)(2w)(\alpha+w)}
+X^{-\alpha}\frac{\alpha-3w}{(-2\alpha)(2w)(-\alpha+w)}
\\&\hspace{1cm}+X^{-w}\frac{\alpha+w}{(2\alpha)(2w)(\alpha-w)}
+X^{-\alpha-w}\frac{-\alpha+w}{(-2\alpha)(2w)(-\alpha-w)}.
\end{split}\end{equation}
The integral is regular at $\delta=0$ then, if we denote $\alpha=\frac{a}{\log X}$, we get
\begin{equation}\notag
\mathcal I(\tfrac{a}{\log X})=\frac{2}{(\log X)^2}\int_{-\infty}^{+\infty}f(x)g_{X}(\tfrac{a}{\log X},\tfrac{2\pi ix}{\log X})dx+O\Big(\frac{1}{\log X}\Big)
\end{equation}
which is regular at $a=0$; indeed, if we take the limit as $a\to0$ we get
\begin{equation}\notag
\mathcal I(0)=\lim_{a\to0}\mathcal I(\tfrac{a}{\log X})=2\int_{-\infty}^{+\infty}f(x)g(2\pi ix)dx+O\Big(\frac{1}{\log X}\Big)
\end{equation}
where
\begin{equation}\begin{split}\notag
g(w)&:=\lim_{a\to0}\bigg(
\frac{-a-3w}{(2a)(2w)(a+w)}+e^{-a}\frac{a-3w}{(-2a)(2w)(-a+w)}
\\&\hspace{2cm}+e^{-w}\frac{a+w}{(2a)(2w)(a-w)}+e^{-a-w}\frac{-a+w}{(-2a)(2w)(-a-w)}\bigg)
\\&=\frac{-we^{-w}-3w-4e^{-w}+4}{4w^2}.
\end{split}\end{equation}
Then
\begin{equation}\begin{split}\notag
\mathcal D_{1}^{\boldsymbol L_\chi}(f)&=\lim_{\alpha\to0}\frac{1}{\frac{\mathcal A}{2}X^*\log X}\sum_{d\leq X}\sum_{\gamma_d}f\Big(\frac{\log X}{2\pi}\gamma_d\Big)L(\tfrac{1}{2}+\alpha,\chi_d)
\\&=-\mathcal J(0)+2\mathcal I(0)+O\Big(\frac{1}{\log X}\Big)
\\&=\int_{-\infty}^{+\infty}f(x)\Big(1+4g(2\pi ix)\Big)dx+O\Big(\frac{1}{\log X}\Big)
\end{split}\end{equation}
and since $f$ is even the main term above equals
\begin{equation}\notag
%\int_{-\infty}^{+\infty}f(x)\Big(1+\frac{\sin(2\pi x)}{2\pi x}+4\frac{\cos(2\pi x)-1}{(2\pi x)^2}\Big)dx =
 \int_{-\infty}^{+\infty}f(x)\Big(1+\frac{\sin(2\pi x)}{2\pi x}-\frac{2\sin^2(\pi x)}{(\pi x)^2}\Big)dx .
\end{equation}
\endproof

Analogously, we can compute the weighted one-level density, tilted by the second power of $L(\frac{1}{2},\chi_d)$, i.e.
\begin{equation}\label{5.3.10}
\mathcal D_{2}^{\boldsymbol L_\chi}(f):=\frac{1}{\sum_{d\leq X}L(\frac{1}{2},\chi_d)^2}\sum_{d\leq X}\sum_{\gamma_d}f\Big(\frac{\log X}{2\pi}\gamma_d\Big)L(\tfrac{1}{2},\chi_d)^2
\end{equation}
under the assumption of Conjecture~\ref{ratioconjsympl}, in the case $K=3,Q=1$.

\begin{prop}\label{W1LD-S2}
Assume GRH and Conjecture~\ref{ratioconjsympl} for $K=3, Q=1$. For any function $f$ holomorphic in the strip $\Im(z)<2$, even, real on the real line and such that $f(x)\ll 1/(1+x^2)$ as $x\to \infty$, we have
\begin{equation}\notag
\mathcal D_{2}^{\boldsymbol L_\chi }(f)
=\int_{-\infty}^{+\infty}f(x)W_{U\!Sp}^{2}(x)dx+O\Big(\frac{1}{\log X}\Big)
\end{equation}
 as $X\to\infty$, where
\begin{equation}\begin{split}\notag
W_{U\!Sp}^{2}(x)&:=
1-\frac{\sin(2\pi x)}{2\pi x}
-\frac{24(1-\sin^2(\pi x))}{(2\pi x)^2}
+\frac{48\sin(2\pi x)}{(2\pi x)^3}
-\frac{96\sin^2(\pi x)}{(2\pi x)^4}. 
%\\& \textcolor{red}{
%=1-\frac{\sin(2\pi x)}{(2\pi x)}
%-\frac{12(1+\cos(2\pi x))}{(2\pi x)^2}
%+\frac{48\sin(2\pi x)}{(2\pi x)^3}
%-\frac{48(1-\cos(2\pi x))}{(2\pi x)^4}
%}
\end{split}\end{equation}
%In particular we notice that $W_{USp}^{(2)}(x)\sim \frac{2\pi^6}{1575}x^6$, as $x\to0$.
\end{prop}

\proof
The proof works like that of Proposition~\ref{W1LD-S1}; first, for $\alpha=\frac{a}{\log X}\asymp(\log X)^{-1}$ and $\beta=\frac{b}{\log X}\asymp(\log X)^{-1}$, we analyze
\begin{equation}\label{5.3.10}
\frac{1}{\frac{\mathcal B}{24}X^*(\log X)^3}\sum_{d\leq X}\sum_{\gamma_d}f\Big(\frac{\log X}{2\pi}\gamma_d\Big)L(\tfrac{1}{2}+\alpha,\chi_d)L(\tfrac{1}{2}+\beta,\chi_d)
\end{equation}
which can be written as
\begin{equation}\label{5.3.11} -\mathcal J(\alpha,\beta)+2\mathcal I(\alpha,\beta)+O\Big(\frac{1}{\log X}\Big)\end{equation}
where
\begin{equation}\begin{split}\label{5.3.12}
\mathcal J(\alpha,\beta)&:=\frac{-24\log X}{\mathcal BX^*(\log X)^3}\frac{1}{2\pi}\int_{-\infty}^{+\infty}f\Big(\frac{y\log X}{2\pi}\Big)\sum_{d\leq X}L(\tfrac{1}{2}+\alpha,\chi_d)L(\tfrac{1}{2}+\beta,\chi_d)
\\& =-\int_{-\infty}^{+\infty}f(x)dx+O\Big(\frac{1}{\log X}\Big)
\end{split}\end{equation}
and
\begin{equation}\begin{split}\notag
\mathcal I(\alpha,\beta):=&
\frac{24}{\mathcal BX^*(\log X)^4}\int_{-\infty}^{+\infty}f(x)
\\&\quad\times\sum_{d\leq X}\frac{L'}{L}(\tfrac{1}{2}+\delta+\tfrac{2\pi ix}{\log X},\chi_d)L(\tfrac{1}{2}+\alpha,\chi_d)L(\tfrac{1}{2}+\beta,\chi_d)dx
\end{split}\end{equation}
where $\delta\asymp(\log T)^{-1}$, as usual. With the usual machinery, the ratio conjecture (see Equations~\eqref{5.2.11} and~\eqref{5.2.12}) allows us to evaluate the sum over $d$; the resulting quantity is regular at $\delta=0$ and at $\alpha=\frac{a}{\log X}=0$, $\beta=\frac{b}{\log X}=0$, thus taking the limit we get
\begin{equation}\label{5.3.13}
\mathcal I(0,0)=24\int_{-\infty}^{+\infty}f(x)h(2\pi ix)dx
\end{equation}
with
\begin{equation}\notag
h(y):=\frac{y^3e^{-y}-5y^3+12y^2e^{-y}+12y^2+48ye^{-y}+48e^{-y}-48}{48y^4}.
\end{equation}
Putting all together, from~\eqref{5.3.10},~\eqref{5.3.11},~\eqref{5.3.12} and~\eqref{5.3.13}, we finally get
\begin{equation}\begin{split}\notag
\mathcal D_2^{\boldsymbol L_\chi}&%=\lim_{\alpha,\beta\to0}\frac{1}{\frac{\mathcal B}{24}X^*(\log X)^3}\sum_{d\leq X}\sum_{\gamma_d}f\Big(\frac{\log X}{2\pi}\gamma_d\Big)L(\tfrac{1}{2}+\alpha,\chi_d)L(\tfrac{1}{2}+\beta,\chi_d)
%\\&
=-\mathcal J(0,0)+2\mathcal I(0,0)+O\Big(\frac{1}{\log X}\Big)
\\&=\int_{-\infty}^{+\infty}f(x)\Big(1+48h(2\pi ix)\Big)dx+O\Big(\frac{1}{\log X}\Big).
\end{split}\end{equation}
Moreover, since $f$ is even, the main term equals
\begin{equation}\notag
\int_{-\infty}^{+\infty}f(x)
\Big(1-\frac{\sin(2\pi x)}{2\pi x}
-\frac{24(1-\sin^2(\pi x))}{(2\pi x)^2}
+\frac{48\sin(2\pi x)}{(2\pi x)^3}
-\frac{96\sin^2(\pi x)}{(2\pi x)^4}
\Big)dx.
\end{equation}
\endproof

In the same way, we study 
\begin{equation}\label{nuova5.3.10}
\mathcal D_{3}^{\boldsymbol L_\chi}(f):=\frac{1}{\sum_{d\leq X}L(\frac{1}{2},\chi_d)^3}\sum_{d\leq X}\sum_{\gamma_d}f\Big(\frac{\log X}{2\pi}\gamma_d\Big)L(\tfrac{1}{2},\chi_d)^3
\end{equation}
assuming Conjecture~\ref{ratioconjsympl}, in the case $K=4,Q=1$.

\begin{prop}\label{W1LD-S3}
Assume GRH and Conjecture~\ref{ratioconjsympl} for $K=4, Q=1$. For any function $f$ holomorphic in the strip $\Im(z)<2$, even, real on the real line and such that $f(x)\ll 1/(1+x^2)$ as $x\to \infty$, we have
\begin{equation}\notag
\mathcal D_{3}^{\boldsymbol L_\chi}(f)
=\int_{-\infty}^{+\infty}f(x)W_{U\!Sp}^{3}(x)dx+O\Big(\frac{1}{\log X}\Big)
\end{equation}
 as $X\to\infty$, where
\begin{equation}\begin{split}\notag
W_{U\!Sp}^{3}(x):=
1&+\frac{\sin(2\pi x)}{2\pi x}-\frac{12\sin^2(\pi x)}{(\pi x)^2}-\frac{240\sin(2\pi x)}{(2\pi x)^3}
\\&-\frac{15(6-10\sin^2(\pi x))}{(\pi x)^4}+\frac{2880\sin(2\pi x)}{(2\pi x)^5}-\frac{90\sin^2(\pi x)}{(\pi x)^6}.
\end{split}\end{equation}
%In particular we notice that $W_{USp}^{(2)}(x)\sim \frac{2\pi^6}{1575}x^6$, as $x\to0$.
\end{prop}

\proof
We consider $\alpha,\beta,\nu\in\mathbb R$ of size $\asymp1/\log X$, we denote
$$ \boldsymbol L_{\alpha,\beta,\nu}(\tfrac{1}{2},\chi_d):=L(\tfrac{1}{2}+\alpha,\chi_d)L(\tfrac{1}{2}+\beta,\chi_d)L(\tfrac{1}{2}+\nu,\chi_d) $$
and we look at
$$ \frac{1}{\sum_{d\leq X}L(\frac{1}{2},\chi_d)^3}\sum_{d\leq X}\sum_{\gamma_d}f\Big(\frac{\log X}{2\pi}\gamma_d\Big)\boldsymbol L_{\alpha,\beta,\nu}(\tfrac{1}{2},\chi_d) .$$
With the usual machinery we get that the above equals
\begin{equation}\notag
\int_{-\infty}^{+\infty}f(x)\bigg(1+2\frac{1}{\sum_{d\leq X}L(\frac{1}{2},\chi_d)^3}\sum_{d\leq X}\frac{L'}{L}(\tfrac{1}{2}+\delta+\tfrac{2\pi ix}{\log X},\chi_d)\boldsymbol L_{\alpha,\beta,\nu}(\tfrac{1}{2},\chi_d)\bigg)dx
\end{equation}
up to an error $O(1/\log X)$, with $\delta\asymp 1/\log X$. The remaining sum over $d$ can be evaluated asymptotically by using the ratio conjecture (i.e. Conjecture~\ref{ratioconjsympl} for $K=4,Q=1$). This can be done by using Sage to carry out the easy but very long computations. Doing so, letting $\alpha,\beta,\nu\to 0$, we obtain
\begin{equation}\notag
\mathcal D_3^{\boldsymbol L_{\chi}}(f)=\int_{-\infty}^{+\infty}f(x)\bigg(1+2\cdot2880\cdot h(2\pi i x)\bigg)dx+O\Big(\frac{1}{\log X}\Big)
\end{equation}
with
\begin{equation}\begin{split}\notag
h(y):=&\frac{-7y^5+24y^4-240y^2+2880}{5760y^6}\\&+\frac{e^{-y}(-y^5-24y^4-240y^3-1200y^2-2880y-2880)}{5760y^6}.
\end{split}\end{equation}
The claim follows, since $f$ is even.
\endproof

Finally, we look at 
\begin{equation}\label{nuova5.3.11}
\mathcal D_{4}^{\boldsymbol L_\chi}(f):=\frac{1}{\sum_{d\leq X}L(\frac{1}{2},\chi_d)^4}\sum_{d\leq X}\sum_{\gamma_d}f\Big(\frac{\log X}{2\pi}\gamma_d\Big)L(\tfrac{1}{2},\chi_d)^4
\end{equation}
assuming Conjecture~\ref{ratioconjsympl}, in the case $K=5,Q=1$.

\begin{prop}\label{W1LD-S4}
Assume GRH and Conjecture~\ref{ratioconjsympl} for $K=5, Q=1$. For any function $f$ holomorphic in the strip $\Im(z)<2$, even, real on the real line and such that $f(x)\ll 1/(1+x^2)$ as $x\to \infty$, we have
\begin{equation}\notag
\mathcal D_{4}^{\boldsymbol L_\chi}(f)
=\int_{-\infty}^{+\infty}f(x)W_{U\!Sp}^{4}(x)dx+O\Big(\frac{1}{\log X}\Big)
\end{equation}
 as $X\to\infty$, where
\begin{equation}\begin{split}\notag
W_{U\!Sp}^{4}(x):=
1&-\frac{\sin(2\pi x)}{2\pi x}-\frac{10(1+\cos(2\pi x))}{(\pi x)^2}+\frac{90\sin(2\pi x)}{(\pi x)^3}
\\&-\frac{15(3-31\cos(2\pi x))}{(\pi x)^4}-\frac{1470\sin(2\pi x)}{(\pi x)^5}
\\&-\frac{315(1+9\cos(2\pi x))}{(\pi x)^6}+\frac{3150\sin(2\pi x)}{(\pi x)^7}-\frac{1575(1-\cos(2\pi x))}{(\pi x)^8}.
\end{split}\end{equation}
%In particular we notice that $W_{USp}^{(2)}(x)\sim \frac{2\pi^6}{1575}x^6$, as $x\to0$.
\end{prop}

\proof
The proof works in the same way as the previous ones. We consider $\alpha,\beta,\nu,\eta\in\mathbb R$ of size $\asymp1/\log X$, we denote
$$ \boldsymbol L_{\alpha,\beta,\nu,\eta}(\tfrac{1}{2},\chi_d):=L(\tfrac{1}{2}+\alpha,\chi_d)L(\tfrac{1}{2}+\beta,\chi_d)L(\tfrac{1}{2}+\nu,\chi_d)L(\tfrac{1}{2}+\eta,\chi_d) $$
and we look at
$$ \frac{1}{\sum_{d\leq X}L(\frac{1}{2},\chi_d)^4}\sum_{d\leq X}\sum_{\gamma_d}f\Big(\frac{\log X}{2\pi}\gamma_d\Big)\boldsymbol L_{\alpha,\beta,\nu,\eta}(\tfrac{1}{2},\chi_d) .$$
By the usual manipulations, the above equals
\begin{equation}\notag
\int_{-\infty}^{+\infty}f(x)\bigg(1+2\frac{1}{\sum_{d\leq X}L(\frac{1}{2},\chi_d)^4}\sum_{d\leq X}\frac{L'}{L}(\tfrac{1}{2}+\delta+\tfrac{2\pi ix}{\log X},\chi_d)\boldsymbol L_{\alpha,\beta,\nu,\eta}(\tfrac{1}{2},\chi_d)\bigg)dx
\end{equation}
up to an error $O(1/\log X)$, with $\delta\asymp 1/\log X$. Thanks to Conjecture~\ref{ratioconjsympl} with $K=5,Q=1$, the above can be computed asymptotically. As $\alpha,\beta,\nu,\eta\to 0$, with the help of Sage, we then obtain
\begin{equation}\notag
\mathcal D_4^{\boldsymbol L_{\chi}}(f)=\int_{-\infty}^{+\infty}f(x)\bigg(1+2\cdot 4838400 \cdot h(2\pi i x)\bigg)dx+O\Big(\frac{1}{\log X}\Big)
\end{equation}
with
\begin{equation}\begin{split}\notag
h(y):=&\frac{-9y^7+40y^6-720y^4+20160y^2-403200}{9676800y^8}\\& + \frac{e^{-y}(y^7 + 40y^6 + 720y^5 + 7440y^4)}{9676800y^8}\\&+\frac{e^{-y}( 47040y^3 + 181440y^2 + 403200y + 403200)}{9676800y^8}.
\end{split}\end{equation}
Again, being $f$ even, the claim follows.
\endproof

\subsection{Proof of Theorem \ref{analogoRMTsimplettico}}
The same remark as in Section \ref{RMTaggiuntaUnitary} applies here, relying on \cite[Theorem 4.2]{CFZratioconj} instead of Conjecture \ref{ratioconjsympl}.

%%%%%%%%%%%%%%%%%%%%%%%%%%%%%%%%%%%%%%%%%%%%%%%%%%%%%%

\section{Proof of Theorem~\ref{thm1LDorthogonal} and Theorem \ref{analogoRMTortogonale}}\label{C5S2}

As a last example, we analyze the orthogonal case of the family of quadratic twists of the $L$-functions associated with the discriminant modular form $\Delta$. 
%
%
%
%
%
%
%
%
%
%
%
%
%
%
%
%
%\subsection{Quadratic twists of $L_\Delta(s)$}\label{C5S2S1}
%
%We start by denoting with $SL_2(\mathbb Z)$ the modular group, that is the group of $2\times2$ matrices with integer coefficients and determinant 1; let $k$ be a positive integer, then a modular form of weight $k$ is a complex-valued function $f$ on the upper half-plane $\mathbb H:=\{z\in\mathbb C:\Im(z)>0\}$ which is holomorphic on $\mathbb H$ and \lq\lq at the cusp\rq\rq (i.e. as $z\to i\infty$), such that for any $z\in\mathbb H$ we have
%\begin{equation}\notag 
%f\Big(\frac{az+b}{cz+d}\Big)=(cz+d)^kf(z)\quad\text{for any}\quad\begin{pmatrix} a & b \\  c & d \end{pmatrix}\in SL_2(\mathbb Z).
%\end{equation}
% Since $f(z+1)=f(z)$, a modular form is also periodic with period 1, then it has a Fourier series of the form
%$f(z)=\sum a_f(n)q^n$, with $q=e^{2\pi iz}$. Moreover a modular form is called cusp form if $a_f(0)=0.$ For any cusp form $f$, the Dirichlet generating series for the Fourier coefficients of $f$ is defined as
%\begin{equation}\notag 
%L_f(s):=\sum_{n=1}^\infty\frac{a_f(n)}{n^{s+k/2}}
%\end{equation}
%(see \cite[Chapter 14]{Iwaniec-Kowalski} for a complete account on the definition of Hecke $L$-functions and their standard properties). Here we are interested in the discriminant modular form $\Delta$, 
which is the unique normalized cusp form of weight 12. Its Fourier coefficients define the Ramanujan tau function $\tau(n)$, being
\begin{equation}\notag 
\Delta(z)=q \prod_{n\geq 1}(1-q^n)^{24}=\sum_{n=1}^\infty\tau(n)q^n
\end{equation} 
with $q=e^{2\pi iz}$. Thus the $L$-function associated with $\Delta$ is defined by
\begin{equation}\notag 
L_\Delta(s):=\sum_{n=1}^\infty\frac{\tau^*(n)}{n^s}
\end{equation} 
where $\tau^*(n)=\tau(n)/n^{11/2}$ %(the normalization $n^{(k-1)/2}$ is standard in literature)
. The family we want to describe is the collection of the quadratic twists of $L_\Delta$, that are
\begin{equation}\notag 
L_\Delta(s,\chi_d)=\sum_{n=1}^\infty \frac{\chi_d(n)\tau^*(n)}{n^s}
=\prod_p\bigg(1-\frac{\tau^*(p)\chi_d(p)}{p^s}+\frac{\chi_d(p^2)}{p^{2s}}\bigg)^{-1}
\end{equation} 
and, for $d>0$, they satisfy the functional equation
%\begin{equation}\notag 
%\xi_\Delta(s,\chi_d):=\Big(\frac{d^2}{4\pi^2}\Big)^{s/2}\Gamma(s+11/2)L_\Delta(s,\chi_d)
%=\xi_\Delta(1-s,\chi_d).
%\end{equation}
\begin{equation}\notag 
\Big(\frac{d^2}{4\pi^2}\Big)^{\frac{s}{2}}\Gamma(s+11/2)L_\Delta(s,\chi_d)
=\Big(\frac{d^2}{4\pi^2}\Big)^{\frac{1-s}{2}}\Gamma(1-s+11/2)L_\Delta(1-s,\chi_d).
\end{equation}
Finally we also record that 
\begin{equation}\notag 
\frac{1}{L_\Delta(s,\chi_d)}
=\prod_p\bigg(1-\frac{\tau^*(p)\chi_d(p)}{p^s}+\frac{\chi_d(p^2)}{p^{2s}}\bigg)
=:\sum_{n=1}^\infty \frac{\chi_d(n)\mu_\Delta(n)}{n^s}
\end{equation}
where $\mu_\Delta$ is the multiplicative function defined by $\mu_\Delta(p)=-\tau^*(p)$, $\mu_\Delta(p^2)=1$ and $\mu_\Delta(p^\alpha)=0$ if $\alpha\geq 3$.\\

The family $\{L_\Delta(1/2,\chi_d): d>0, f.d.\}$ is an even orthogonal family, modeled by the group $SO(2N)$ with the identification $2N\approx\log\frac{d^2}{4\pi^2}$.\\

%\subsection{The ratio conjecture for $L_\Delta(s,\chi_d)$}\label{C5S2S2}

The moments at the central value of $L$-functions associated with quadratic twists of a modular form have been studied extensively in recent years, but only the first moment \cite{BFH, IwaniecOrder, MurtyMurty} and partially the second \cite{SoundYoung, RSMomentsanddistr} have been obtained. It is known that such a family can be either symplectic or orthogonal, depending on the specific $L$-function we twist; in particular, if we start with the $L$-function associated with the discriminant modular form $\Delta$, then we are in the latter case. For an orthogonal family $\mathcal F$, ordered by the conductor $C(f)$, Conrey-Farmer \cite{ConreyFarmer} and Keating-Snaith \cite{KSboh} predict that 
\begin{equation}\label{5.2.2.1}
\frac{1}{X^*}\sum_{\substack{f\in\mathcal F \\ C(f)\leq X}}L_f(1/2)^k\sim  \frac{f_O(k)}{2}a(k)(\log X^A)^{k(k-1)/2}
\end{equation}
where the above sum is over the $X^*$ elements of the family $\mathcal F$ such that $C(f)\leq X$; $f_O(k)$ is the leading order coefficient of the moments of characteristic polynomials of matrices in $SO(2N)$; $a(k)$ is a constant depending on the particular family involved; $A$ is a constant depending on the functional equation satisfied by the $L$-functions in the family, in particular on the degree of the relevant parameter in the functional equation for $L_f(s)$ (see \cite[Equation (1.3)]{ConreyFarmer} for further details and examples). Moreover, in this case the recipe \cite{CFKRSrecipe} provides a precise formula with all the main terms for any integral moment, extended by \cite{CFZratioconj} to ratios. The ratio conjecture for the orthogonal family of quadratic twists of the discriminant modular form can be stated as follows.

\begin{conj}[\cite{CFZratioconj}, Conjecture 5.3]\label{RCorthogonal}
Let $K,Q$ two positive integers, $\alpha_1,\dots,\alpha_K$ and $\gamma_1,\dots,\gamma_Q$ complex shifts with real part $\asymp(\log X)^{-1}$ and imaginary part $\ll_\varepsilon X^{1-\varepsilon}$ for every $\varepsilon>0$, then
\begin{equation}\begin{split}\notag
&\sum_{d\leq X}\frac{\prod_{k=1}^KL_\Delta(1/2+\alpha_k,\chi_d)}{\prod_{q=1}^QL_\Delta(1/2+\gamma_q,\chi_d)}
\\&\hspace{.8cm}=\sum_{d\leq X}\sum_{\epsilon\in\{-1,1\}^K}\Big(\frac{d^2}{4\pi^2}\Big)^{\frac{1}{2}\sum_k(\epsilon_k\alpha_k-\alpha_k)}\prod_{k=1}^Kg_O\Big(\frac{1}{2}+\frac{\alpha_k-\epsilon_k\alpha_k}{2}\Big)Y_O\mathcal A_O(\cdots)
\\&\hspace{10.5cm}+O(X^{1/2+\varepsilon})
\\&\hspace{1cm}with\; (\cdots)=(\epsilon_1\alpha_1,\dots,\epsilon_K\alpha_K;\gamma)
\end{split}\end{equation}
where
\begin{equation}\notag Y_O(\alpha;\gamma):=\frac{\prod_{j< k\leq K}\zeta(1+\alpha_j+\alpha_k)\prod_{q<r\leq Q}\zeta(1+\gamma_q+\gamma_r)\prod_{q\leq Q}\zeta(1+2\gamma_q)}{\prod_{k=1}^K\prod_{q=1}^Q\zeta(1+\alpha_k+\gamma_q)} \end{equation}
and $\mathcal A_O$ is an Euler product, absolutely convergent for all of the variables in small disks around 0, which is given by
\begin{equation}\begin{split}\notag
\mathcal A_O(\alpha;\gamma):=&\prod_p\frac{\prod_{j< k\leq K}(1-\frac{1}{p^{1+\alpha_j+\alpha_k}})\prod_{q<r\leq Q}(1-\frac{1}{p^{1+\gamma_q+\gamma_r}})\prod_{q\leq Q}(1-\frac{1}{p^{1+2\gamma_q}})}{\prod_{k=1}^K\prod_{q=1}^Q(1-\frac{1}{p^{1+\alpha_k+\gamma_q}})}
\\& \bigg(1+(1+1/p)^{-1}\sum_{0<\sum_ka_k+\sum_qc_q \text{ is even}}\frac{\prod_k\tau^*(p^{a_k})\prod_q\mu_\Delta(p^{c_q})}{p^{\sum_ka_k(1/2+\alpha_k)+\sum_qc_q(1/2+\gamma_q)}}\bigg)
\end{split}\end{equation}
while \begin{equation}\notag g_O(s):=\frac{\Gamma(\frac{1}{2}-s+6)}{\Gamma(s-\frac{1}{2}+6)} .\end{equation}
\end{conj}

In the following, we will analyze the applications of this conjecture to the weighted one-level density, as we did in Section~\ref{C5S1S3} for a symplectic family. To do so, we first look at what Conjecture~\ref{RCorthogonal} gives in a few specific examples.

\subsection{Conjecture \ref{RCorthogonal} in the case K=1, Q=0.}
This is the easiest situation possible, corresponding to the first moment of $L_\Delta(1/2,\chi_d)$; for $A$ a complex number which satisfies the hypotheses prescribed by Conjecture~\ref{RCorthogonal}, the ratio conjecture yields
\begin{equation}\notag
\frac{1}{X^*}\sum_{d\leq X}L_\Delta(1/2+A,\chi_d)=\mathcal A(A)+\Big(\frac{d}{2\pi}\Big)^{-2A}\frac{\Gamma(6-A)}{\Gamma(6+A)} \mathcal A(-A)+O(X^{-1/2+\varepsilon})
\end{equation}
with 
\begin{equation}\notag
\mathcal A(A):=\prod_p\bigg(1+\frac{p}{p+1}\bigg[-1+\sum_{m=0}^{\infty}\frac{\tau^*(p^{2m})}{p^{(\frac{1}{2}+A)2m}}\bigg]\bigg)
\end{equation}
We note that $\mathcal A(A)$ is regular at $A=0$; indeed the $m=0$ and $m=1$ terms give $1$ and $\tau^*(p^2)p^{-1-2A}$ respectively, therefore an approximation for $\mathcal A(A)$ would be $\prod_p(1+\frac{\tau^*(p^2)}{p^{1+2A}}+\cdots)$. Differently from the unitary and symplectic cases, where the first term in the corresponding Euler products gives the polar factor $\zeta(1+2A)$, here we would have $L_\Delta(\text{sym}^2,1+2A)$ the symmetric square of $L_\Delta$, which is well-known to be regular and nonzero at $1$ (see \cite[Chapter 13]{IwaniecTopics} for a complete overview about the symmetric square and its properties). However, for the sake of brevity, we prefer not to factor out $L_\Delta(\text{sym}^2,1+2A)$ and we leave the contribution of the symmetric square encoded in the arithmetical factor $\mathcal A(A)$, which converges in a small disk around 0. Thus, for $A=0$, we immediately get
\begin{equation}\label{5.2.2.2}
\frac{1}{X^*}\sum_{d\leq X}L_\Delta(1/2,\chi_d)=2\mathcal A+O(X^{-1/2+\varepsilon})
\end{equation}
where
\begin{equation}\label{5.2.2.3}
\mathcal A:=\mathcal A(0)=\prod_p\bigg(1+\frac{p}{p+1}\bigg[-1+\sum_{m=0}^{\infty}\frac{\tau^*(p^{2m})}{(\sqrt p)^{2m}}\bigg]\bigg).
\end{equation}
%The sum over $m$ in the above formula will recur often in the following, thus we denote it by $\mathcal P$; a closed formula for $\mathcal P$ is given by \cite[Equation (2.49)]{CSapplications}
%\begin{equation}\notag
%\mathcal P:=\sum_{m=0}^{\infty}\frac{\tau^*(p^{2m})}{(\sqrt p)^{2m}}
%=\frac{1}{2}\bigg\{\bigg(1-\frac{\tau^*(p)}{\sqrt p}+\frac{1}{p}\bigg)^{-1}+\bigg(1+\frac{\tau^*(p)}{\sqrt p}+\frac{1}{p}\bigg)^{-1}\bigg\}.
%\end{equation}

\subsection{Conjecture \ref{RCorthogonal} in the case K=2, Q=1.}
We consider
\begin{equation}\label{5.2.2.5} 
\sum_{d\leq X}\frac{L_\Delta(1/2+A,\chi_d)L_\Delta(1/2+C,\chi_D)}{L_\Delta(1/2+D,\chi_d)}
\end{equation}
with $A,C,D$ shifts satisfying the usual hypotheses prescribed by the ratio conjecture; by Conjecture~\ref{RCorthogonal}, up to a negligible error, this is a sum of four terms and the first is
\begin{equation}\notag
\sum_{d\leq X}
\frac{\zeta(1+A+C)\zeta(1+2D)}{\zeta(1+A+D)\zeta(1+C+D)}\mathcal A(A,C;D)
\end{equation}
where
\begin{equation}\begin{split}\notag
&\mathcal A(A,C;D)=\prod_p(1-\tfrac{1}{p^{1+A+C}})(1-\tfrac{1}{p^{1+2D}})(1-\tfrac{1}{p^{1+A+D}})^{-1}(1-\tfrac{1}{p^{1+C+D}})^{-1}
\\&\hspace{3.7cm}\cdot\bigg(1+\tfrac{p}{p+1}\sum_{0<a+c+d \text{ even}}\frac{\tau^*(p^a)\tau^*(p^c)\mu(p^d)}{p^{a(1/2+A)+c(1/2+C)+d(1/2+D)}}\bigg)
\end{split}\end{equation}
As usual, we note that $\mathcal A(A,C;D)\sim\mathcal A(0,0;0)=\mathcal A$ defined in~\eqref{5.2.2.3} as $A,C,D\to0$; this can be proved by a modification of the proof of \cite[Corollary 6.4]{CFZratioconj} or by direct computation. 
All the other terms can be easily recovered from the first one, then we get a formula for~\eqref{5.2.2.5}, written as a sum of four pieces; by computing the derivative $\frac{d}{dC}[\cdots]_{C=D}$, we get
\begin{equation}\begin{split}\label{5.2.2.10}
&\sum_{d\leq X}\frac{L'_\Delta}{L_\Delta}(1/2+D,\chi_d)L(1/2+A,\chi_d)
\\&\hspace{.5cm}=
\sum_{d\leq X}\bigg( 
Q_1+(\tfrac{d}{2\pi})^{-2A}g_O(\tfrac{1}{2}+A)Q_2+(\tfrac{d}{2\pi})^{-2D}g_O(\tfrac{1}{2}+D)Q_3
+
\\&\hspace{4cm}(\tfrac{d}{2\pi})^{-2A-2D}g_O(\tfrac{1}{2}+A+D)Q_4
\bigg)+O(X^{1/2+\varepsilon})
\end{split}\end{equation}
with
\begin{equation}\notag
Q_1=\mathcal A(A,D;D)\Big(\tfrac{\zeta'}{\zeta}(1+A+D)-\tfrac{\zeta'}{\zeta}(1+2D)\Big)+\mathcal A'(A,D;D)
\end{equation}
\begin{equation}\notag
Q_2=\mathcal A(-A,D;D)\Big(\tfrac{\zeta'}{\zeta}(1-A+D)-\tfrac{\zeta'}{\zeta}(1+2D)\Big)+\mathcal A'(-A,D;D)
\end{equation}
\begin{equation}\notag
Q_3=-\mathcal A(A,-D;D)\frac{\zeta(1+A-D)\zeta(1+2D)}{\zeta(1+A+D)}
\end{equation}
\begin{equation}\notag
Q_4=-\mathcal A(-A,-D;D)\frac{\zeta(1-A-D)\zeta(1+2D)}{\zeta(1-A+D)}
\end{equation}
Moreover, we notice that if the shifts are of order $\asymp(\log X)^{-1}$, then we can approximate the formula~\eqref{5.2.2.10}, getting
\begin{equation}\begin{split}\label{5.2.2.11}
&\sum_{d\leq X}\frac{L'_\Delta}{L_\Delta}(1/2+D,\chi_d)L(1/2+A,\chi_d)
\\&\hspace{1cm}=
\mathcal AX^*\bigg(\frac{A-D}{(A+D)2D}
+X^{-2A}\frac{-A-D}{(-A+D)2D}
\\&\hspace{1.5cm}+X^{-2D}\frac{-A-D}{(A-D)2D}
+X^{-2A-2D}\frac{A-D}{(-A-D)2D}
\bigg)
+O(1)%\\&\hspace{1.5cm}\cdot\bigg(1+O\Big(\frac{1}{\log X}\Big)\bigg)
\end{split}\end{equation}
being $\mathcal A(\pm A,\pm D,D)=\mathcal A+O(1/\log X)$ and $\zeta(1+z)=\frac{1}{z}+O(1)$ as $z\to0$.

\subsection{Conjecture \ref{RCorthogonal} in the case K=2, Q=0.}
We now analyze closely the second moment of $L_\Delta(1/2,\chi_d)$; we take two complex shifts $A,B$ such that $A,B\asymp (\log X)^{-1}$ and we look at
\begin{equation}\notag
\frac{1}{X^*}\sum_{d\leq X}L_\Delta(1/2+A,\chi_d)L_\Delta(1/2+B,\chi_d).
\end{equation}
By Conjecture~\ref{RCorthogonal}, ignoring the negligible error term $O(X^{1/2+\varepsilon})$, the above is 
\begin{equation}\label{5.aux1}
f(A,B)+X^{-2A}f(-A,B)+X^{-2B}f(A,-B)+X^{-2A-2B}f(-A-B)
\end{equation}
with 
\begin{equation}\notag
f(A,B):=\zeta(1+A+B)\prod_p\bigg(1-\frac{1}{p^{1+A+B}}\bigg)\bigg(1+\frac{p}{p+1}\sum_{\substack{m+n>0 \\ \text{even}}}\frac{\tau^*(p^m)\tau^*(p^n)}{p^{(\frac{1}{2}+A)m+(\frac{1}{2}+B)n}}\bigg)
\end{equation}
Since $A,B\asymp(\log X)^{-1}$, we set $a=A\log X\asymp 1$ and $b=B\log X\asymp 1$, so that~\eqref{5.aux1} becomes
\begin{equation}\label{5.aux2}
\mathcal B\log X\bigg(\frac{1}{a+b}+\frac{e^{-2a}}{-a+b}+\frac{e^{-2b}}{a-b}+\frac{e^{-2a-2b}}{-a-b}\bigg)\bigg(1+O\Big(\frac{1}{\log X}\Big)\bigg)
\end{equation}
where
\begin{equation}\begin{split}\label{5.aux3}
\mathcal B:&=\prod_p\bigg(1-\frac{1}{p}\bigg)\bigg(1+\frac{p}{p+1}\sum_{\substack{m+n>0 \\ \text{even}}}\frac{\tau^*(p^m)\tau^*(p^n)}{p^{(m+n)/2}}\bigg).
%\\& =\prod_p\bigg(1-\frac{1}{p}\bigg)\bigg(1+\frac{p}{p+1}(-1+\mathcal P^2+\mathcal D^2)\bigg)
\end{split}\end{equation}
%with $\mathcal P$ and $\mathcal D$ defined in~\eqref{5.2.2.7} and~\eqref{5.2.2.8} respectively. 
The expression in~\eqref{5.aux2} is regular at $a=0$ and $b=0$, since the limit of the first parentheses as $a,b\to0$ equals 4, therefore we finally get
\begin{equation}\begin{split}\label{5.aux4}
\frac{1}{X^*}\sum_{d\leq X}L_\Delta(1/2,\chi_d)^2\sim 4\mathcal B\log X.
\end{split}\end{equation}

\subsection{Conjecture \ref{RCorthogonal} in the case K=3, Q=1.}
Finally we look at
\begin{equation}\label{5.aux10}
\sum_{d\leq X}\frac{L_\Delta(1/2+A,\chi_d)L_\Delta(1/2+B,\chi_d)L_\Delta(1/2+C,\chi_d)}{L_\Delta(1/2+D,\chi_d)}
\end{equation}
with $A,B,C,D$ as Conjecture~\ref{RCorthogonal} prescribes. The first of the eight terms given by the recipe is
\begin{equation}\notag
\sum_{d\leq X}\frac{\zeta(1+2D)\zeta(1+A+B)\zeta(1+A+C)\zeta(1+B+C)}{\zeta(1+A+D)\zeta(1+B+D)\zeta(1+C+D)}\mathcal A(A,B,C;D)
\end{equation}
where 
\begin{equation}\begin{split}\notag
\mathcal A(A,B,C;D)=&\prod_p \frac{(1-\frac{1}{p^{1+2D}})(1-\frac{1}{p^{1+A+B}})(1-\frac{1}{p^{1+A+C}})(1-\frac{1}{p^{1+B+C}})}{(1-\frac{1}{p^{1+A+D}})(1-\frac{1}{p^{1+B+D}})(1-\frac{1}{p^{1+C+D}})}
\\&\quad \bigg(1+\frac{p}{p+1}\sum_{\substack{0<a+b+c+d \\ \text{even}}}\frac{\tau^*(p^a)\tau^*(p^b)\tau^*(p^c)\mu_\Delta(p^d)}{p^{(\frac{1}{2}+A)a+(\frac{1}{2}+B)b+(\frac{1}{2}+C)c+(\frac{1}{2}+D)d}}\bigg)
\end{split}\end{equation}
is the arithmetical coefficient, absolutely convergent in small disks around 0, such that
$\mathcal A(0,0,0;0) = \mathcal B$.
%\begin{equation}\begin{split}\notag
%\mathcal A(0,0,0;0)&=\prod_p \Big(1-\frac{1}{p}\Big)\bigg(1+\frac{p}{p+1}\sum_{\substack{0<a+b+c+d \\ \text{even}}}\frac{\tau^*(p^a)\tau^*(p^b)\tau^*(p^c)\mu_\Delta(p^d)}{p^{(a+b+c+d)/2}}\bigg)
%\\& = \prod_p\bigg(1-\frac{1}{p}\bigg)\bigg(1+\frac{p}{p+1}\bigg(\Big(1-\frac{1}{p}\Big)\sum_{\substack{a+b+c \\ \text{even}}}\frac{\tau^*(p^a)\tau^*(p^b)\tau^*(p^c)}{(\sqrt p)^{a+b+c}}
%\\&\hspace{4.3cm}-\frac{\tau^*(p)}{\sqrt p}\sum_{\substack{a+b+c \\ \text{odd}}}\frac{\tau^*(p^a)\tau^*(p^b)\tau^*(p^c)}{(\sqrt p)^{a+b+c}}-1\bigg)\bigg)
%\\& =\prod_p\bigg(1-\frac{1}{p}\bigg)\bigg(1+\frac{p}{p+1}(-1+\mathcal P^2+\mathcal D^2)\bigg)=\mathcal B
%\end{split}\end{equation}
%where from the second line to the third, we used the formulae 
%\begin{equation}\notag
%\sum_{a+b+c  \text{ even}}\frac{\tau^*(p^a)\tau^*(p^b)\tau^*(p^c)}{p^{a/2}p^{b/2}p^{c/2}}
%=\mathcal P(\mathcal P^2+3\mathcal D^2)
%\end{equation}
%\begin{equation}\notag
%\sum_{a+b+c  \text{ odd}}\frac{\tau^*(p^a)\tau^*(p^b)\tau^*(p^c)}{p^{a/2}p^{b/2}p^{c/2}}
%=\mathcal D(\mathcal D^2+3\mathcal P^2)
%\end{equation}
%and the multiplicative law of the Ramanujan tau function. 
As in all the previous examples, this gives a formula for~\eqref{5.aux10} with all the main terms and error $O(X^{1/2+\varepsilon})$ and differentiating this formula with respect to $C$ at $C=D$, we get
\begin{equation}\begin{split}\label{5.aux11}
\sum_{d\leq X}&\frac{L'_\Delta}{L_\Delta}(1/2+D,\chi_d)L_\Delta(1/2+A,\chi_d)L_\Delta(1/2+B,\chi_d)
\\&=\sum_{d\leq X}\bigg( 
R_1+(\tfrac{d}{2\pi})^{-2A}g_O(\tfrac{1}{2}+A)R_2+(\tfrac{d}{2\pi})^{-2B}g_O(\tfrac{1}{2}+B)R_3
\\&\hspace{1.3cm}+(\tfrac{d}{2\pi})^{-2D}g_O(\tfrac{1}{2}+D)R_4
+(\tfrac{d}{2\pi})^{-2A-2B}g_O(\tfrac{1}{2}+A+B)R_5
\\&\hspace{1.3cm}+(\tfrac{d}{2\pi})^{-2A-2D}g_O(\tfrac{1}{2}+A+D)R_6
+(\tfrac{d}{2\pi})^{-2B-2D}g_O(\tfrac{1}{2}+B+D)R_7
\\&\hspace{1.3cm}+(\tfrac{d}{2\pi})^{-2A-2B-2D}g_O(\tfrac{1}{2}+A+B+D)R_8
\bigg)+O(X^{1/2+\varepsilon})
\end{split}\end{equation}
with
\begin{equation}\begin{split}\notag
&R_1=R_1(A,B,D)=\mathcal A(A,B,D;D)\zeta(1+A+B)\Big(\tfrac{\zeta'}{\zeta}(1+A+D)
\\&\hspace{2cm}+\tfrac{\zeta'}{\zeta}(1+B+D)-\tfrac{\zeta'}{\zeta}(1+2D)\Big)+\zeta(1+A+B)\mathcal A'(A,B,D;D)
\\&R_2=R_1(-A,B,D)
\\&R_3=R_1(A,-B,D)
\\&R_4=R_4(A,B,D)=-\frac{\zeta(1+2D)\zeta(1+A+B)\zeta(1+A-D)\zeta(1+B-D)}{\zeta(1+A+D)\zeta(1+B+D)}\\&\hspace{10cm}\cdot\mathcal A(A,B,-D;D)
\\&R_5=R_1(-A,-B,D)
\\&R_6=R_4(-A,B,D)
\\&R_7=R_4(A,-B,D)
\\&R_8=R_4(-A,-B,D).
\end{split}\end{equation}
If $A,B,D\asymp(\log X)^{-1}$ the above formula simplifies a lot, since
\begin{equation}\begin{split}\notag
R_1&=\frac{AB-AD-BD-3D^2}{2D(A+B)(A+D)(B+D)}\mathcal B+O\Big(\frac{(\log X)^4}{(\log X)^3}\Big)%\bigg(1+O\Big(\frac{1}{\log X}\Big)\bigg)
\\&=:f(A,B,D)+O(\log X)%\bigg(1+O\Big(\frac{1}{\log X}\Big)\bigg)
\end{split}\end{equation}
and
\begin{equation}\begin{split}\notag
R_4&=\frac{(A+D)(B+D)}{(-2D)(A+B)(A-D)(B-D)}\mathcal B+\Big(\frac{(\log X)^4}{(\log X)^3}\Big)%\bigg(1+O\Big(\frac{1}{\log X}\Big)\bigg)
\\&=:g(A,B,D)+O(\log X)%\bigg(1+O\Big(\frac{1}{\log X}\Big)\bigg)
\end{split}\end{equation}
giving
\begin{equation}\begin{split}\label{5.aux20}
\sum_{d\leq X}&\frac{L'_\Delta}{L_\Delta}(1/2+D,\chi_d)L_\Delta(1/2+A,\chi_d)L_\Delta(1/2+B,\chi_d)
\\&=\mathcal BX^*\bigg( 
f(A,B,D)+X^{-2A}f(-A,B,D)+X^{-2B}f(A,-B,D)
\\&+X^{-2D}g(A,B,D)+X^{-2A-2B}f(-A,-B,D)+X^{-2A-2D}g(-A,B,D)
\\&+X^{-2B-2D}g(A,-B,D)+X^{-2A-2B-2D}g(-A,-B,D)
\bigg)+O(\log X).%\bigg(1+O\Big(\frac{1}{\log X}\Big)\bigg).
\end{split}\end{equation}

Analogous formulae can be obtained in the cases $K=4,Q=1$ and $K=5,Q=1$. Again, with the same technique one can get formulae also in the case $K>5$, $Q=1$.

\subsection{The weighted one-level density for $\{L_\Delta(\tfrac{1}{2},\chi_d)\}_d$}\label{C5S2S3}

In analogy to what we did in Section~\ref{C5S1S3}, we now compute the weighted one-level density for the orthogonal family of quadratic twists of $L_\Delta$. We assume the Riemann Hypothesis for the $L$-functions we are considering and we denote with $\gamma_{\Delta,d}$ the imaginary part of a generic zero of $L_\Delta(s,\chi_d)$. In the classical case, assuming the ratio conjecture, Conrey and Snaith \cite{CSapplications} proved that  
\begin{equation}\label{5.6.1}
\lim_{X\to\infty}\frac{1}{X^*}\sum_{d\leq X}\sum_{\gamma_{\Delta,d}}f\Big(\frac{\log X}{\pi}\gamma_{\Delta,d}\Big)=\int_{-\infty}^{+\infty}f(x)\bigg(1+\frac{\sin(2\pi x)}{2\pi x}\bigg)dx
\end{equation}
for any test function $f$, satisfying the usual properties as in Theorem~\ref{thm1LDorthogonal}. We now use the formulae of the previous section to derive the weighted one-level density; we denote 
\begin{equation}\label{5.6.2}
\mathcal D_1^{\boldsymbol L_{\Delta,\chi}}(f):=\frac{1}{\sum_{d\leq X}L_\Delta(\frac{1}{2},\chi_d)}\sum_{d\leq X}\sum_{\gamma_{\Delta,d}}f\Big(\frac{\log X}{\pi}\gamma_{\Delta,d}\Big)L_\Delta(\tfrac{1}{2},\chi_d)
\end{equation}
and we prove the following result.

\begin{prop}\label{W1LD-O1}
Assume GRH and Conjecture~\ref{RCorthogonal} for $K=2, Q=1$. For any function $f$ holomorphic in the strip $\Im(z)<2$, even, real on the real line and such that $f(x)\ll 1/(1+x^2)$ as $x\to \infty$, we have
\begin{equation}\notag
\mathcal D_1^{\boldsymbol L_{\Delta,\chi}}(f)
=\int_{-\infty}^{+\infty}f(x)W_{SO^+}^{1}(x)dx+O\Big(\frac{1}{\log X}\Big)
\end{equation}
as $X\to\infty$, where
\begin{equation}\label{aggiuntaperkowalski}
W_{SO^+}^{1}(x):=1-\frac{\sin(2\pi x)}{2\pi x}. 
\end{equation}
%In particular we notice that $W_O^{(1)}(x)\sim \frac{2\pi^2}{3}x^2$, as $x\to0$. 
\end{prop}

\proof
The strategy of the proof is the same as in the unitary and symplectic cases, thus we will just sketch how the proof works, highlighting the differences with the other cases. For $a\asymp 1/\log X$ a real parameter, we consider the quantity
\begin{equation}\notag
\frac{1}{2\mathcal AX^*}\sum_{d\leq X}\sum_{\gamma_{\Delta,d}}f\Big(\frac{\log X}{\pi}\gamma_{\Delta,d}\Big)L_\Delta(\tfrac{1}{2}+\tfrac{a}{\log X},\chi_d)
\end{equation}
which can be written as ($\delta\asymp(\log X)^{-1}$)
\begin{equation}\begin{split}\notag
&\frac{\log (X^2)}{2\pi}\int_{-\infty}^{+\infty}f\Big(\frac{\log X}{\pi}y\Big)dy
\\&\quad\quad+2\frac{1}{2\pi}\int_{-\infty}^{+\infty}f\Big(\frac{\log X}{\pi}y\Big)\frac{1}{2\mathcal AX^*}\sum_{d\leq X}\frac{L'_\Delta}{L_\Delta}(\tfrac{1}{2}+\delta+iy)L_\Delta(\tfrac{1}{2}+\tfrac{a}{\log X},\chi_d)dy
\end{split}\end{equation}
with an error $O((\log X)^{-1})$, by using the Cauchy's theorem and the functional equation $\frac{L'_\Delta}{L_\Delta}(1-s,\chi_d)=\frac{Y'_\Delta}{Y_\Delta}(s,\chi_d)-\frac{L'_\Delta}{L_\Delta}(s,\chi_d)$ with $\frac{Y_\Delta'}{Y_\Delta}(s,\chi_d)=-\log d^2+O(1)$ (note that the square here is due to the conductor of $L_\Delta(s,\chi_d)$, which is $\frac{d^2}{4\pi^2}$). With the change of variable $\frac{\log X}{\pi}y=x$ the above equals
\begin{equation}\notag
\int_{-\infty}^{+\infty}f(x)\bigg(1+\frac{1}{2\mathcal AX^*\log X}\sum_{d\leq X}\frac{L'_\Delta}{L_\Delta}\Big(\frac{1}{2}+\delta+\frac{\pi ix}{\log X}\Big)L_\Delta(\tfrac{1}{2}+\tfrac{a}{\log X},\chi_d)\bigg)dx.
\end{equation}
Now we use the assumption of the ratio conjecture in the form of~\eqref{5.2.2.10} and~\eqref{5.2.2.11} to evaluate the sum over $d$, getting
\begin{equation}\notag
\int_{-\infty}^{+\infty}f(x)\bigg(1+\frac{1}{2\log X}h_X\Big(\frac{a}{\log X},\frac{\pi ix}{\log X}\Big)\bigg)dx+O\Big(\frac{1}{\log X}\Big)
\end{equation}
with
\begin{equation}\begin{split}\notag
h_X(\alpha,w):=
\frac{\alpha-w}{(\alpha+w)2w}&
+X^{-2\alpha}\frac{-\alpha-w}{(-\alpha+w)2w}
\\&+X^{-2w}\frac{-\alpha-w}{(\alpha-w)2w}
+X^{-2\alpha-2w}\frac{\alpha-w}{(-\alpha-w)2w}.
\end{split}\end{equation}
Letting $a\to0$, and since $h_X(0,w)=\frac{X^{-2w}-1}{w}$, we get
\begin{equation}\notag
\int_{-\infty}^{+\infty}f(x)\bigg(1+\frac{e^{-2\pi ix}-1}{2\pi ix}\bigg)dx+O\Big(\frac{1}{\log X}\Big).
\end{equation}
Putting all together, since $f$ is even, we finally have
\begin{equation}\notag
\mathcal D_1^{\boldsymbol L_{\Delta,\chi}}(f)=\int_{-\infty}^{+\infty}f(x)\bigg(1-\frac{\sin(2\pi x)}{2\pi x}\bigg)dx+O\Big(\frac{1}{\log X}\Big).
\end{equation}
\endproof

Similarly we compute the analogue of~\eqref{5.6.1}, tilting by the second power of $L_\Delta(\frac{1}{2},\chi_d)$, i.e.
\begin{equation}\label{5.6.10}
\mathcal D_2^{\boldsymbol L_{\Delta,\chi}}(f):=\frac{1}{\sum_{d\leq X}L_\Delta(\frac{1}{2},\chi_d)^2}\sum_{d\leq X}\sum_{\gamma_{\Delta,d}}f\Big(\frac{\log X}{\pi}\gamma_{\Delta,d}\Big)L_\Delta(\tfrac{1}{2},\chi_d)^2
\end{equation}
under the assumption of Conjecture~\ref{RCorthogonal}, in the case $K=3,Q=1$. This is achieved in the following proposition.

\begin{prop}\label{W1LD-O2}
Assume GRH and Conjecture~\ref{RCorthogonal} for $K=3, Q=1$. For any function $f$ holomorphic in the strip $\Im(z)<2$, even, real on the real line and such that $f(x)\ll 1/(1+x^2)$ as $x\to \infty$, we have
\begin{equation}\notag
\mathcal D_2^{\boldsymbol L_{\Delta,\chi}}(f)
=\int_{-\infty}^{+\infty}f(x)W_{SO^+}^{2}(x)dx+O\Big(\frac{1}{\log X}\Big)
\end{equation}
as $X\to\infty$, where
\begin{equation}\notag
W_{SO^+}^{2}(x):=1+\frac{\sin(2\pi x)}{\pi x}-\frac{2\sin^2(\pi x)}{(\pi x)^2}. 
\end{equation}
%In particular we notice that $W_O^{(2)}(x)\sim \frac{2\pi^4}{45}x^4$, as $x\to0$.
\end{prop}

\proof
Again we start with
\begin{equation}\notag
\frac{1}{4\mathcal BX^*\log X}\sum_{d\leq X}\sum_{\gamma_{\Delta,d}}f\Big(\frac{\log X}{\pi}\gamma_{\Delta,d}\Big)L_\Delta(\tfrac{1}{2}+\tfrac{a}{\log X},\chi_d)L_\Delta(\tfrac{1}{2}+\tfrac{b}{\log X},\chi_d)
\end{equation}
and with the usual machinery we write it as
\begin{equation}\notag
\int_{-\infty}^{+\infty}f(x)\bigg(1+\frac{1}{4\mathcal BX^*(\log X)^2}\mathcal I\Big(\frac{a}{\log X},\frac{b}{\log X}, \delta+\frac{\pi ix}{\log X}\Big)\bigg)dx+O\Big(\frac{1}{\log X}\Big)
\end{equation}
with
\begin{equation}\notag
\mathcal I(\alpha,\beta,w):=\sum_{d\leq X}\frac{L'_\Delta}{L_\Delta}(\tfrac{1}{2}+w,\chi_d)L_\Delta(\tfrac{1}{2}+\alpha,\chi_d)L_\Delta(\tfrac{1}{2}+\beta,\chi_d).
\end{equation}
Thanks to the assumption of the ratio conjecture (see Equations~\eqref{5.aux11} and~\eqref{5.aux20}) we are able to evaluate asymptotically the above sum, which is regular at $\alpha,\beta,\delta=0$. More specifically we have that
\begin{equation}\notag
\lim_{\substack{a\to0 \\ b\to0 \\ \delta\to0}}\;\mathcal I\Big(\frac{a}{\log X},\frac{b}{\log X}, \delta+\frac{y}{\log X}\Big)
=\mathcal BX^*(\log X)^2 h(y)+ O(\log X)
\end{equation}
with
\begin{equation}\notag
h(y):=\frac{-2ye^{-2y}-6y-4e^{-2y}+4}{y^2}.
\end{equation}
Then we get
\begin{equation}\begin{split}\notag
\mathcal D_2^{\boldsymbol L_{\Delta,\chi}}(f)&=\int_{-\infty}^{+\infty}f(x)\bigg(1+\frac{1}{4}h(\pi ix)\bigg)dx+O\Big(\frac{1}{\log X}\Big)
\\& =\int_{-\infty}^{+\infty}f(x)\bigg(1+\frac{\sin(2\pi x)}{2\pi x}-\frac{\sin^2(\pi x)}{(\pi x)^2}\bigg)dx+O\Big(\frac{1}{\log X}\Big)
\end{split}\end{equation}
since $f$ is even.
\endproof

We go on and define
\begin{equation}\notag
\mathcal D_3^{\boldsymbol L_{\Delta,\chi}}(f):=\frac{1}{\sum_{d\leq X}L_\Delta(\frac{1}{2},\chi_d)^3}\sum_{d\leq X}\sum_{\gamma_{\Delta,d}}f\Big(\frac{\log X}{\pi}\gamma_{\Delta,d}\Big)L_\Delta(\tfrac{1}{2},\chi_d)^3,
\end{equation}
analyzing the third-moment case.

\begin{prop}\label{W1LD-O3}
Assume GRH and Conjecture~\ref{RCorthogonal} for $K=4, Q=1$. For any function $f$ holomorphic in the strip $\Im(z)<2$, even, real on the real line and such that $f(x)\ll 1/(1+x^2)$ as $x\to \infty$, we have
\begin{equation}\notag
\mathcal D_3^{\boldsymbol L_{\Delta,\chi}}(f)
=\int_{-\infty}^{+\infty}f(x)W_{SO^+}^{3}(x)dx+O\Big(\frac{1}{\log X}\Big)
\end{equation}
as $X\to\infty$, where
\begin{equation}\notag
W_{SO^+}^{3}(x):=1-\frac{\sin(2\pi x)}{2\pi x}
-\frac{24(1-\sin^2(\pi x))}{(2\pi x)^2}
+\frac{48\sin(2\pi x)}{(2\pi x)^3}
-\frac{96\sin^2(\pi x)}{(2\pi x)^4}. 
\end{equation}
\end{prop}

\proof
We introduce the usual real parameters $\alpha,\beta,\nu$ of size $\asymp1/\log X$, we denote
$$ \boldsymbol L_{\Delta}^{\alpha,\beta,\nu}(\tfrac{1}{2},\chi_d):=L_{\Delta}(\tfrac{1}{2}+\alpha,\chi_d)L_{\Delta}(\tfrac{1}{2}+\beta,\chi_d)L_{\Delta}(\tfrac{1}{2}+\nu,\chi_d) $$
and we consider
$$ \frac{1}{\sum_{d\leq X}L_{\Delta}(\frac{1}{2},\chi_d)^3}\sum_{d\leq X}\sum_{\gamma_d}f\Big(\frac{\log X}{2\pi}\gamma_d\Big)\boldsymbol L_{\Delta}^{\alpha,\beta,\nu}(\tfrac{1}{2},\chi_d) .$$
With the usual strategy we get that the above equals
\begin{equation}\notag
\int_{-\infty}^{+\infty}f(x)\bigg(1+\frac{1}{\sum_{d\leq X}L(\frac{1}{2},\chi_d)^3}\sum_{d\leq X}\frac{L_{\Delta}'}{L_{\Delta}}(\tfrac{1}{2}+\delta+\tfrac{2\pi ix}{\log X},\chi_d)\boldsymbol L_{\Delta}^{\alpha,\beta,\nu}(\tfrac{1}{2},\chi_d)\bigg)dx
\end{equation}
up to an error $O(1/\log X)$, with $\delta\asymp 1/\log X$. We evaluate asymptotically the remaining sum over $d$ thanks to Conjecture~\ref{RCorthogonal} for $K=4,Q=1$), using Sage to carry out the computations. Doing so, letting $\alpha,\beta,\nu\to 0$, we obtain
\begin{equation}\notag
\mathcal D_3^{\boldsymbol L_{\Delta,\chi}}(f)=\int_{-\infty}^{+\infty}f(x)\bigg(1+\frac{1}{2} h(\pi i x)\bigg)dx+O\Big(\frac{1}{\log X}\Big)
\end{equation}
with
\begin{equation}\begin{split}\notag
h(y):=\frac{-5y^3+6y^2-6+e^{-2y}(y^3+6y^2+12y+6)}{y^4}.
\end{split}\end{equation}
The claim follows, since $f$ is even.
\endproof

Finally, in the following result we study the case $k=4$, given by
\begin{equation}\notag
\mathcal D_4^{\boldsymbol L_{\Delta,\chi}}(f):=\frac{1}{\sum_{d\leq X}L_\Delta(\frac{1}{2},\chi_d)^4}\sum_{d\leq X}\sum_{\gamma_{\Delta,d}}f\Big(\frac{\log X}{\pi}\gamma_{\Delta,d}\Big)L_\Delta(\tfrac{1}{2},\chi_d)^4.
\end{equation}

\begin{prop}\label{W1LD-O4}
Assume GRH and Conjecture~\ref{RCorthogonal} for $K=5, Q=1$. For any function $f$ holomorphic in the strip $\Im(z)<2$, even, real on the real line and such that $f(x)\ll 1/(1+x^2)$ as $x\to \infty$, we have
\begin{equation}\notag
\mathcal D_4^{\boldsymbol L_{\Delta,\chi}}(f)
=\int_{-\infty}^{+\infty}f(x)W_{SO^+}^{4}(x)dx+O\Big(\frac{1}{\log X}\Big)
\end{equation}
as $X\to\infty$, where
\begin{equation}\begin{split}\notag
W_{SO^+}^{4}(x):=
1&+\frac{\sin(2\pi x)}{2\pi x}-\frac{12\sin^2(\pi x)}{(\pi x)^2}-\frac{240\sin(2\pi x)}{(2\pi x)^3}
\\&-\frac{15(6-10\sin^2(\pi x))}{(\pi x)^4}+\frac{2880\sin(2\pi x)}{(2\pi x)^5}-\frac{90\sin^2(\pi x)}{(\pi x)^6}.
\end{split}\end{equation}
\end{prop}

\proof
As usual, if we set
$$ \boldsymbol L_{\Delta}^{\alpha,\beta,\nu,\eta}(\tfrac{1}{2},\chi_d):=L_{\Delta}(\tfrac{1}{2}+\alpha,\chi_d)L_{\Delta}(\tfrac{1}{2}+\beta,\chi_d)L_{\Delta}(\tfrac{1}{2}+\nu,\chi_d)L_{\Delta}(\tfrac{1}{2}+\eta,\chi_d),  $$
then we express $\mathcal D_4^{\boldsymbol L_{\Delta,\chi}}(f)$ as the limit for $\alpha,\beta,\nu,\eta\to0$ of 
$$ \int_{-\infty}^{+\infty}f(x)\bigg(1+\frac{1}{\sum_{d\leq X}L(\frac{1}{2},\chi_d)^4}\sum_{d\leq X}\frac{L_{\Delta}'}{L_{\Delta}}(\tfrac{1}{2}+\delta+\tfrac{2\pi ix}{\log X},\chi_d)\boldsymbol L_{\Delta}^{\alpha,\beta,\nu,\eta}(\tfrac{1}{2},\chi_d)\bigg)dx $$
up to an error $O(1/\log X)$, with $\delta\asymp 1/\log X$. 
The above can be evaluated asymptotically (again Sage is of help in carrying out the computation) and we get
$$ \mathcal D_4^{\boldsymbol L_{\Delta,\chi}}(f) =\int_{-\infty}^{+\infty}f(x)\bigg(1+ \frac{1}{2} h(\pi i x)\bigg)dx+O\Big(\frac{1}{\log X}\Big)$$
with
\begin{equation}\begin{split}\notag 
h(y):=&\frac{-7y^5+12y^4-30y^2+90}{y^6}
\\&-\frac{e^{-2y}(y^5+12y^4+60y^3+150y^2+180y+90)}{y^6}. 
\end{split}\end{equation}
Since $f$ is even the claim follows.
\endproof

Theorem~\ref{thm1LDorthogonal} follows by Propositions~\ref{W1LD-O1}$-$\ref{W1LD-O4}.

\subsection{Proof of Theorem \ref{analogoRMTortogonale}} Citing \cite[Theorem 4.3]{CFZratioconj} instead of assuming Conjecture \ref{RCorthogonal}, the same strategy gives an unconditional proof in random matrix theory.

%%%%%%%%%%%%%%%%%%%%%%%%%%%%%%%%%%%%%%%%%%%%%%%%%%%

\section{Proof of Theorem~\ref{conj1implicaconj2}}

%\begin{prop}\label{conj1implicaconj2}  Let us assume Conjecture \ref{congetturadecisiva1}. Then for any $k\in\mathbb N$ we have
%$$ W_{U\!Sp}^{k}(x)=\sum_{m=1}^\infty \beta_m x^{2m} $$
%with
%$$ \beta_m=(-1)^{m+1}\frac{(2\pi)^{2m}}{(2m+1)!}\bigg((-1)^k+\frac{k(k+1)}{m+1} \pFq{3}{2}{1-k,k+2,m+1}{m+2,2}{1}\bigg) $$ 
%where ${}_3F_{2}$ denotes the generalized hypergeometric function.\\ In particular, Conjecture \ref{congetturadecisiva2} follows.
%\end{prop}

By Fourier inversion, we have that
$$W_{U\!Sp}^{k}(x) = \int_{-\infty}^{+\infty}\widehat W_{U\!Sp}^{k}(y)e^{2\pi ixy}dy =\sum_{n=0}^\infty \bigg(\frac{(2\pi i)^n}{n!}\int_{-\infty}^{+\infty}\widehat W_{U\!Sp}^{k}(y)y^n dy\bigg) x^n.$$
Moreover, since $\widehat W_{U\!Sp}^{k}$ is even, then $\int_{-\infty}^{+\infty}\widehat W_{U\!Sp}^{k}(y)y^n dy=0$ if $n$ is odd. Hence, by definition of $\widehat W_{U\!Sp}^{k}$
$$W_{U\!Sp}^{k}(x) = \sum_{m=0}^\infty \beta_{m,k} x^{2m}$$
with 
\begin{equation}\begin{split}\notag
\beta_{m,k}= \frac{(2\pi i)^{2m}}{(2m)!}\int_{-\infty}^{+\infty}\bigg(\delta_0(y)&+\chi_{[-1,1]}(y)\bigg(-\frac{2k+1}{2}\\&-k(k+1)\sum_{j=1}^k (-1)^j c_{j,k}\frac{|y|^{2j-1}}{2j-1}\bigg)\bigg)y^{2m}dy,
\end{split}\end{equation}
where $\chi_{[-1,1]}$ denotes the indicator function of the interval $[-1,1]$.
By computing the integral, being $\int_{-1}^1 y^{2m}dy=\frac{2}{2m+1}$ and $\int_{-1}^{1}y^{2m}|y|^{2j-1}dy=\frac{1}{m+j}$, the above yields
\begin{equation}\notag
\beta_{m,k}:= \delta_0(m)-\frac{(2\pi i)^{2m}}{(2m)!}\bigg[\frac{2k+1}{2m+1}+k(k+1)\sum_{j=1}^k\frac{(-1)^j}{(2j-1)(j+m)}c_{j,k} \bigg].
\end{equation}
Since 
$$ c_{j,k}=\frac{1}{j}\binom{k-1}{j-1}\binom{k+j}{j-1}=\frac{j}{k(k+1)}\binom{k+j}{j}\binom{k}{j} $$
then we get
\begin{equation}\label{AppB.1} \beta_{m,k}=\delta_0(m)-\frac{(2\pi i)^{2m}}{(2m)!}\bigg[ \frac{2k+1}{2m+1}+S_{m}(1) \bigg]\end{equation}
where
$$S_{m,k}(x)
%=\sum_{j=1}^k\frac{(-1)^jj}{(2j-1)(j+m)} \binom{k+j}{j}\binom{k}{j}
:=\sum_{j=1}^\infty\frac{(-1)^jj}{(2j-1)(j+m)} \binom{k+j}{j}\binom{k}{j}x^j.$$

Now we write the factors in the above sum in terms of the Pochhammer symbol, defined as
$$ (a)_0:=1 \quad \text{ and }\quad (a)_n:=a(a+1)(a+2)\cdots(a+n-1) \quad\text{for }n\geq 1;$$
namely
\begin{equation}\begin{split}\notag
&\binom{k+j}{j}=\frac{(k+1)_j}{j!}
\\& (-1)^j\binom{k}{j}=\frac{(-k)_j}{(1)j}
\\& \frac{j}{(2j-1)(j+m)}=\frac{1}{2m+1}\bigg(\frac{1}{2j-1}+\frac{m}{j+m}\bigg)
\\&\hspace{3cm}=\frac{1}{2m+1}\bigg(-\frac{(-1/2)_j}{(1/2)_j}+\frac{(m)_j}{(m+1)_j}\bigg),
\end{split}\end{equation}
so that we have
\begin{equation}\label{AppB.2} S_{m,k}(x)=\frac{1}{2m+1}\sum_{j=1}^\infty\bigg(-\frac{(-1/2)_j}{(1/2)_j}+\frac{(m)_j}{(m+1)_j}\bigg)\frac{(-k)_j(k+1)_j}{(1)_j}\frac{x^j}{j!} .\end{equation}
Reparametrising the sum and using $(a)_{j+1}=a(a+1)_j$, this gives
$$ S_{m,k}(x)=S_{m,k}^1(x)+S_{m,k}^2(x),$$
where
$$S_{m,k}^1(x):=\frac{-k(k+1)x}{2m+1}\sum_{j=0}^\infty\frac{(1/2)_j(-k+1)_j(k+2)_j}{(3/2)_j(2)_j}\frac{x^j}{j!}\frac{1}{j+1}$$
and
$$S_{m,k}^2(x):=\frac{-mk(k+1)x}{(2m+1)(m+1)}\sum_{j=0}^\infty\frac{(m+1)_j(-k+1)_j(k+2)_j}{(m+2)_j(2)_j}\frac{x^j}{j!}\frac{1}{j+1}.$$
By writing
$$ \frac{1}{j+1}=2\bigg(1-\frac{2j+1}{2j+2}\bigg) $$
and
$$ \frac{(1/2)_j}{(3/2)_j}=\frac{1}{2j+1} $$
then we get
$$S_{m,k}^1(x)=-\frac{2xk(k+1)}{m+1} \pFq{3}{2}{1-k,k+2,\frac{1}{2}}{\frac{3}{2},2}{x}-\frac{\pFq{2}{1}{-k,k+1}{1}{x}-1}{2m+1},$$
where ${}_pF_q$ denotes the generalized hypergeometric function, defined as
$$\pFq{p}{q}{a_1,\dots,a_p}{b_1,\dots,b_q}{x}:=\sum_{n=0}^\infty\frac{(a_1)_n\cdots (a_p)_n}{(b_1)_n\cdots (b_q)_n}\frac{x^n}{n!}.$$
Similarly, since $$ -\frac{1}{j+1}=\frac{1}{m}-\frac{j+m+1}{(j+1)m} $$
then we have
$$ S_{m,k}^2(x)=\frac{xk(k+1)}{(2m+1)(m+1)} \pFq{3}{2}{1-k,k+2,m+1}{m+2,2}{x}+\frac{\pFq{2}{1}{-k,k+1}{1}{x}-1}{2m+1}. $$
Therefore, substituting in~\eqref{AppB.2} yields
\begin{equation}\begin{split}\label{AppB.3} S_{m,k}(x)=&-\frac{2xk(k+1)}{m+1} \pFq{3}{2}{1-k,k+2,\frac{1}{2}}{\frac{3}{2},2}{x}
\\&\hspace{1cm}+\frac{xk(k+1)}{(2m+1)(m+1)} \pFq{3}{2}{1-k,k+2,m+1}{m+2,2}{x}.\end{split}\end{equation}
Plugging~\eqref{AppB.3} in~\eqref{AppB.1}, we obtain
\begin{equation}\begin{split}\label{AppB.4}
\beta_{m,k}=\delta_0(m)-&\frac{(-1)^m(2\pi)^{2m}}{(2m+1)!}\bigg(2k+1-2k(k+1)\pFq{3}{2}{1-k,k+2,\frac{1}{2}}{\frac{3}{2},2}{1}
\\&\hspace{2.8cm}+\frac{k(k+1)}{(m+1)}\pFq{3}{2}{1-k,k+2,m+1}{m+2,2}{1}\bigg).
\end{split}\end{equation}

Now we need a few lemmas, in order to be able to compute the remaining hypergeometric functions.

\begin{lemma}\label{AppB.Lemma1}
For any $k\in\mathbb N$ we have
$$\pFq{3}{2}{1-k,k+2,\frac{1}{2}}{\frac{3}{2},2}{1}=\begin{cases} \frac{1}{k+1} & \text{if } k \text{ even}\\ \frac{1}{k} & \text{if } k \text{ odd.}\end{cases}$$
\end{lemma}
\proof
We recall the reduction formula for the generalized hypergeometric function (see e.g. \cite{GM-reduction}, Equation (17) in the case $n=1$), being 
\begin{equation}\begin{split}\notag
&\pFq{A+1}{B+1}{a_1,\dots,a_A,c+n}{b_1,\dots,b_B,c}{x}
\\&\hspace{3.5cm}=\sum_{j=0}^{n}\binom{n}{j}\frac{1}{(c)_j}\frac{\prod_{i=1}^A(a_i)_j}{\prod_{i=1}^B(b_i)_j}\pFq{A}{B}{a_1+j,\dots,a_A+j}{b_1+j,\dots,b_B+j}{x}
\end{split}\end{equation}
for any $A,B$ positive integers, $n\in\mathbb N$.
The left hand side can be then written as
\begin{equation}\begin{split}\label{AppB.Lemma1.1}
\pFq{3}{2}{1-k,\frac{1}{2},k+2}{\frac{3}{2},2}{1}&=\sum_{j=0}^k\binom{k}{j}\frac{1}{(2)_j}\frac{(1-k)_j(\frac{1}{2})_j}{(\frac{3}{2})_j}\pFq{2}{1}{1-k+j,\frac{1}{2}+j}{\frac{3}{2}+j}{1}
\\&=\sum_{j=0}^{k-1}\binom{k}{j}\frac{1}{(2)_j}\frac{(1-k)_j(\frac{1}{2})_j}{(\frac{3}{2})_j}\pFq{2}{1}{1-k+j,\frac{1}{2}+j}{\frac{3}{2}+j}{1}
\end{split}\end{equation}
as $(1-k)_k=0$. The remaining hypergeometric function can be computed by applying Gauss' summation theorem (see e.g. \cite{Koepf}, Equation (3.1)), i.e. the formula
$$\pFq{2}{1}{a,b}{c}{1}=\frac{\Gamma(c)\Gamma(c-a-b)}{\Gamma(c-a)\Gamma(c-b)}, \hspace{1cm}\Re(c)>\Re(a+b).$$
We recall that if $a=-n$, $n\in\mathbb N$, this is the Chu-Vandermonde identity (see again \cite{Koepf}, immediately below Equation (3.1))
$$ \pFq{2}{1}{-n,b}{c}{1}=\frac{(c-b)_n}{(c)_n}. $$
This yields
\begin{equation}
\label{AppB.Lemma1.2}
\pFq{2}{1}{1-k+j,\frac{1}{2}+j}{\frac{3}{2}+j}{1}%=\frac{(k-j-1)!}{(j+\frac{3}{2})_{k-j-1}}\mathcal P_{k-j-1}^{(j+\frac{1}{2},j-k)}(-1)
=\frac{(k-j-1)!(j+\frac{1}{2})!}{(k-\frac{1}{2})!}
\end{equation}
for $k>j$.
Plugging Equation~\eqref{AppB.Lemma1.2} into~\eqref{AppB.Lemma1.1}, we get
\begin{equation}\begin{split}\label{AppB.Lemma1.3}
\pFq{3}{2}{1-k,\frac{1}{2},k+2}{\frac{3}{2},2}{1}
&=\frac{(k-1)!}{2(k-\frac{1}{2})!}\sum_{j=0}^{k-1}\binom{k}{j}(-1)^j\frac{(j-\frac{1}{2})!}{(j+1)!}
\\&=\frac{(k-1)!}{2(k-\frac{1}{2})!}\sum_{j=0}^{k}\binom{k}{j}(-1)^j\frac{(j-\frac{1}{2})!}{(j+1)!}-\frac{(-1)^k}{2k(k+1)}.
\end{split}\end{equation}
Moreover, since $(\frac{1}{2})_j=\frac{1}{\sqrt{\pi}}(j-\frac{1}{2})!$, $(-1)^j\binom{k}{j}=(-k)_j/j!$ and $(2)_j=(j+1)!$, we have
$$ \sum_{j=0}^{k}\binom{k}{j}(-1)^j\frac{(j-\frac{1}{2})!}{(j+1)!}=\sqrt{\pi}\;\pFq{2}{1}{-k,\frac{1}{2}}{2}{1}=2\frac{(k+\frac{1}{2})!}{(k+1)!}$$
by applying the Chu-Vandermonde identity. % (see e.g. \cite{Koepf}, equation (3.1)), i.e. the formula
%$$\pFq{2}{1}{a,b}{c}{1}=\frac{\Gamma(c)\Gamma(c-a-b)}{\Gamma(c-a)\Gamma(c-b)}, \hspace{1cm}\Re(c)>\Re(a+b).$$ 
Putting this into~\eqref{AppB.Lemma1.3}, we finally get
$$\pFq{3}{2}{1-k,\frac{1}{2},k+2}{\frac{3}{2},2}{1}
=\frac{k+\frac{1}{2}}{k(k+1)}-\frac{(-1)^k}{2k(k+1)}$$
and the claim follows.
\endproof

By using Lemma~\ref{AppB.Lemma1}, Equation~\eqref{AppB.4} becomes
\begin{equation}\begin{split}\label{AppB.5}
\beta_{m,k}=\delta_0(m)-&\frac{(-1)^m(2\pi)^{2m}}{(2m+1)!}\bigg((-1)^k+
\frac{k(k+1)}{(m+1)}\pFq{3}{2}{1-k,k+2,m+1}{m+2,2}{1}\bigg).
\end{split}\end{equation}

The coefficient $\beta_{0,k}$ can be then computed thanks to the following lemma.
\begin{lemma}\label{AppB.Lemma2}
For any $k\in\mathbb N$ we have
$$\pFq{3}{2}{1-k,k+2,1}{2,2}{1}=\begin{cases} \; 0 & \text{if } k \text{ even}\\  \frac{2}{k(k+1)} & \text{if } k \text{ odd.}\end{cases}$$
\end{lemma}
\proof
By definition we have
\begin{equation}\begin{split}\notag
\pFq{3}{2}{1-k,k+2,1}{2,2}{1}
&=\sum_{j=0}^\infty \frac{(1-k)_j(k+2)_j(1)_j}{(2)_j(2)_j}\frac{1}{j!}
\\&=-\frac{1}{k(k+1)}\sum_{j=0}^\infty \frac{(-k)_{j+1}(k+1)_{j+1}}{(1)_{j+1}}\frac{1}{j!} 
\end{split}\end{equation}
since $(1)_j/(2)_j=1/(j+1)$, $(1-k)_j=(-k)_{j+1}/(-k)$ and $(2)_j=(1)_{j+1}$.
Reparametrising the series with $l=j+1$, the above yields
$$ -\frac{1}{k(k+1)}\bigg(\sum_{l=0}^\infty \frac{(-k)_l(k+1)_l}{(1)_l}\frac{1}{l!}-1\bigg)
=\frac{1-\pFq{2}{1}{-k,k+1}{1}{1}}{k(k+1)}.$$
The claim is then proven, by noticing that $\pFq{2}{1}{-k,k+1}{1}{1}=(-1)^k$ thanks to the Chu-Vandermonde identity. 
\endproof

This implies that $\beta_{0,k}=0$ for any $k\in\mathbb N$, proving the first part of Theorem~\ref{conj1implicaconj2}. To complete our proof, we need to show that 
\begin{equation}\label{AppB.6} \beta_{k+1,k}=\frac{2\pi^{2(k+1)}}{(2k+1)!!(2k+1)!!} \end{equation} 
is the first nonzero coefficient. As a first step, the following lemma shows that $\beta_{i,k}=0$ for all $1\leq i\leq k$.

\begin{lemma}\label{AppB.Lemma3} 
For any $k\in\mathbb N$ and for any $1\leq m \leq k$ we have 
$$\pFq{3}{2}{1-k,k+2,m+1}{m+2,2}{1}=\frac{(m+1)(-1)^{k+1}}{k(k+1)}.$$
\end{lemma}
\proof
We begin by applying the reduction formula, which yields
\\ $\pFq{3}{2}{1-k,k+2,m+1}{m+2,2}{1}=$ 
\begin{equation}\label{AppB.L3.1}
\sum_{j=0}^{m-1}\binom{m-1}{j}\frac{1}{(2)_j}\frac{(1-k)_j(k+2)_j}{(m+2)_j}\pFq{2}{1}{1-k+j,k+2+j}{m+2+j}{1}.
\end{equation}
Moreover, the Chu-Vandermonde identity gives
\begin{equation}\notag
\pFq{2}{1}{1-k+j,k+2+j}{m+2+j}{1}% =\frac{(k-j-1)!}{(m+j+2)_{k-j-1}}\mathcal P_{k-j-1}^{(m+j+1,-m+j+1)}(-1) 
%\\&\hspace{4cm} =\frac{(k-j-1)!(m+j+1)!}{(m+k)!}(-1)^{k-j-1}\binom{k-m}{k-j-1}.
= \frac{(m-k)_{k-j-1}}{(m+j+2)_{k-j-1}}
\end{equation}
Since $(m-k)_{k-j-1}=0$ for all $j<m-1$, only the term $j=m-1$ survives in the sum in Equation~\eqref{AppB.L3.1}. Hence we get
\begin{equation}\begin{split}\label{AppB.L3.2}
&\pFq{3}{2}{1-k,k+2,m+1}{m+2,2}{1}
%=\sum_{j=0}^{m-1}\binom{m-1}{j}\frac{1}{(2)_j}\frac{(1-k)_j(k+2)_j}{(m+2)_j}\pFq{2}{1}{1-k+j,k+2+j}{m+2+j}{1}.
%=\binom{m-1}{m-1}\frac{1}{(2)_{m-1}}\frac{(1-k)_{m-1}(k+2)_{m-1}}{(m+2)_{m-1}}\frac{(k-m+1-1)!(m+m-1+1)!}{(m+k)!}(-1)^{k-m+1-1}\binom{k-m}{k-m+1-1}.
=\frac{(1-k)_{m-1}}{(2)_{m-1}}\frac{(k+2)_{m-1}}{(m+2)_{m-1}}\frac{(m-k)_{k-m}}{(2m+1)_{k-m}}
\\& \hspace{2cm} =\frac{(-1)^{m-1}(k-1)!}{(k-m)!m!}\frac{(k+m)!(m+1)!}{(k+1)!(2m)!}\frac{(-1)^{k-m}(k-m)!(2m)!}{(m+k)!}
\\& \hspace{2cm} =\frac{(-1)^{k-1}(k-1)!}{m!}\frac{(m+1)!}{(k+1)!}=\frac{(-1)^{k-1}(m+1)}{k(k+1)},
\end{split}\end{equation}
where in the first line we applied the equalities $(m+2)_{m-1}=(2m)!/(m+1)!$, $(k+2)_{m-1}=(k+m)!/(k+1)!$ and $(1-k)_{m-1}=(-1)^{m-1}(k-1)!/(k-m)!$. Similarly also $(m-k)_{k-m}=(-1)^{k-m}(k-m)!$ and $(2m+1)_{k-m}=(k+m)!/(2m)!$.
\endproof

Finally, with the following lemma, we can also compute $\beta_{k+1,k}$.

\begin{lemma}\label{AppB.Lemma4}
For any $k\in\mathbb N$ we have
$$\pFq{3}{2}{1-k,k+2,k+2}{k+3,2}{1}=\frac{2(-1)^{k+1}(k-1)!(k+2)!}{(2k+2)!}\bigg(\binom{2k+1}{k+1}-1\bigg).$$
\end{lemma}
\proof
The idea of the proof is similar the the one of Lemma~\ref{AppB.Lemma3}. First we apply the reduction formula in order to write $\pFq{3}{2}{1-k,k+2,k+2}{k+3,2}{1}$ as a finite sum of terms involving ${}_2F_{1}$, namely
\begin{equation}\begin{split}\label{AppB.L4.1}
&\pFq{3}{2}{1-k,k+2,k+2}{k+3,2}{1}
\\&=\sum_{j=0}^{k-1}\binom{k}{j}\frac{1}{(2)_j}\frac{(1-k)_j(k+2)_j}{(k+3)_j}\pFq{2}{1}{1-k+j,k+2+j}{k+3+j}{1}.
\end{split}\end{equation}
Note that the term $j=k$ vanishes, as $(1-k)_k=0$. Now we use Gauss's summation theorem and compute the remaining hypergeometric function, i.e.
$$ \pFq{2}{1}{1-k+j,k+2+j}{k+3+j}{1} 
%= \frac{(k+j+2)!}{(2k+1)!}(-1)^{k-j-1}\binom{-1}{k-j-1}
= \frac{(k+j+2)!}{(2k+1)!}. $$
Plugging this into Equation~\eqref{AppB.L4.1}, since $(-1)^j\binom{k}{j}=(-k)_j/j!$ and $(k+2)_j=(k+j+1)!/(k+1)!$, we have
\begin{equation}\begin{split}\label{AppB.L4.2}
&\pFq{3}{2}{1-k,k+2,k+2}{k+3,2}{1}=\frac{(k-1)!(k+2)!}{(2k+1)!}\sum_{j=0}^{k-1}\frac{(-k)_j(k+2)_j}{(2)_j}\frac{1}{j!}
\\&\hspace{1.5cm}=\frac{(k-1)!(k+2)!}{(2k+1)!}\bigg(\pFq{2}{1}{-k,k+2}{2}{1}-\frac{(-k)_k(k+2)_k}{(2)_k}\frac{1}{k!} \bigg).
\end{split}\end{equation}
Therefore, since $\pFq{2}{1}{-k,k+2}{2}{1}=\frac{(-1)^k}{k+1}$, we get
\begin{equation}\begin{split}\notag
\pFq{3}{2}{1-k,k+2,k+2}{k+3,2}{1}
&= \frac{(k-1)!(k+2)!}{(2k+1)!}\bigg(\frac{(-1)^k}{k+1}-\frac{(-1)^k(2k+1)!}{((k+1)!)^2} \bigg) 
\\& = \frac{2(k-1)!(k+2)!(-1)^k}{(2k+2)!}\bigg(1-\frac{(2k+1)!}{(k+1)!k!} \bigg)
\end{split}\end{equation}
and the claim follows.
\endproof

To conclude the proof of Theorem~\ref{conj1implicaconj2}, we just combine Equation~\eqref{AppB.5} with Lemma~\ref{AppB.Lemma4}, getting
\begin{equation}\begin{split}\notag
\beta_{k+1,k}&=\frac{(2\pi)^{2k+2}}{(2k+3)!}\bigg(1-\frac{k(k+1)}{k+2}\frac{2(k-1)!(k+2)!}{(2k+2)!}\bigg(\binom{2k+1}{k+1}-1\bigg)\bigg)
\\& = \frac{(2\pi)^{2k+2}}{(2k+3)!}\bigg(1-\frac{2[(k+1)!]^2}{(2k+2)!}\bigg(\frac{(2k+1)!}{(k+1)!k!}-1\bigg)\bigg)
\\& = \frac{(2\pi)^{2k+2}}{(2k+3)!}\bigg(1-\frac{2(k+1)}{2k+2}+\frac{2[(k+1)!]^2}{(2k+2)!}\bigg)
=\frac{(2\pi)^{2k+2}}{(2k+3)!}\frac{(k+1)!k!}{(2k+1)!}.
\end{split}\end{equation}

Equation~\eqref{AppB.6} follows by the identities $(2k+1)!=2^kk!(2k+1)!!$ and $(2k+3)!=2^{k+1}(k+1)!(2k+3)!!$.\\

\textbf{Acknowledgments}. I would like to thank Sandro Bettin for many interesting and helpful discussions, for all his help while I was working on this paper and also for several corrections and improvements that he suggested. I am also grateful to Joseph Najnudel for inspiring this project. I also wish to thank the referees for a very careful reading of the paper and for indicating several inaccuracies and mistakes. The author was supported by Czech Science Foundation, grant 21-00420M.

%%%%%%%%%%%%%%%%%%%%%%%%%%%%%%%%%%%%%%%%%%%%%%%%%%%

{\small

\end{document}